\documentclass[10pt,a4paper]{article}

\pdfoutput=1

\usepackage{resizegather}		
\usepackage[section]{placeins}
\usepackage{amsmath,amscd}
\usepackage{amsthm}
\usepackage{amsrefs}
\usepackage{amssymb}
\usepackage{amsfonts}
\usepackage{bm}
\usepackage{mathtools}
\usepackage{mathdots}
\usepackage[T1]{fontenc}
\usepackage{framed, color}
\usepackage{dsfont}
\usepackage{graphicx}
\usepackage{hyperref}
\usepackage{latexsym}
\usepackage{mathrsfs}
\usepackage{subfig}
\usepackage[section]{placeins}
\usepackage[customcolors]{hf-tikz}
\usepackage{tikz-cd}
\usepackage{verbatim}
\usepackage[a4paper, left=3.5cm, 
		     right=3.5cm, top=2.97cm, 
		    ]{geometry}
		    
\tikzset{  
	style gray/.style={
    set fill color=lightgray,
    set border color=white, opacity=0.4,
  },
  style gray_hor/.style={
    set fill color=lightgray,
    set border color=white, opacity=0.4,
  },
  hor/.style={
    above left offset={-0.4,0.4},
    below right offset={0.28,-0.2},
    #1
  },
  ver/.style={
    above left offset={-0.35,0.4},
    below right offset={0.35,-0.20},
    #1
  }
}
		    
\theoremstyle{plain}
\newtheorem{theorem}{Theorem}[section]
\newtheorem{proposition}[theorem]{Proposition}
\newtheorem{lemma}[theorem]{Lemma}

\newtheorem{remark}[theorem]{Remark}

\numberwithin{equation}{section}
\numberwithin{figure}{section}
\numberwithin{table}{section}

\theoremstyle{definition}
\newtheorem{definition}[theorem]{Definition}

\newtheorem{application}{Application}

\newcommand{\RR}{\mathbb R}
\newcommand{\ZZ}{\mathbb Z}
\newcommand{\CC}{\mathbb C}
\newcommand{\NN}{\mathbb N}

\newcommand{\cp}{\bm{\ast}}

\title{Validated computations for connecting orbits\\ in polynomial vector fields}
\author{Jan Bouwe van den Berg\thanks{
Department of Mathematics, 
VU Amsterdam, 
1081 HV Amsterdam, 
The Netherlands, {\tt janbouwe@few.vu.nl};
partially supported by NWO-VICI grant 639033109.
} 
\and 
Ray Sheombarsing
}
\date{\today}

\begin{document}

\maketitle

\begin{abstract}
In this paper we present a computer-assisted procedure for proving the existence of transverse heteroclinic orbits
	connecting hyperbolic equilibria of polynomial vector fields. The idea is to compute high-order Taylor approximations
	of local charts on the (un)stable manifolds by using the Parameterization Method and to use Chebyshev series to
	parameterize the orbit in between, which solves a boundary value problem. The existence of a heteroclinic orbit can then be established by setting up an
	appropriate fixed-point problem amenable to computer-assisted analysis. The fixed point problem simultaneously solves for the local (un)stable manifolds and the orbit which connects these. We obtain explicit rigorous control on the distance between the numerical approximation and the heteroclinic orbit. Transversality of the stable and unstable manifolds is also proven.
\end{abstract}

\section{Introduction}
\label{sec:introduction}
Connecting orbits play a central role in the study of dynamical systems. They provide a detailed
picture of how a dynamical system can evolve from one ``state'' (e.g. an equilibrium, a periodic orbit
or another type of recurrent set) into another. 
Furthermore, their existence can often be used to establish more
complicated dynamical phenomena through forcing theorems. 
However, proving the existence of a connecting orbit for a given nonlinear ODE is in general a difficult 
(if not impossible) task to accomplish by hand. 
For this reason, one often resorts to numerical methods. 
While numerical methods can provide valuable insight into quantitive properties of a connecting orbit, 
which would otherwise be out of reach 
with merely a pen and paper analysis, the results are usually non-rigorous and cannot be used in mathematical arguments. In particular, 
a standard numerical method does not yield a proof for the existence of a connecting orbit. 

In this paper we present a general computer-assisted method for proving the existence of 
transverse connecting orbits between hyperbolic equilibria for nonlinear ODEs.
The method is based on solving the \emph{finite time} boundary value problem
\begin{align}
	\label{eq:main_BVP}
	\begin{cases}
		\dfrac{du}{dt} = g \left( u \right), & t \in [0,L], \\[2ex]
		u(0) \in W^{u}_{\text{loc}} \left( p_{0} \right) , \\[2ex]
		u(L) \in W^{s}_{\text{loc}} \left( q_{0} \right),
	\end{cases}
\end{align}
where $g: \RR^{n} \rightarrow \RR^{n}$ is a \emph{general} polynomial vector field,  $p_{0}, q_{0} \in \RR^{n}$ are hyperbolic equilibria and 
$L>0$ is the time needed to travel between the local (un)stable manifolds.
We assume that $\dim \left( W^{u} \left( p_{0} \right) \right) + \dim \left( W^{s} \left( q_{0} \right) \right)  = n +1$, which is a necessary condition for a transverse connecting orbit to exist.
The idea is to solve \eqref{eq:main_BVP} by computing Taylor expansions for charts on the local (un)stable manifolds via
the parameterization method \cite{ManifoldTaylor, MR3437754}, which are used to  supplant the boundary conditions in \eqref{eq:main_BVP} with explicit equations, 
and to use Chebyshev series and domain decomposition techniques \cite{domaindecomposition} to parameterize the orbit in between. 
We remark that the assumption that $g$ is polynomial 
is not as restrictive as it initially might seem, since many nonlinearities which consist of elementary functions 
can be brought into polynomial form by using automatic differentiation techniques, see \cite{MR3545977} for instance. 

Before we proceed with a more detailed description of our method, a few remarks concerning the development of numerical methods
for connecting orbits are in order. Many numerical methods (both rigorous and
non-rigorous) are based on approximating (un)stable manifolds and solving finite time boundary value problems. We mention the numerical (non-rigorous) methods implemented in the
continuation packages \textsc{Matcont} \cite{MR2000880} and \textsc{AUTO} \cite{AUTO} in particular. Furthermore, many validated numerical methods have been developed over the last decade. 
It is beyond the scope of this text to give an overview.
Nevertheless, we mention the functional analytic methods developed in 
\cite{Maxime, domaindecomposition, MR3281845, MR3022075, suspensionbridge, MR3353132, MR3207723, MR3148084, MR2821596, Kepley, JayAMS}, which are based on 
solving fixed point problems, the topological methods developed in \cite{MR3463691, MR2173545, MR2505658, MR1961956},
which are based on covering relations, cone conditions \cite{MR2494688, MR2060531,MR2060532 } and rigorous integration 
of the flow via Lohner-type algorithms \cite{MR1930946}, and the methods in \cite{MR2302059, MR2339601} based on shadowing techniques.

The first step in the development of our validated numerical method is to recast \eqref{eq:main_BVP} into  an equivalent zero finding problem $F(x)=0$.
The unknowns in this problem are the Taylor coefficients of the
parameterizations of the local (un)stable manifolds, which include the equilibria and the associated eigendata, the coordinates 
of the endpoints $u(0)$ and $u(L)$ on the associated charts, and the Chebyshev coefficients of the orbit. Next, we use the computer 
to determine an approximate zero of $F$. The numerical computations are then combined with 
analysis on paper to construct a Newton-like map $T$ whose fixed points correspond to zeros of $F$. 
Finally, we use pen and paper estimates to derive a \emph{finite} number of inequalities, which can be used 
to determine a neighborhood around the approximate solution on which~$T$ is a contraction. An essential property of these inequalities
is that they can be rigorously verified with the aid of a computer. 
This approach is in the literature often referred to as a parameterized Newton-Kantorovich method or the radii-polynomial approach
(see \cite{MR1639986, MR2338393}). 
 
Our construction builds on the work in \cite{suspensionbridge, MR3353132, MR3148084}. In those papers heteroclinic orbits for particular vector fields were studied using the parameterization method, Chebyshev series and the radii-polynomial approach to validate solutions of the boundary value problem~\eqref{eq:main_BVP}. 
Whereas in \cite{suspensionbridge, MR3353132, MR3148084} the charts on the local (un)stable manifolds were validated separately from the proof of the connecting orbit between them, the main contribution of the current paper is that we solve all ingredients of the problem \emph{simultaneously}, and in great generality. While this means we have to introduce a fair amount of notation, the advantage is that our framework is very flexibly. Indeed, the mathematical analysis and the
code are developed in full detail for general polynomial vector fields. 
Moreover, setting up this framework opens the door for various extensions. In particular, the formulation as a zero finding problem in an appropriately chosen (product) Banach space provides a foundation to study continuation and bifurcation problems, see also Section~\ref{s:extensions}.
Furthermore, the presented method deals in a 
unified and systematic manner with the cases of real and complex eigenvalues associated to the parameterizations of the (un)stable manifolds (see also the examples in Section~\ref{s:examples}). We note that in part of the paper we impose a technical ``non-resonance'' condition on the eigenvalues of the linearized problems at the equilibria, see Section~\ref{s:charts}. This condition is by no means fundamental, but it reduces the (already heavy) notational burden. As explained in Section~\ref{s:extensions}, the general case follows by combining the current work with~\cite{ManifoldTaylor}. 

Finally, we note that it is a feature of the Newton-Kantorovich fixed point method that any heteroclinic orbit that we find using our methodology is the \emph{transverse} intersection of the stable and unstable manifolds of the equilibria it connects, see Proposition~\ref{prop:transversality}. Such transversality information is often essential when one aims to use the connecting orbits as ingredients for forcing theorems.

In this paper we at times refer to previous work for certain proofs, bit we give sufficient details of the constructions to understand the algorithm. 
Hence the duo of paper and code is self-contained.
The fully documented Matlab code is available at~\cite{heteroclinicscode}.

\subsection{Applications}
\label{s:examples}

To illustrate the efficacy of the method,
we have chosen two travelling wave problems originating from partial differential equations 
Here we present some example results. Additional rigorously validated orbits and more details can be found in Section~\ref{sec:applications}.

\begin{application}[Traveling fronts in the Lotka-Volterra equations]
	The Lotka-Volterra equations are a system of reaction-diffusion 
	equations given by 
	\begin{align}
		\label{eq:LK}
		\begin{cases}
			\dfrac{ \partial v }{ \partial t} = D \dfrac{ \partial^{2} v}{ \partial x^{2} } + v \left( 1 - v - w \right), \\[2ex]
			\dfrac{ \partial w }{ \partial t} = \dfrac{ \partial^{2} w}{ \partial x^{2} } + a w \left( v - b \right),
		\end{cases}
	\end{align}
	where $D >0,a>0, b\in (0,1)$ and $\left(t,x \right) \in \RR^{2}$. This system has three homogeneous equilibrium states:
	$(v,w)=(0,0)$, $(v,w) = (1,0)$ and $(v,w) = (b,1-b)$. 	
	We have used our method to prove the existence of 
	solutions of \eqref{eq:LK} of the form $v(t,x) = \zeta_{1}(x-\kappa t)$ and $w(t,x) = \zeta_{2}(x-\kappa t)$, where
	$\zeta_{1}, \zeta_{2} : \RR \rightarrow \RR$ and $\kappa <0$, which satisfy  
	\begin{align*}
		\lim_{\tau \rightarrow - \infty} \left( \zeta_{1} \left( \tau \right), \zeta_{2} \left( \tau \right) \right) = (1,0), \quad
		\lim_{\tau \rightarrow \infty} \left( \zeta_{1} \left( \tau \right), \zeta_{2} \left( \tau \right) \right) = (b,1-b).
	\end{align*}
	Such solutions are often referred to as \emph{traveling fronts} with wave speed $\kappa$. 
	
	Substitution of the traveling wave Ansatz $\left( \zeta_{1}, \zeta_{2} \right)$ into \eqref{eq:LK}  shows that connecting orbits from 
	$(b,0,1-b,0)$ to $(1,0,0,0)$ for the four dimensional system of ODEs 
	\begin{align}
		\dfrac{du}{dt} = 
			\begin{bmatrix}
				- u_{2} \\[1ex]
				D^{-1} \left( \kappa u_{2} + u_{1} \left( 1 - u_{1} - u_{3} \right) \right) \\[1ex]
				- u_{4} \\[1ex]
				 \left( \kappa u_{4} + a u_{3} \left( u_{1} - b \right) \right)
			\end{bmatrix}
		\label{eq:LK_ODE}
	\end{align}
	correspond to traveling wave profiles 
	$\left( \zeta_{1}(t), \zeta_{2}(t) \right) = \left( u_{1}(-t), u_{3}(-t) \right)$ and vice versa. 
	We have successfully validated connecting orbits in \eqref{eq:LK_ODE} for various values of $\kappa \in \left[-1,-0.5938 \right]$.
	For these parameter values, the equilibria $(b,0,1-b,0)$ and $(1,0,0,0)$ have a two dimensional unstable and three dimensional stable
	manifold, respectively. In particular, the stable eigenvalues of the linearization at $\left(1,0,0,0 \right)$ 
	consist of one complex conjugate pair of eigenvalues and one real eigenvalue. 
	We have depicted a validated traveling wave profile and the corresponding connecting orbit for a particular wave speed
	in Figures \ref{fig:LK_u1_u3} and \ref{fig:LK_phase_space}, respectively. 
	The reader is referred to Section \ref{sec:LK} for the details. 
	\begin{figure}[tbh]
		\begin{center}
			\includegraphics[scale = 0.5]{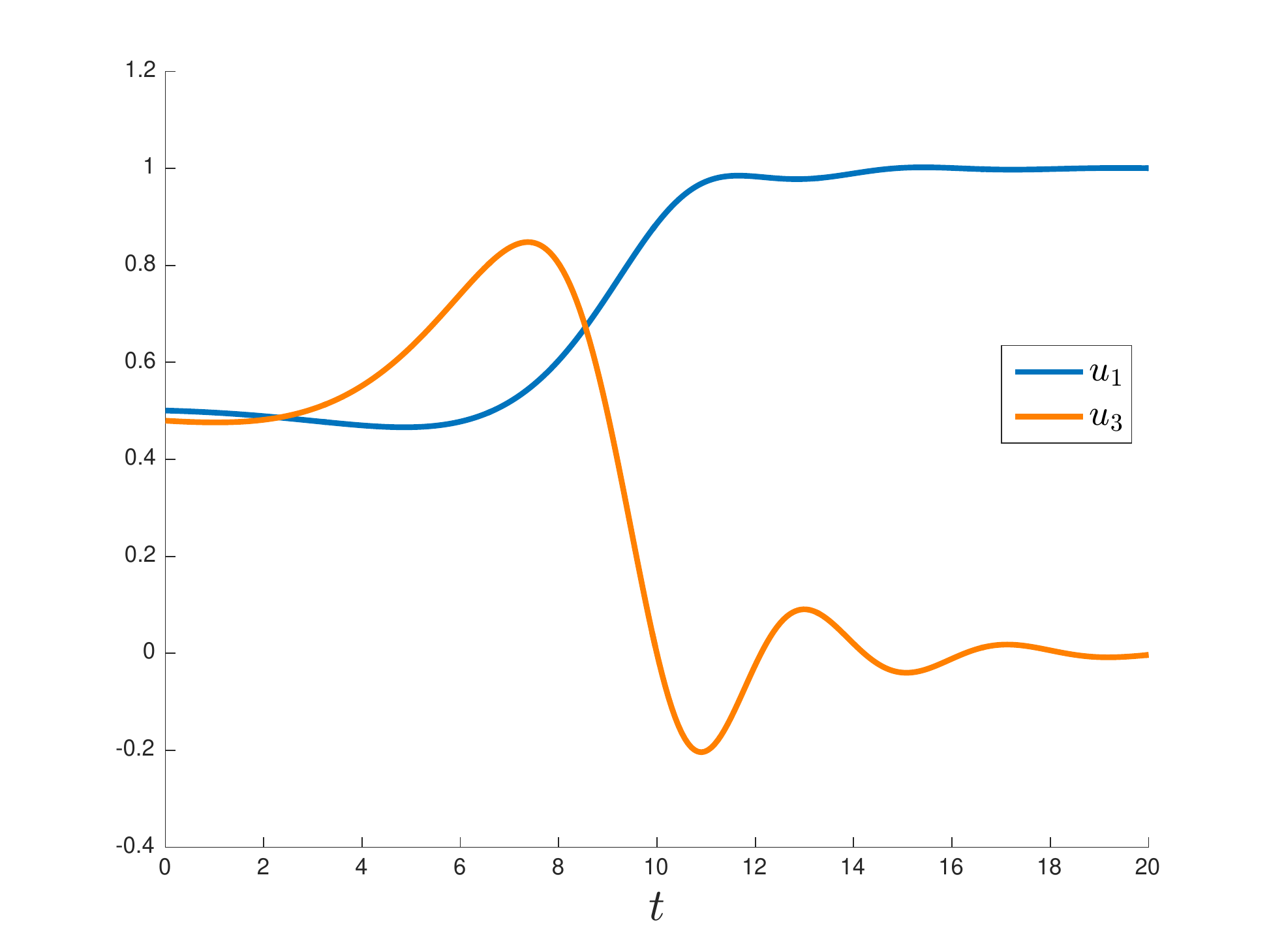}
		\end{center}
		\caption{\label{fig:LK_u1_u3} Validated traveling wave profiles of \eqref{eq:LK} for
		$a=5$, $D=3$, $b= \frac{1}{2}$ and $\kappa = -0.7767$. The depicted parts of the traveling wave profiles
		correspond to the first and third component of the connecting orbit between the unstable and stable manifold of 
		$\left(b,0,1-b,0 \right)$ and $\left(1,0,0,0\right)$, respectively. 
		The time of flight of the connection between the (un)stable manifolds was $L=20$. 
	}
	\end{figure}	
	\begin{figure}[tbh]
		\begin{center}
			\includegraphics[scale = 0.1]{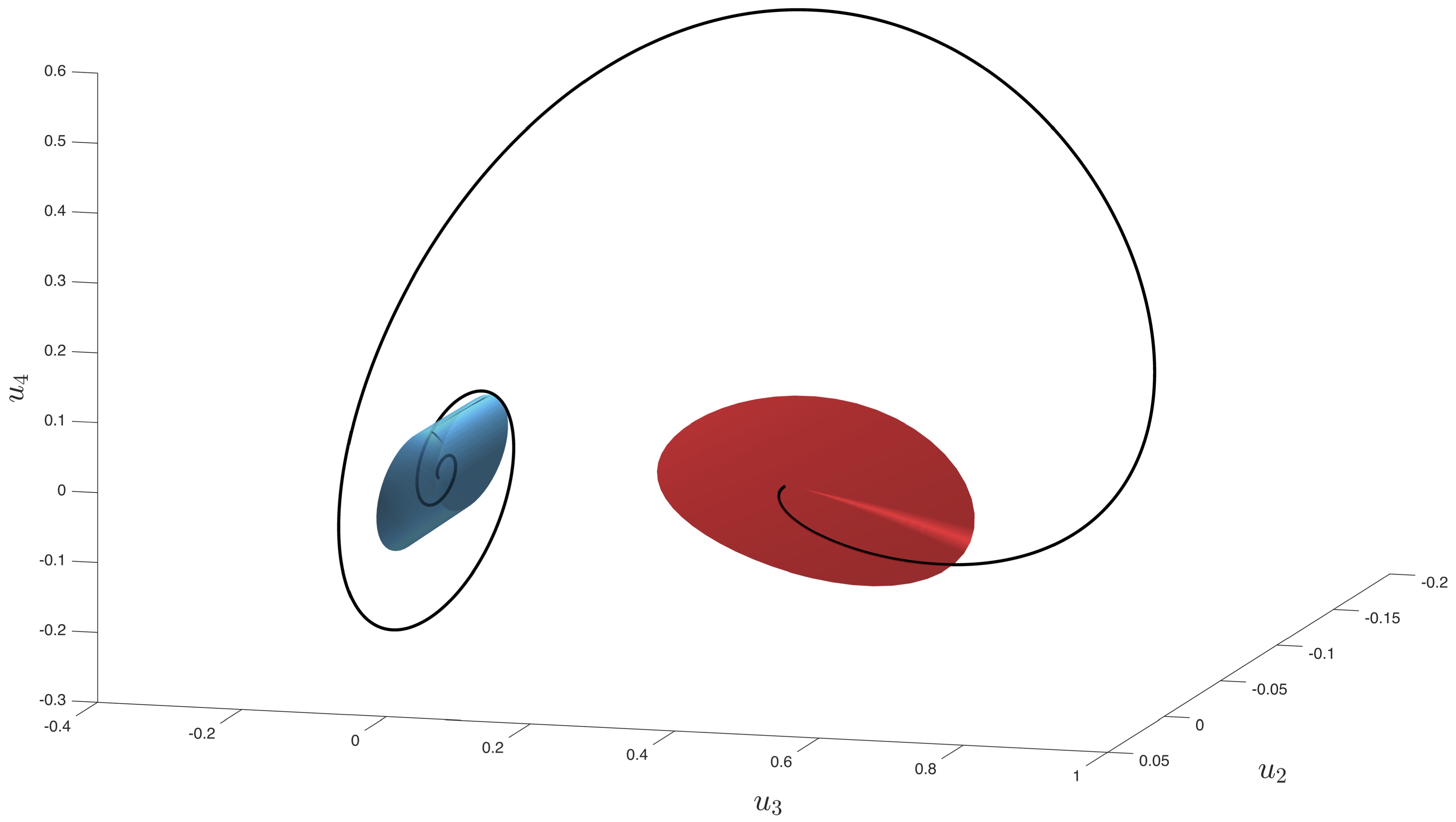}
		\end{center}
		\caption{\label{fig:LK_phase_space} A three dimensional projection of a validated connecting orbit of \eqref{eq:LK_ODE} for
		$a=5$, $D=3$, $b= \frac{1}{2}$ and $\kappa = -0.7767$. The geometric objects colored in red and blue correspond to 
		parameterizations of the local unstable and local stable manifold of $\left(b,0,1-b,0 \right)$ and $\left(1,0,0,0\right)$, respectively, 
		which were computed (and validated) by using the parameterization method. The curve in black corresponds to the 
		piece of the connecting orbit which was validated by using Chebyshev series. The time of flight of the connection 
		between the (un)stable manifolds was $L=20$. We remark that the integration time was not optimized, i.e., it is 
		possible to decrease $L$ and ``arrive'' in the parameterized local (un)stable manifolds earlier.  
	}
	\end{figure}
\end{application}

\begin{application}[Traveling fronts in a fourth order parabolic PDE]
		We have proven the existence of traveling fronts $v(t,x) = \zeta \left( x - \kappa t \right)$ 
		for the following fourth order parabolic PDE: 
		\begin{align}
				\label{eq:Fourth_order_PDE}
				\frac{ \partial v }{ \partial t} = - \gamma \frac{ \partial^{4} v }{ \partial x^{4} } + 
				\frac{ \partial^{2} v }{ \partial x^{2} } + \bigl( v-a \bigr) \bigl( 1- v^{2} \bigr)
		\end{align}
		where $-1 < a \leq 0$ and $\gamma >0$.
Such travelling waves have been studied in
		\cite{MR1629027} used geometric singular perturbation theory for small $\gamma$,
		and in~\cite{BHV} by Conley index techniques.
		 In this paper we focus
		on traveling fronts between the homogeneous states $v \equiv -1$ and $v \equiv a$ for 
		(relatively) large $\gamma$. 
		
		We have successfully established the existence of connecting orbits from $(-1,0,0,0)$ to $(a,0,0,0)$ 
		for the four dimensional system of ODEs
		\begin{align}
				\label{eq:Fourth_order_ODE}
				\frac{du}{dt} = 
				- \begin{bmatrix}
						\gamma u_{2} \\
						\gamma u_{3} \\
						\gamma u_{4} \\
						\kappa u_{2} + u_{3} + \left( u_{1} - a \right) \left( 1 - u_{1}^{2} \right)
				\end{bmatrix},
		\end{align}
		which correspond to traveling wave profiles $\zeta \left( t \right) = u_{1} \left( -\frac{t}{\gamma} \right)$, 
		for fixed values of $a$ and the wave speed $\kappa$, and various values of $\gamma \in [0.4557, 10.50]$. We 
		rescaled time with a factor $\gamma$ so that the system 
		in \eqref{eq:Fourth_order_ODE} is well-defined at $\gamma = 0$. 
		We have depicted a validated traveling wave profile and the corresponding connecting orbit
		in Figures \ref{fig:FK_u1} and \ref{fig:FK_phase_space}, respectively. The reader is referred to Section \ref{sec:Fourth_order}
		for the details. 
		\begin{figure}[tbh]
			\begin{center}
				\includegraphics[scale = 0.5]{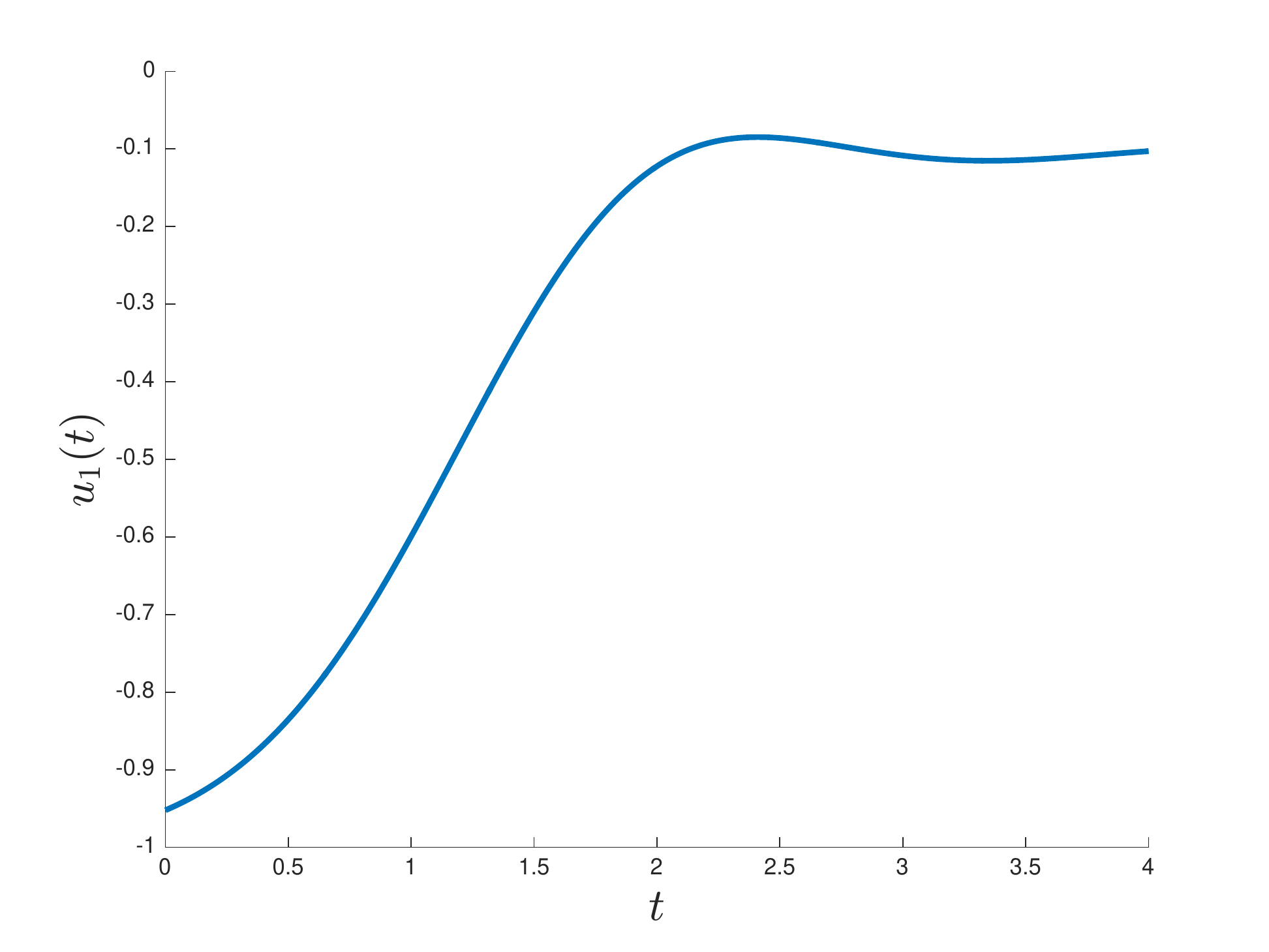}
			\end{center}
			\caption{\label{fig:FK_u1} A validated traveling wave profile of \eqref{eq:Fourth_order_PDE} for
			$a=-0.1$, $\kappa=-2$ and $\gamma = 4.202$. The depicted part of the traveling wave profile
			corresponds to the first component of the connecting orbit between the unstable and stable manifold of $\left(-1,0,0,0 \right)$ and 
			$\left(a,0,0,0\right)$, respectively. 
			The time of flight of the connection between the (un)stable manifolds was $L=4$. 
		}
		\end{figure}	
		\begin{figure}[hbt]
			\begin{center}
				\includegraphics[scale = 0.1]{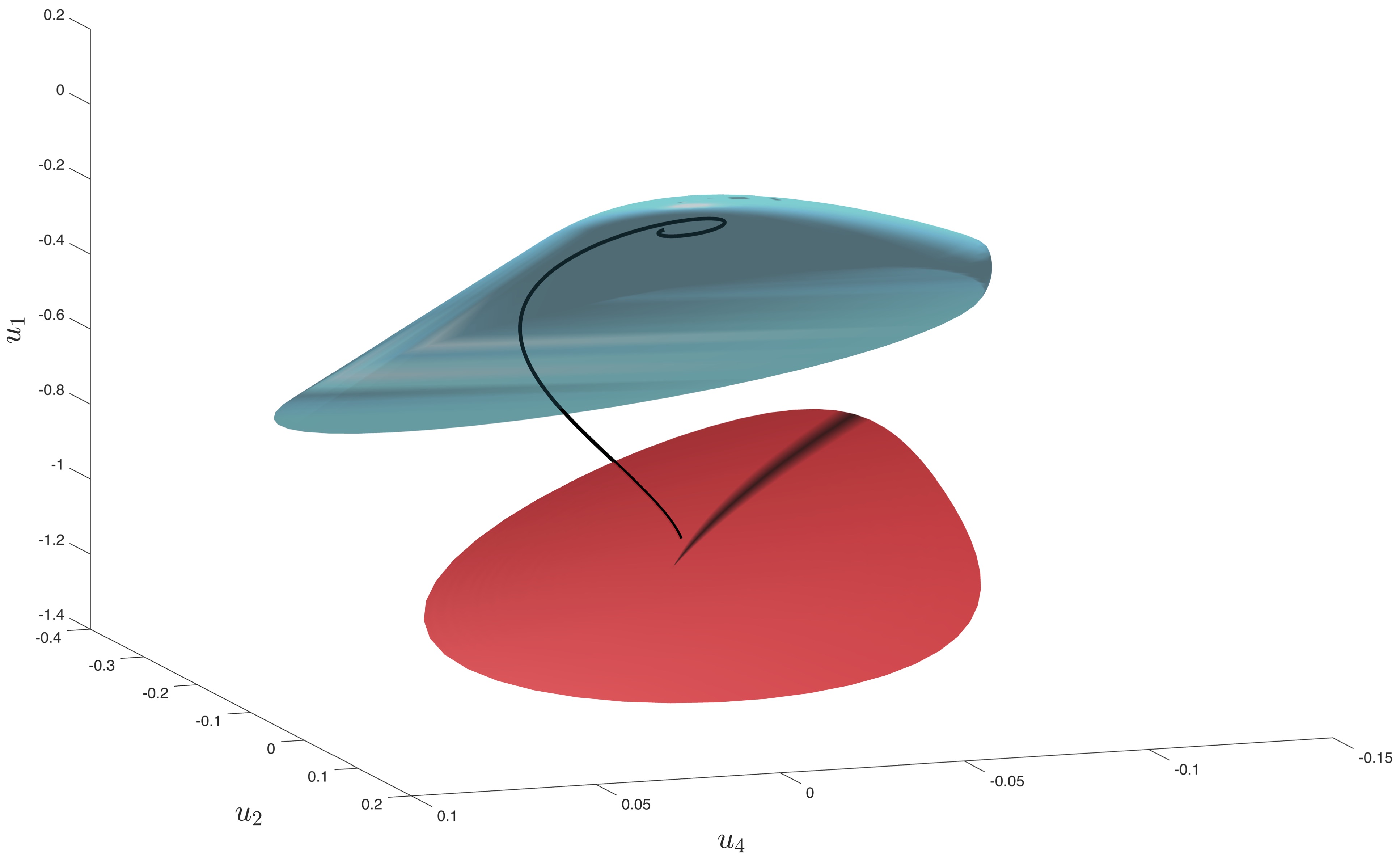}
			\end{center}
			\caption{\label{fig:FK_phase_space} A three dimensional projection of a validated connecting orbit of \eqref{eq:Fourth_order_ODE} for
			$a=-0.1$, $\kappa=-2$ and $\gamma = 4.202$. The geometric objects colored in red and blue correspond to 
			parameterizations of the local unstable and local stable manifold of $\left(-1,0,0,0 \right)$ and $\left(a,0,0,0\right)$, respectively, 
			which were computed (and validated) by using the parameterization method. The curve in black corresponds to the 
			piece of the connecting orbit which was validated by using Chebyshev series.
			The connecting orbit was relatively ``short'' and the time of flight was $L=4$. We remark that
			the integration time was not optimized. 
		}
		\end{figure}
\end{application}

\subsection{Extensions and future work}
\label{s:extensions}

We now discuss possible extensions. 
The first extension is to let the vector field explicitly depend on a parameter and to perform rigorous 
(pseudo-arclength) continuation of connecting orbits.
This involves a relatively straightforward application of the 
uniform contraction principle and a slight modification of the estimates developed in this paper 
(see \cite{Queirolo,MR2630003, MR3125637, suspensionbridge} for instance). 
Furthermore, in order to carry out continuation efficiently, we need to develop algorithms (heuristics) which 
automatically determine near-optimal parameter values for the validation of the charts on the local (un)stable manifolds and the connecting orbit.
More specifically, during continuation it might become necessary to modify the number of Taylor coefficients, 
the size of the charts on the local (un)stable manifolds, the grid on which the connecting orbit is computed, the number of Chebyshev coefficients, or the integration time. 

The second extension involves the incorporation of \emph{resonances}. 
In this paper, we assume that the (un)stable eigenvalues
associated to the (un)stable manifolds satisfy a so-called \emph{non-resonance} condition. This condition is related to the regularity of the chart mappings obtained via the parameterization method. In short, the parameterization method is based on 
constructing a smooth conjugacy (analytic in our case) between the nonlinear flow on the (un)stable manifold and an ``easier'' fully understood model system. One can
choose this model system to be linear, which we do in this paper, if the (un)stable eigenvalues satisfy a non-resonance condition. 
If there are resonant eigenvalues, however, one needs
to use a nonlinear model system instead. This is explained in detail in \cite{ManifoldTaylor}. Generically, one will encounter resonances during continuation. 
Therefore, in order to successfully perform validated continuation, we need to allow for the possibility of resonant eigenvalues 
and modify the current method accordingly as explained in \cite{ManifoldTaylor}.  In particular, we need to develop an algorithm which automatically detects when 
to ``switch'' between the linear and nonlinear model flow during continuation. Furthermore, a careful analysis of the case in which 
a pair of complex conjugate eigenvalues become real (or vice versa) is needed as well. 
After the above extensions have been implemented, one can start developing tools 
for the rigorous study of bifurcations of connecting orbits for nonlinear ODEs. 

The computer-assisted method presented in this paper is implemented in an object oriented framework in \textsc{Matlab} using the
\textsc{Intlab} package \cite{Intlab} for interval arithmetic. A third and useful extension would be to incorporate an extra
degree of freedom into the classes for the connecting orbit so that additional equations and variables can be added (or removed)
in a convenient manner.\ This would facilitate the required modifications for dealing with 
non-polynomial vector fields via automatic differentiation techniques, analyzing connecting orbits in vector fields with symmetry, proving the
existence of homoclinic instead of heteroclinic orbits, and performing bifurcation analysis. 

Finally, we remark that the computational efficiency 
of the current implementation can be improved. For instance, the equations for the parameterizations of the local (un)stable manifolds and 
the connecting orbit in between are to a large extent uncoupled. As a consequence, the derivative of the zero finding map $F$ has a block 
structure, which can be exploited to reduce the computational costs of the computation of an approximate inverse (we need an explicit finite 
dimensional approximate inverse to construct a Newton-like map~$T$). 
Furthermore, in applications it might not be necessary to resolve the full local stable manifold equally well in all directions,
but a rather more focussed parameterization centered around the slow eigendirections of the equilibria is appropriate, 
since a connecting orbit generically tends to enter the stable manifold via these directions. 

\subsection{Outline of the paper}

This paper is organized as follows. In Section \ref{sec:preliminaries} we review some basic facts about Chebyshev series,
Taylor series and sequence spaces, which will be used extensively throughout this paper.
In Section \ref{sec:connection_operator} we set up an equivalent zero finding problem for~\eqref{eq:main_BVP}
by using domain decomposition, the parameterization method, Chebyshev series and Taylor series. 
In Section~\ref{sec:setup} we set up an equivalent fixed-point problem and explain how the existence of a zero can be 
established with the aid of a computer.
This involves the construction of computable bounds which are developed in full detail in Section \ref{sec:bounds}. 
Finally, in Section \ref{sec:applications} we demonstrate the effectiveness of the method by proving the 
existence of traveling fronts in parabolic PDEs. 
We also discuss some algorithmic aspects.

\section{Preliminaries}
\label{sec:preliminaries}

In this section we develop a functional analytic framework for analyzing maps
which arise from the study of connecting orbits. We start in Section \ref{sec:Chebyshev}
by recalling basic results from Chebyshev approximation theory. 
In Sections \ref{sec:seq} and \ref{sec:arrays} we introduce spaces of geometrically 
decaying sequences and multivariate arrays, respectively. In addition, we review 
methods for analyzing bounded linear operators on them.

\subsection{Chebyshev series}
\label{sec:Chebyshev}

In this section we recall basic notions and results from Chebyshev approximation theory. 
The reader is referred to \cite{MR3012510} for the proofs and a more comprehensive introduction into 
the theory of Chebyshev approximations. 
\begin{definition}
	\label{def:chebpoly}
	The Chebyshev polynomials $T_{k}: [-1,1] \rightarrow \RR$ are defined by the relation $T_{k} \left( \cos \left( \theta \right) \right) = \cos \left( k \theta \right)$,
	where $k \in \NN_{0}$ and $\theta \in \left [0,\pi \right]$. 
\end{definition}

Chebyshev series constitute a non-periodic analog of Fourier cosine series and have similar convergence
properties. For instance, any Lipschitz continuous function admits a unique 
Chebyshev expansion. In this paper, we will consider Chebyshev expansions of analytic functions. 
The Chebyshev coefficients of such regular functions decay (in analogy with Fourier series) at a geometric rate to zero. 
A more precise statement is given in the next proposition. 
\begin{proposition}
	\label{prop:geometricDecay}
	Suppose $u: \left[ -1, 1\right] \rightarrow \RR$ is analytic and let 
	\begin{align*}
		u = a_{0} + 2 \sum_{k=1}^{\infty} a_{k} T_{k}
	\end{align*}
	be its Chebyshev expansion. Let $\mathcal{E}_{\nu} \subset \CC$ denote an open ellipse with foci 
	$\pm 1$ to which $u$ can be analytically extended, where $\nu>1$ is the sum of the semi-major 
	and semi-minor axis of $\mathcal{E}_{\nu}$. If $u$ is bounded on $\mathcal{E}_{\nu}$, then  $\left \vert a_{k} \right \vert \leq M \nu^{-k}$ 
	for all $k \in \NN_{0}$, where $M = \sup_{z \in \mathcal{E}_{\nu} } \left \vert u(z) \right \vert$. 
\end{proposition}

The Chebyshev coefficients of the product of two Chebyshev series is (in direct analogy with Fourier cosine series) 
given by the \emph{symmetric} discrete convolution:
\begin{proposition}
	\label{prop:convolutionTheorem}
	Suppose $u,v: \left[-1,1\right] \rightarrow \RR$ are Lipschitz continuous and let 
	\begin{align*}
		u = a_{0} + 2 \sum_{k=1}^{\infty} a_{k} T_{k}, & \quad
		v = b_{0} + 2 \sum_{k=1}^{\infty} b_{k} T_{k}, 
	\end{align*}
	be the associated Chebyshev expansions. Then 
	\begin{align*}
		u \cdot v = \left(a \ast b\right)_{0} + 2 \sum_{k=1}^{\infty} \left(a \ast b\right)_{k} T_{k}, \quad 
		\qquad\text{where}\qquad
		\left( a \ast b \right)_{k} :=
		\mathop{ \sum_{k_{1} + k_{2} = k } }_{ k_{1},k_{2} \in \ZZ}  a_{ \left \vert k_{1} \right \vert} b_{ \left \vert k_{2}\right \vert} .
	\end{align*}
\end{proposition}

\subsection{Geometrically decaying sequences}
\label{sec:seq}

In this section we introduce a sequence space suitable for analyzing analytic functions,
and elementary operations on them, via their Chebyshev coefficients. Recall that the Chebyshev 
coefficients of an analytic function $u : [a,b] \rightarrow \RR^{n}$ decay exponentially fast to zero by 
Proposition \ref{prop:geometricDecay}. In light of this observation we define 
\begin{align*}
	\ell^{1}_{\nu,n} := 
	\left \{ 
		\left( a_{k} \right)_{k \in \NN_{0}} : a_{k} \in \CC^{n},  
		\ \sum_{k=0}^{\infty} \left \vert \left [ a_{k} \right]_{j} \right \vert \nu^{k} < \infty, 
		\ 1 \leq j \leq n 
	\right \},
\end{align*}
where $\left[ a_{k} \right]_{j}$ denotes the $j$-th component 
of $a_{k}$ and $\nu >1$ is some prescribed weight, endowed with the norm 
\begin{align*}	
	\left \Vert a \right \Vert_{\nu,n} := 
	\max_{1 \leq j \leq n} 
	\left \{ 
		\left \vert \left [ a_{0} \right ]_{j} \right \vert + 2 \sum_{k=1}^{\infty} \left \vert \left [ a_{k} \right]_{j} \right \vert \nu^{k}
	\right \}.
\end{align*}
In the special case that $n=1$ we shall write $\ell^{1}_{\nu} := \ell^{1}_{ \nu,1}$
and $ \left \Vert \cdot \right \Vert_{\nu} := \left \Vert \cdot \right \Vert_{\nu,1}$. 
It is a straightforward task to verify that $\ell^{1}_{\nu,n}$ equipped with this norm is a Banach space over $\CC$. 
\begin{remark}
	In this paper we are exclusively concerned with Chebyshev expansions of real-valued functions.
	From this perspective it is more natural to consider sequence spaces over $\RR$ instead of $\CC$ . 
	The reason for using a space of complex valued sequences is that we wish to couple the Chebyshev 
	expansions with chart maps for (un)stable manifolds, which might be complex-valued (see Section 
	\ref{sec:connection_operator}).
	We will proof a-posteriori that the Chebyshev coefficients are real by using arguments based on symmetry. 
\end{remark} 

The operation of multiplying two Chebyshev series can be lifted to the level of sequences, giving rise to
the symmetric discrete convolution $\ast$, as shown in Proposition \ref{prop:convolutionTheorem}. 
This additional product structure on $\ell^{1}_{\nu}$ yields a particularly nice space:

\begin{proposition}
	\label{prop:Banach_algebra}
	The space $\left( \ell^{1}_{\nu}, \ast \right)$ is a commutative Banach algebra over $\CC$. 
	\begin{proof}
		This follows directly from Proposition \ref{prop:convolutionTheorem} and the triangle inequality. 
	\end{proof}
\end{proposition}

One of the reasons for using the space $\ell^{1}_{\nu}$ is to have a relatively simple and sharp convolution estimate.
Another important reason is that it is easy to compute the norm of bounded linear operators. 
To explain how to compute the norm of a bounded linear operator on $\ell^{1}_{\nu}$ we introduce
the notion of the \emph{corner points}. Let $\left( e_{k} \right)_{k \in \NN_{0}}$ denote the canonical 
Schauder basis for $\ell^{1}_{\nu}$, i.e. $\left(e_{k}\right)_{l} := \delta_{kl}$ for $l \in \NN_{0}$, so that 
\begin{align*}
	a = \sum_{k=0}^{\infty} a_{k} e_{k},
\end{align*}
for any $a \in \ell^{1}_{\nu}$. 
\begin{remark}
	We shall frequently use the Schauder basis $\left(e_{k}\right)_{k \in \NN_{0}}$
	to identify an element $a \in \ell^{1}_{\nu}$ with the infinite column vector 
	$\begin{bmatrix} a_{0} & a_{1} & \ldots & \end{bmatrix}^{T}$.
\end{remark}
\begin{definition} 
	\label{def:cornerpoints}
	The corner points $\left \{ \xi_{k,\nu} \right \}_{k \in \NN_{0}} \subset \ell^{1}_{\nu}$ 
	of the unit ball in $\ell^{1}_{\nu}$ are defined by $\xi_{k,\nu}  := \varepsilon_{k,\nu} e_{k}$, where 
	\begin{align*}
		\varepsilon_{k,\nu} = 
		\begin{cases}
			1 & k = 0, \\
			\frac{1}{2} \nu^{-k}, & k \in \NN. 
		\end{cases}
	\end{align*}
	We shall write $\xi_{k,\nu} = \xi_{k}$ and $ \varepsilon_{k,\nu}  = \varepsilon_{k}$ whenever there is
	no chance of confusion. 
\end{definition}

The norm of a bounded linear operator on $\ell^{1}_{\nu}$ can be computed by
simply evaluating it at the corner points as shown in the next proposition:
\begin{proposition}
	\label{prop:l1_operator_norm_X}
	Let $\left(X, \left \Vert \cdot \right \Vert_{X} \right)$ be a normed vector space. If  
	$\mathcal{L} \in \mathcal{B} \left( \ell^{1}_{\nu}, X \right)$, then 
	\begin{align*}
		\left \Vert \mathcal{L} \right \Vert_{ \mathcal{B} \left( \ell^{1}_{\nu}, X \right) } = 
		\sup_{k \in \NN_{0}} \left \Vert \mathcal{L} \left( \xi_{k} \right) \right \Vert_{X} . 
	\end{align*}	
	\begin{proof}
		It is clear that  
		\begin{align*}
			\left \Vert \mathcal{L} \right \Vert_{ \mathcal{B} \left( \ell^{1}_{\nu}, X \right) } \geq
			\sup_{k \in \NN_{0}} \left \Vert \mathcal{L} \left( \xi_{k} \right) \right \Vert_{X}, 
		\end{align*}
		since $\left \Vert \xi_{k} \right \Vert_{\nu} =1$ for all $k \in \NN_{0}$ by definition. 
		
	 	Conversely, let $a \in \ell^{1}_{\nu}$ be arbitrary and observe that 
		\begin{align*}
			a = a_{0} \xi_{0} + 2 \sum_{k=1}^{\infty} a_{k} \xi_{k} \nu^{k}.
		\end{align*}
		Therefore, since $\mathcal{L}$ is bounded, 
		\begin{align*}
			\mathcal{L} \left( a \right)  = 
			a_{0} \mathcal{L} \left( \xi_{0} \right) + 
			2 \sum_{k=1}^{\infty} a_{k} \mathcal{L} \left( \xi_{k} \right) \nu^{k}.
		\end{align*}
		Consequently, 
		\begin{align*}
			\left \Vert \mathcal{L} \left( a \right) \right \Vert_{X} \leq 
			\sup_{k \in \NN_{0}} \left \Vert \mathcal{L} \left( \xi_{k} \right) \right \Vert_{X}
			\left \Vert a \right \Vert_{\nu}
		\end{align*}
		for any $a \in \ell^{1}_{\nu}$, which proves the claim. 
	\end{proof}
\end{proposition}

Now, suppose $\mathcal{L} \in \mathcal{B} \left( \ell^{1}_{\nu_{1}}, \ell^{1}_{\nu_{2}} \right)$, where $\nu_{1}, \nu_{2} >1$. 
Then $\mathcal{L}$ can be identified with an infinite dimensional matrix with respect to the basis $\left(e_{k}\right)_{k \in \NN_{0}}$. 
More precisely, there exists unique coefficients $\left \{ \mathcal{L}_{ij} \in \CC : i,j \in \NN_{0} \right \}$ 
such that  
\begin{align*}
	\mathcal{L} \left( e_{j} \right) = 
	\sum_{i=0}^{\infty} \mathcal{L}_{ij} e_{i}  \simeq
	\begin{bmatrix}
		 \mathcal{L}_{0j} &  \mathcal{L}_{1j} & \ldots &
	\end{bmatrix}^{T}, \quad j \in \NN_{0}. 
\end{align*}
Hence
\begin{align}
	\label{eq:matrix_rep}
	\mathcal{L}(a) = 
	\begin{bmatrix}
		\mathcal{L}_{00} & \mathcal{L}_{01} & \ldots \\[2ex]
		\mathcal{L}_{10} & \mathcal{L}_{11} & \ldots \\[2ex]	
		\vdots & & 	
	\end{bmatrix}
	\begin{bmatrix}
		a_{0} \\[2ex]
		a_{1} \\[2ex]
		\vdots
	\end{bmatrix},
\end{align}
for any $a \in \ell^{1}_{\nu_{1}}$. In this particular setting, Proposition \ref{prop:l1_operator_norm} 
can be interpreted as the statement that $\left \Vert \mathcal{L} \right \Vert_{ \mathcal{B} \left( \ell^{1}_{\nu_{1}}, \ell^{1}_{\nu_{2}} \right) }$ 
is a weighted supremum of the $\ell^{1}_{\nu_{2}}$-norms of the columns of $\mathcal{L}$. Moreover,
in this case the converse of Proposition \ref{prop:l1_operator_norm} holds as well: 
\begin{proposition}
	\label{prop:l1_operator_norm}
	Let $\nu_{1}, \nu_{2} >1$ and 
	suppose $\left \{ \mathcal{L}_{ij} \in \CC : i,j \in \NN_{0} \right \}$ are coefficients such
	that the expression in \eqref{eq:matrix_rep} yields a well-defined linear operator $\mathcal{L} : \ell^{1}_{\nu_{1}}
	\rightarrow \CC^{\NN_{0}}$, i.e., $\left( \mathcal{L}(a) \right)_{k}$ is finite for all $a \in \ell^{1}_{\nu_{1}}$ and $k \in \NN_{0}$.
	Then  
	$\mathcal{L} \in \mathcal{B} \left( \ell^{1}_{\nu_{1}}, \ell^{1}_{\nu_{2}} \right)$
	if and only if 
	$
		\sup_{l \in \NN_{0}} \varepsilon_{l, \nu_{1}}
		\left \Vert \mathcal{L}_{\left( \cdot, l \right) } \right \Vert_{\nu_{2}} < \infty.
	$
	Moreover, if $\mathcal{L} \in \mathcal{B} \left( \ell^{1}_{\nu_{1}}, \ell^{1}_{\nu_{2}} \right)$, then
	\begin{align}
		\label{eq:L_operator_norm}
		\left \Vert \mathcal{L} \right \Vert_{\mathcal{B} \left( \ell^{1}_{\nu_{1}}, \ell^{1}_{\nu_{2}} \right)} = 
		\sup_{l \in \NN_{0}} \varepsilon_{l, \nu_{1}}
		\left \Vert \mathcal{L}_{\left( \cdot, l \right) } \right \Vert_{\nu_{2}}.
	\end{align}
	\begin{proof}
		It follows directly from Proposition \ref{prop:l1_operator_norm} that 
		$
			\sup_{l \in \NN_{0}} \varepsilon_{l, \nu_{1}}
			\left \Vert \mathcal{L}_{\left( \cdot, l \right) } \right \Vert_{\nu_{2}} < \infty
		$
		whenever $\mathcal{L} \in \mathcal{B} \left( \ell^{1}_{\nu_{1}}, \ell^{1}_{\nu_{2}} \right)$ and that
		in this case the operator norm is given by \eqref{eq:L_operator_norm}.
		Conversely, suppose
		$\sup_{l \in \NN_{0}} \varepsilon_{l, \nu_{1}} \left \Vert \mathcal{L}_{\left( \cdot, l \right) } 
		\right \Vert_{\nu_{2}} < \infty$. Let $a \in \ell^{1}_{\nu_{1}}$ be arbitrary, then 
		\begin{align*}
			\left \Vert \mathcal{L} \left( a \right) \right \Vert_{\nu_{2}} & \leq
			 \sum_{l=0}^{\infty} \left \vert \mathcal{L}_{0l} \right \vert \left \vert a_{l} \right \vert + 
			 2 \sum_{k=1}^{\infty} \sum_{l=0}^{\infty} \left \vert \mathcal{L}_{kl} \right \vert 
			 \left \vert a_{l} \right \vert \nu_{2}^{k} \\[2ex] &= 
			 \sum_{l=0}^{\infty} \left \vert a_{l} \right \vert \left \Vert \mathcal{L}_{\left(\cdot, l\right)} \right \Vert_{\nu_{2}}
			 \\[2ex] &= 
			 \left \vert a_{0} \right \vert \left \Vert \mathcal{L}_{\left(\cdot, 0\right)} \right \Vert_{\nu_{2}}
			 + 2 \sum_{l=1}^{\infty}  \varepsilon_{\nu_{1},l} \left \vert a_{l} \right \vert 
			  \left \Vert \mathcal{L}_{\left(\cdot, l\right)} \right \Vert_{\nu_{2}} \nu_{1}^{l} \\[2ex] &\leq
			  \left( \sup_{l \in \NN_{0}} \varepsilon_{l, \nu_{1}} \left \Vert \mathcal{L}_{\left( \cdot, l \right) } \right \Vert_{\nu_{2}} \right)
			  \left \Vert a \right \Vert_{\nu_{1}}. 
		\end{align*}
		Since the first factor in the right-hand side is  finite, $\mathcal{L} : \ell^{1}_{\nu_{1}} \rightarrow \ell^{1}_{\nu_{2}}$ is thus a 
		bounded linear operator. 
	\end{proof}	
\end{proposition}

\subsection{Multivariate sequences}
\label{sec:arrays}
In this section we introduce a space of sequences indexed by $d$-dimensional multi-indices, where $d \in \NN$. 
This space will be used to analyze Taylor series of analytic functions $P: \left \{ z \in \CC : \left \vert z \right \vert \leq \nu \right \}^{d} \rightarrow \CC^{n}$,
where $\nu > 0$. Such functions arise in the the analysis of local charts on (un)stable manifolds via the parameterization method developed in \cite{MR2177465}. 

Formally, a sequence indexed by $d$-dimensional multi-indices is a function $p: \NN_{0}^{d} \rightarrow \CC$. 
The function $p$ is usually referred to as a $d$-dimensional \emph{array} or \emph{multivariate sequence}. In analogy
with ordinary sequences, we shall write (as usual) 
\begin{align*}
	p_{k} := p_{k_{1} \ldots k_{d}} := p \left( k \right), \quad k \in \NN_{0}^{d}.
\end{align*}
Furthermore, for any multi-index $k \in \NN_{0}^{d}$, we write $\left \vert k \right \vert  = \sum_{i=1}^{d} k_{i}$, which is 
not to be confused with the absolute value of a (complex) number. In addition, we introduce a partial ordering $\preceq$ on 
$\NN_{0}^{d}$ by 
\begin{align*}
	k \preceq l \quad \mathop{\Leftrightarrow}^{\bold{def}} \quad k_{j} \leq l_{j}, \ \forall \ 1 \leq j \leq d. 
\end{align*}

We will now follow the same approach as in the previous section to set up a functional analytic framework
for analyzing geometrically decaying arrays.
Let $\nu > 0$ and define 
\begin{align*}
	W^{1}_{\nu,n,d} := 
	\left \{ 
		p : \NN_{0}^{d} \rightarrow \CC^{n} 
		\ \middle | \ \sum_{k \in \NN_{0}^{d}} 
		\left \vert \left[ p_{k} \right]_{j} \right \vert \nu^{\left \vert k \right \vert} < \infty, \ 1 \leq j \leq n
	 \right \}
\end{align*}
endowed with the norm 
\begin{align*}
	\left \Vert p \right \Vert_{W^{1}_{\nu,n,d}} := \max_{1 \leq j \leq n} 
	\sum_{k \in \NN_{0}^{d}} \left \vert \left[ p_{k} \right]_{j} \right \vert \nu^{\left \vert k \right \vert}.
\end{align*}
In the case that the dimension $d$ can be easily inferred from the context it will be omitted from 
the notation. In addition, if $n=1$ and it is clear from the context whether $p \in \ell^{1}_{\nu}$ or $p \in W^{1}_{\nu}$,  
we shall write $\left \Vert p \right \Vert_{W^{1}_{\nu,n,d}} = \left \Vert p \right \Vert_{\nu}$.

Next, recall that the Taylor coefficients of the product of two Taylor series is given by the \emph{one-sided} discrete convolution, also referred to 
as the \emph{Cauchy product}. More precisely, if $f,g: \left \{ z \in \CC : \left \vert z \right \vert < \nu \right \}^{d} \rightarrow \CC$ admit
power series expansions  
\begin{align*}
	f(z) = \sum_{k \in \NN_{0}^{d}} f_{k} z^{k}, \quad g(z) = \sum_{k \in \NN_{0}^{d}} g_{k} z^{k},
\end{align*}
where $\tilde f= \left( f_{k} \right)_{k \in \NN_{0}^{d}}$ and $\tilde g= \left( g_{k} \right)_{k \in \NN_{0}^{d}}$ are $d$-dimensional arrays, then 
\begin{align}
	(fg)(z) = \sum_{k \in \NN_{0}^{d}} \left( \tilde f \cp \tilde g \right)_{k} z^{k}, \quad
	\left( \tilde f \cp \tilde g \right)_{k} := \sum_{ \substack{ \alpha + \beta = k, \\ \alpha, \beta \in \NN_{0}^{d}} }
	\tilde f_{\alpha} \tilde g_{\beta},
	\label{eq:multiconv}
\end{align}
on $\left \{ z \in \CC : \left \vert z \right \vert < \nu \right \}^{d}$. In particular, the Cauchy product $\cp$ yields
a natural product structure on $W^{1}_{\nu}$. This is summarized in the following proposition. 
\begin{proposition}
	The space $\left( W^{1}_{\nu}, \cp \right)$ is a commutative Banach algebra.
	\begin{proof}
			This follows directly from the definition of $\cp$ in \eqref{eq:multiconv} and the triangle inequality. 
	\end{proof}
\end{proposition}

Next, we derive an expression for the norm of a bounded linear operator on $W^{1}_{\nu}$. For this purpose 
we introduce a multivariate analog of the corner-points: 

\begin{definition} 
	\label{def:cornerpoints_array}
	The corner-points $\left \{ \bm{\xi}_{k,d,\nu} \right \}_{k \in \NN^{d}_{0}} \subset W^{1}_{\nu}$ 
	of the unit ball in $W^{1}_{\nu}$ are defined by 
	$\left( \bm{\xi}_{k,d,\nu} \right)_{l}  :=  \nu^{-\left \vert k \right \vert} \delta_{kl}$, where $l \in \NN_{0}^{d}$. 
	We shall write $\bm{\xi}_{k,d,\nu} = \bm{\xi}_{k}$ whenever there is no chance of confusion. 
\end{definition}

As before, the norm of a bounded linear operator on $W^{1}_{\nu}$ can be computed by evaluating it at the corner-points. 
\begin{proposition}
	\label{prop:operator_norm_mult}
	Let $\left(X, \left \Vert \cdot \right \Vert_{X} \right)$ be a normed vector space. If  $\mathcal{L} \in \mathcal{B} \left( W^{1}_{\nu}, X \right)$, then 
	\begin{align*}
			\left \Vert \mathcal{L} \right \Vert_{ \mathcal{B} \left( W^{1}_{\nu}, X \right) } = 
			\sup_{k \in \NN^{d}_{0}} \left \Vert \mathcal{L} \left( \bm{\xi}_{k} \right) \right \Vert_{X} . 
	\end{align*}
	\begin{proof}
			See Proposition \ref{prop:l1_operator_norm_X}. 
	\end{proof}
\end{proposition}

\section{An equivalent zero finding problem}
\label{sec:connection_operator}
In this section we set up a zero finding problem for establishing the existence of connecting orbits. Let us 
start by giving a precise description of the problem. Suppose $\tilde p_{0}, \tilde q_{0} \in \RR^{n}$ are \emph{hyperbolic} 
equilibria of $g$. The objective is to validate an \emph{isolated} connecting 
orbit $u$ from $\tilde p_{0}$ to $\tilde q_{0}$, which is robust with respect to ``small'' perturbations in $g$, by solving a boundary value
problem (BVP) on a \emph{finite} time domain. The method is based on the observation that a connecting orbit from 
$\tilde p_{0}$ to $\tilde q_{0}$ is characterized by 
\begin{align}
	\label{eq:bvp1}
	\begin{cases}
		\dfrac{du}{dt} = g(u), & t \in [0,L], \\[2ex]
		u(0) \in W_{\text{loc}}^{u} \left( \tilde p_{0} \right), \\[1ex]
		u(L) \in W_{\text{loc}}^{s} \left( \tilde q_{0} \right),
	\end{cases}
\end{align} 
where $L>0$ is the time of flight needed to travel from $W_{\text{loc}}^{u} \left( \tilde p_{0} \right)$ to $W_{\text{loc}}^{s} \left( \tilde q_{0} \right)$.

If $W^{u} \left( \tilde p_{0} \right)$ and $W^{s} \left( \tilde q_{0} \right)$ intersect transversally along $u$,
then the connecting orbit is robust, i.e., it will persist for sufficiently ``small'' perturbations in $g$. 
In this case, the intersection $W^{u} \left( \tilde p_{0} \right) \cap W^{s} \left( \tilde q_{0} \right) \cap U$,
where $U$ is a neighborhood of the connecting orbit in which it is unique, is necessarily an one dimensional manifold. 
Hence, by counting dimensions, a necessary condition for the 
existence of a transverse isolated connecting orbit is 
\begin{align*}
	n_{u} + n_{s} - n = 1, \quad n_{u} := \dim W^{u} \left( \tilde p_{0} \right), \ n_{s} := \dim W^{s} \left( \tilde q_{0} \right).
\end{align*}
This condition is often referred to as a \emph{non-degeneracy} condition for connecting orbits.
We shall henceforth assume that this condition is satisfied. In particular, we do not assume
a-priori that the connecting orbit is isolated and transverse. Instead, we will obtain these properties from the proof of existence
(a contraction argument), see Proposition \ref{prop:symmetry_connection}.  

We start by setting up equations for local charts on the (un)stable manifolds by using the parameterization method 
\cite{MR2177465} and the methodology presented in 
\cite{ManifoldTaylor}. These charts will be used to supplant the boundary conditions in \eqref{eq:bvp1} with explicit equations. Next, we set up an equivalent system of equations for the differential equation by using Chebyshev series and domain decomposition as explained in 
\cite{domaindecomposition}. Finally, in order for the resulting zero finding problem to be well posed, we complete the system of equations by 
imposing appropriate phase conditions. 

\subsection{Charts on the (un)stable manifolds}
\label{s:charts}
In this section we give a brief overview of the method developed in \cite{ManifoldTaylor} to compute
local charts on the (un)stable manifolds. The reader is referred to \cite{ManifoldTaylor} for a more detailed 
exposition of the theory. 
We consider the computation of a local chart on the stable manifold of $\tilde q_{0}$. A 
chart on the unstable manifold of $\tilde p_{0}$ can be computed in the same way by reversing the sign of the vector field. 

\paragraph{The Parameterization Method}
The idea of the parameterization method \cite{MR2177465} 
is to construct a diffeomorphism which conjugates the nonlinear dynamics on the stable 
manifold to an easier and fully understood flow $\psi$. For the sake of simplicity, let us assume that $Dg \left( \tilde q_{0} \right)$ is diagonalizable. This assumption is, however, not necessary, as we will explain in a moment. 

Let $\lambda_{1}^{s}, \ldots, \lambda^{s}_{n_{s}} \in \CC$ be the stable eigenvalues of $Dg \left( \tilde q_{0} \right)$. If all eigenvalues are 
real and semisimple, then there exists neighborhoods $U \subset \RR^{n}$ and $V \subset \RR^{n_{s}}$ of $\tilde q_{0}$ and $0$, respectively, 
such that the dynamics on $W^{s} \left( \tilde q_{0} \right) \cap U$ is conjugate to the flow 
\begin{align}
	\label{eq:linear_flow}
	\psi \left(t,\phi \right) := \exp \left( t \cdot \mbox{diag} \left( \lambda_{1}^{s}, \ldots, \lambda^{s}_{n_{s}} \right) \right) \phi,
	\quad t \geq 0, \ \phi \in V. 
\end{align}
If some of the eigenvalues are complex, however, special care has to be taken. Let us for the moment
forget about this technicality and consider the \emph{complex} dynamics generated by $u'=g(u)$ on $\CC^{n}$. Then the dynamics
on the complex local stable manifold, which we denote by $W_{\text{loc}}^{s,c} \left( \tilde q_{0} \right)$, is conjugate to the flow 
$\psi$ restricted to the polydisk 
\begin{align*}
	\mathbb{B}_{\nu_{s}} := \left \{ \phi \in \CC^{n_{s}}: \max_{1 \leq i \leq n_{s} } \left \vert \phi_{i} \right \vert \leq \nu_{s} \right \}, 
\end{align*}
for some sufficiently small $\nu_{s}>0$.

The idea is to find an \emph{analytic} map $Q : \mathbb{B}_{\nu_{s} } \rightarrow \CC^{n}$ which conjugates the nonlinear flow 
$\varphi$ on $W^{s,c}_{\text{loc}} \left( \tilde q_{0} \right)$ to the linear flow $\psi$ on $\mathbb{B}_{\nu_{s}}$ for $t \geq 0$. In other words, we seek a map $Q$ such that 
\begin{center}
	\begin{tikzpicture}[baseline= (a).base]
		\node[scale=1.05] (a) at (0,0){
		\begin{tikzcd}
			\mathbb{B}_{\nu_{s}} \arrow[d,"\psi"] \arrow[r, "Q"] &
			\CC^{n} \arrow[d,"\varphi"] \\[1ex]
			\mathbb{B}_{\nu_{s}} \arrow[r,"Q"] & \CC^{n}
		\end{tikzcd}};
	\end{tikzpicture}
\end{center}
commutes, i.e., $Q \left( \psi \left(t,\phi \right) \right) = \varphi \left(t, Q \left( \phi \right) \right)$ for all $(t,\phi) \in \RR_{\geq 0} \times \mathbb{B}_{\nu_{s}}$. Differentiation of
this relation at $t=0$ yields the so-called \emph{invariance equation}:
\begin{align}
	\label{eq:invariance}
	D Q(\phi) \cdot \mbox{diag} \left( \lambda_{1}^{s}, \ldots, \lambda^{s}_{n_{s}} \right) \phi = g \left( Q(\phi) \right), \quad \phi \in \mathbb{B}_{\nu_{s}}. 
\end{align}
Note that this equation does not depend on time anymore. Moreover, it is easy to see that if the invariance equation holds, then 
\begin{align*}
	t \mapsto u(t) := Q \left( \psi \left(t, \phi \right) \right)
\end{align*}
is an orbit in $W^{s,c} \left( \tilde q_{0} \right)$ for any $\phi \in \mathbb{B}_{\nu_{s}}$, i.e., 
$Q : \mathbb{B}_{\nu_{s} } \rightarrow W^{s,c}_{\text{loc}} \left( \tilde q_{0} \right)$
(see \cite[Lemma $2.6$]{ManifoldTaylor}). Therefore, the problem of computing a chart is now reduced 
to solving \eqref{eq:invariance}. 

\paragraph{Solving the invariance equation}
Since $Q$ is assumed to be analytic on $\mathbb{B}_{\nu_{s}}$, i.e., $Q$ is analytic on a slightly larger open 
neighborhood of $\mathbb{B}_{\nu_{s}}$, there exist coefficients $q \in W^{1}_{\nu_{s},n,n_{s}}$ such that 
\begin{align*}
	Q \left( \phi \right) = \sum_{k \in \NN_{0}^{d}} q_{k} \phi^{k}.
\end{align*}
Observe that the zeroth order Taylor coefficient is necessarily the equilibrium, i.e., $q_{0} = \tilde q_{0}$. 
Furthermore, since $Q$ is assumed to be a diffeomorphism, it must hold that 
\begin{align*}
	D Q \left(0\right) \CC^{n_{s}} = T_{ q_{0} } W_{\text{loc}}^{s,c} \left( q_{0} \right) = E_{s},
\end{align*}
where $E_{s}$ is the stable eigenspace of $Dg \left( q_{0} \right)$. In other words, the first order Taylor coefficients 
$\left \{ q_{k} : \left \vert k \right \vert = 1,\  k \in \NN_{0}^{n_{s}} \right \}$ are the eigenvectors of $Dg \left( q_{0} \right)$. 
Note that these are only determined up to a scaling. 

To determine the higher order Taylor coefficients $\left \{ q_{k} : \left \vert k \right \vert \geq 2, \ k \in \NN_{0}^{n_{s}} \right \}$, we 
first introduce the map $C : W^{1}_{\nu_{s},n} \rightarrow W^{1}_{\nu_{s},n}$ defined by 
\begin{align}
	\label{eq:C(a)}
	C \left( w \right) := \sum_{\alpha \in \mathcal{A} } 
	g_{\alpha} w^{\alpha}, \quad w^{\alpha} := \prod_{j=1}^{n} \left[ w_{j} \right]^{\alpha_{j} },
\end{align}
where the latter product is understood to be the one-sided discrete convolution. Furthermore, $\mathcal{A} \subset \NN_{0}^{n}$ and
\mbox{$\left \{ g_{\alpha} : \alpha \in \mathcal{A} \right \} \subset \RR^{n}$}
are the coefficients of $g$ in the monomial basis. In particular, observe that 
\begin{align}
	\label{eq:conv_taylor}
	g \left( Q \left( \phi \right) \right) = \sum_{k \in \NN_{0}^{n_{s}}} C_{k} \left( q \right) \phi^{k},
\end{align}
since the Taylor coefficients of the product of two Taylor expansions is given by the one-sided discrete convolution. Formally, we should
incorporate the weight $\nu_{s}$ and dimension $n_{s}$ into the notation for $C$. However, since these parameters can usually be
inferred from the context and we wish to use the same notation for the unstable manifold, we have chosen to omit them from the notation. 

Substitution of the Taylor expansion for $Q$ into \eqref{eq:invariance} yields the following system of equations: 
\begin{align}
	\label{eq:invariance2}
	\left \langle \lambda^{s}, k \right \rangle q_{k} - C_{k} \left( q \right)  = 0, \quad
	\lambda^{s} := \begin{bmatrix} \lambda^{s}_{1} \ldots \lambda^{s}_{n_{s}} \end{bmatrix}^{T}, \quad
	\left \vert k \right \vert \geq 2,
\end{align}
where $\left \langle \cdot, \cdot \right \rangle$ denotes the standard Hermitian inner product on $\CC^{n_{s}}$. 
We shall use this system of equations to set up a zero finding problem for computing a chart on 
$W^{s,c}_{\text{loc}} \left( \tilde q_{0} \right)$. Before we proceed, 
observe that \eqref{eq:invariance2} is equivalent to 
\begin{align}
	\label{eq:resonance}
	\left[ Dg \left( q_{0} \right) - \left \langle \lambda^{s}, k \right \rangle I \right]q_{k} = Dg \left( q_{0} \right) q_{k} - C_{k}(q), \quad \left \vert k \right \vert \geq 2. 
\end{align}
Moreover, 
differentiation of \eqref{eq:conv_taylor} at $\phi =0$ shows that $Dg \left( q_{0} \right) q_{k} = C_{k}(q)$ for $\left \vert k \right \vert =1$. 
Hence \eqref{eq:resonance} reduces to the eigenvalue/eigenvector equation for $Dg \left( q_{0} \right)$ for $\left \vert k \right \vert =1$. Similarly, repeated 
differentiation of \eqref{eq:conv_taylor} at $\phi=0$ shows that the right-hand-side of \eqref{eq:resonance} only depends on Taylor coefficients of order strictly below $\left \vert k \right \vert$. 
In conclusion, the Taylor coefficients $q$ can be computed recursively up to any desired order provided 
\begin{equation}\label{e:nonresonance}
\left \langle \lambda^{s}, k \right \rangle \not = \lambda^{s}_{i} 
\quad \text{for all } 1 \leq i \leq n_{s} \text{ and }
\left \vert k \right \vert \geq 2. 
\end{equation}

The latter condition is usually referred to as a \emph{non-resonance} condition and is related to the regularity of $Q$. More precisely, in the presence of a resonance, 
the parameterization method, as applied above with the \emph{linear} ``model'' flow $\psi$, does not yield an analytic conjugation $Q$. 
It is explained in \cite{ManifoldTaylor} how to construct an analytic conjugation in the present of a resonance.  
The idea is to use a \emph{nonlinear} normal form for $\psi$ instead of just the linear flow in \eqref{eq:linear_flow}. 
The interested reader is referred to \cite{ManifoldTaylor} for a detailed exposition of the resonant case. For the sake of presentation, however, we shall 
assume throughout this paper that there are no resonances. 

We are now ready to set up a zero finding problem for computing the Taylor coefficients of $Q$ (which include the equilibrium and 
eigenvectors) and the eigenvalues $\lambda^{s}$:
\begin{definition}[Taylor map for stable manifolds]
	Let $0 < \tilde \nu_{s} < \nu_{s}$ be given weights. The Taylor map $F_{Q}: \CC^{n_{s}} \times W^{1}_{\nu_{s}, n} \rightarrow W^{1}_{\tilde \nu_{s}, n}$ 
	for stable manifolds is defined by 
	\begin{align*}
		\left( F_{Q} \left( \lambda^{s}, q \right) \right)_{k}:= 
		\begin{cases}
			g \left( q_{0} \right), & k = 0, \\[2ex]
			\left[ Dg \left( q_{0} \right) - \left \langle \lambda^{s}, k \right \rangle I \right]q_{k}, & \left \vert k \right \vert = 1, \\[2ex]
			\left \langle \lambda^{s}, k \right \rangle q_{k} - C_{k} \left( q \right), & \left \vert k \right \vert \geq 2. 
		\end{cases}
	\end{align*}
\end{definition}
\begin{remark}
	The latter map is well-defined since 
	$ \left( \left \langle \lambda^{s}, k \right \rangle q_{k} \right)_{k \in \NN_{0}^{n_{s}}} \in W^{1}_{\tilde 	\nu_{s}, n} $ 
	for any $q \in W^{1}_{\nu_{s},n}$ and $0 < \tilde \nu_{s} <\nu_{s}$. 
\end{remark}
\begin{remark}
	\label{remark:scaling}
	If $\left( \lambda^{s}, q \right)$ is a 
	zero of $F_{Q}$, then so is $\left( \lambda^{s}, \mu q \right)$, where $\mu \in \CC^{n_{s}}$ and
	$
		\left( \mu q \right)_{k} := \mu^{k} q_{k} ,
	$ see \cite[Lemma $2.2$]{ManifoldTaylor}. We will get rid of this extra degree
	of freedom by fixing the orientation and length of the eigenvectors. In particular, observe that 
	the scaling of the eigenvectors (and in turn the ``scaling'' of $q$) determines 
	the decay rate of the coefficients $q$ and hence the size of $\mathbb{B}_{\nu_{s}}$. 
	In effect, the length of the eigenvectors determine (roughly speaking) the ``size'' of the patch on $\normalfont W_{\text{loc}}^{s,c} \left( q_{0} \right)$ 
	parameterized by $Q$. The interested reader is referred to \cite{MR3437754} for a more thorough explanation where
	this phenomena is explored in detail. 
\end{remark}

The zero finding problem for the computation of a chart on the unstable manifold is set up in an analogous way. 
For the sake of completeness (and introducing notation) let us explicitly 
state the assumptions and the associated zero finding map. We assume that $Dg \left( \tilde p_{0} \right)$ is diagonalizable and that the associated eigenvalues $\lambda^{u} \in \CC^{n_{u}}$ satisfy the non-resonance condition~\eqref{e:nonresonance}. The goal is to compute a parameterization 
$P: \mathbb{B}_{\nu_{u}} \subset \CC^{n_{u}} \rightarrow W^{u,c}_{\text{loc}} \left( \tilde p_{0} \right)$ of the form 
$P(\theta) = \sum_{k \in \NN_{0}^{n_{u}}} p_{k} \theta^{k}$, where $\nu_{u} >0$,
by finding a zero of the following map: 
\begin{definition}[Taylor map for unstable manifolds]
	Let $0 < \tilde \nu_{u} <  \nu_{u}$ be given weights. The Taylor map $F_{P}: \CC^{n_{u}} \times W^{1}_{\nu_{u}, n} \rightarrow W^{1}_{\tilde \nu_{u}, n}$ 
	for unstable manifolds is defined by 
	\begin{align*}
		\left( F_{P} \left( \lambda^{u}, p \right) \right)_{k}:= 
		\begin{cases}
			g \left( p_{0} \right), & k = 0, \\[2ex]
			\left[ Dg \left( p_{0} \right) - \left \langle \lambda^{u}, k \right \rangle I \right]p_{k}, & \left \vert k \right \vert = 1, \\[2ex]
			\left \langle \lambda^{u}, k \right \rangle p_{k} - C_{k} \left( p \right), & \left \vert k \right \vert \geq 2. 
		\end{cases}
	\end{align*}
\end{definition}

\paragraph{Symmetry}
In the preceding exposition we considered the complex dynamical system $u' = g(u)$ on $\CC^{n}$. 
Our main interest, however, is the computation of  invariant manifolds in the real-valued dynamical system 
on $\RR^{n}$. We will now explain how we can recover charts for the (un)stable 
manifolds in the real system from the complex ones through the use of symmetry. We remark that one could also have set up the 
parameterization method in the real-valued setting from the start. However, in that case, we would have had to separate the cases 
between the presence of complex eigenvalues and a completely real spectrum. It is in our opinion more convenient from both a 
practical and theoretical point of view to develop a unified approach.
 
Let us consider the stable manifold again. Observe that complex eigenvalues will always appear
in conjugate pairs, since $g$ is \emph{real}-analytic. Suppose there are $d_{s}$ complex conjugate pairs of eigenvalues 
and $n_{s} - 2d_{s}$ real ones. Furthermore, assume that we have ordered the eigenvalues $\lambda^{s}$ in such a way that
$\lambda^{s}_{i} = \overline{ \lambda^{s}_{i+1} }$ for $i \in \left \{2l+1: 0 \leq l \leq d_{s}-1 \right \}$.
Next, define the map $\Sigma : \CC^{n_{s}} \rightarrow \CC^{n_{s}}$ by 
\begin{align}
	\Sigma \left( z_{1}, \ldots, z_{2d_{s}}, z_{2d_{s}+1}, \ldots, z_{n_{s}} \right) := 
	\left( \overline{z_{2}}, \overline{z_{1}}, \ldots, \overline{z_{2d_{s}}}, \overline{z_{2d_{s}-1} }, z_{2d_{s}+1}, \ldots, z_{n_{s}} \right),
	\label{e:defSigma}
\end{align}
and note that $\Sigma$ is an involution on $\CC^{n_{s}}$. For this reason we shall frequently write $z^{\star}:= \Sigma(z)$.
In particular, note that we ordered the stable eigenvalues in such a way that
$\left( \lambda^{s} \right)^{\star} = \lambda^{s}$. Finally, we extend this notion of involution to $W^{1}_{\nu_{s},n}$ by defining 
\begin{align*}
	\left( q^{\star} \right)_{k} := \overline{q}_{k^{\star}}, \quad k \in \NN_{0}^{n_{s}}. 
\end{align*}
It can be readily seen from~\eqref{e:defSigma} that $k^{\star} = \Sigma(k)$ is simply a permutation of the multi-index~$k$. 

The key observation for obtaining charts for the real manifolds is stated in the following lemma. The proof can be found in 
\cite[Lemma $2.1$]{ManifoldTaylor}. 
\begin{proposition}
	\label{prop:real_manifold}
	If $q \in W^{1}_{\nu_{s},n}$ is symmetric, i.e., $q^{\star} = q$, then the map $Q: \mathbb{B}_{\nu_{s}} \rightarrow \CC^{n}$ defined by 
	\begin{align*}
		Q \left( \phi \right) := \sum_{k \in \NN_{0}^{n_{s}}} q_{k} \phi^{k}
	\end{align*}
	is real valued on the set 
	$\normalfont \mathbb{B}^{\text{sym}}_{\nu_{s}} := \left \{ \phi \in \mathbb{B}_{\nu_{s}} : \phi^{\star} = \phi \right \}$. In addition, if $\lambda^{s} \in \CC^{n_{s}}$ is
	symmetric and $F \left( \lambda^{s}, q \right) = 0$, then $\normalfont Q \vert_{\mathbb{B}^{\text{sym}}_{\nu_{s}}}$ is a parameterization of the real stable manifold 
	$\normalfont W^{s}_{\text{loc}} \left( q_{0} \right)$. 
\end{proposition}
\begin{remark}
	\label{remark:real_manifold}
	Note that $\normalfont \mathbb{B}^{\text{sym}}_{\nu_{s}}$ is a real manifold of dimension $n_{s}$. More precisely, 
	we can identify $\normalfont \mathbb{B}^{\text{sym}}_{\nu_{s}}$ with the (real) manifold
	\begin{align*}
		\normalfont \mathbb{B}^{\text{sym},\text{re}}_{\nu_{s}} := 
		\left \{ \phi \in \RR^{n_{s}} : \left \vert \phi_{2j-1} \right \vert^{2} + \left \vert \phi_{2j} \right \vert^{2} \leq \nu_{s} ,
		\ 1 \leq j \leq d_{s}, \ \left \vert \phi_{j} \right \vert \leq \nu_{s}, \ 2d_{s} + 1 \leq j \leq n_{s}  \right\}
	\end{align*}	
	by using the (linear) map $\iota_{s} : \RR^{n_{s}} \rightarrow \CC^{n_{s}}$ defined by 
	\begin{align*}
		\iota_{s} \left( \phi \right) := \left( \phi_{1} + i \phi_{2}, \phi_{1} - i \phi_{2}, \ldots, \phi_{2d_{s}-1} + i \phi_{2d_{s}},
								\phi_{2d_{s}-1} - i \phi_{2d_{s}}, \phi_{2d_{s}+1}, \ldots, \phi_{n_{s}}  \right).
	\end{align*}
\end{remark}
\begin{remark}
	Strictly speaking, when $q^{\star}=q$ the assumption that $\lambda^{s}$ is symmetric is not necessary. To see this, 
	let $e^{s}_{i} \in \CC^{n_{s}}$ be the unit vector defined by $\left[ e^{s}_{i} \right]_{j} = \delta_{ij}$, where $1 \leq i,j \leq n_{s}$.
	If $F \left( \lambda^{s}, q \right) = 0$ for some $\lambda^{s} \in \CC^{n_{s}}$, then 
	\begin{align*}
		Dg \left( q_{0} \right) q_{k} = \lambda^{s}_{k} q_{k}, \quad \left \vert k \right \vert =1. 
	\end{align*}
	In particular, if we take the complex conjugate of the lefthand-side of the above expression for $k=e^{s}_{2j-1}$, we obtain
	\begin{align*}
		\overline{ Dg \left( q_{0} \right) q_{e^{s}_{2j-1}} } = 
		Dg \left( q_{0} \right) \overline{ q_{e^{s}_{2j-1}} } = 
		Dg \left( q_{0} \right) q_{e^{s}_{2j}} = 
		\lambda^{s}_{2j} q_{e^{s}_{2j}}, \quad 1 \leq j \leq d_{s},
	\end{align*}
	since $g$ is real-analytic and $q^{\star} = q$ by assumption. On the other hand, 
	\begin{align*}
		\overline{ Dg \left( q_{0} \right) q_{e^{s}_{2j-1}} } = 
		\overline{ \lambda^{s}_{2j-1} q_{e^{s}_{2j-1} } } = 
		\overline{ \lambda^{s}_{2j-1} } q_{e^{s}_{2j}}, 
		\quad 1 \leq j \leq d_{s}.
	\end{align*}Hence it follows that $\left( \lambda^{s} \right)^{\star} = \lambda^{s}$. 
\end{remark}

In conclusion, in order to conclude that a point  lies on the real stable manifold, it suffices to verify that 
$q^{\star} = q$. It is explained in Section \ref{sec:setup} how we can verify this in practice. For now, let us mention that the
verification is based on the following observation whose proof can be found in \cite[Lemma $4.1$]{ManifoldTaylor}:
\begin{lemma}
	\label{lemma:symmetry_F_Q}
	The map $F_{Q}$ is compatible with $\star$, i.e.,  $F_{Q}\left( \left( \lambda^{s} \right)^{\star}, q^{\star} \right) = F_{Q} 		               
	\left( \lambda^{s}, q \right)^{\star}$
	for any $\left(\lambda^{s}, q \right) \in \CC^{n_{s}} \times W^{1}_{\nu_{s},n}$. 
\end{lemma}

Analogous results hold for the parameterization of the unstable manifold of $\tilde p_{0}$. To avoid clutter 
in the notation we shall denote the involution associated to the unstable manifold of $\tilde p_{0}$ by $\star$ as well.

\subsection{Chebyshev series and domain decomposition} 
\label{sec:chebyshev_domaindecomp}
In this section we follow the strategy in \cite{domaindecomposition} to recast the differential equation into an equivalent zero 
finding problem on $\ell^{1}_{\nu,n}$. The reader is referred to \cite{domaindecomposition} for the details. 
Let $\mathcal{P}_{m} := \left \{ t_{0} = 0 < t_{1} < \ldots < \right.$ $\left. t_{m}  = 1 \right \}$, where $m \in \NN$, be any partition of $[0,1]$.
Then the differential equation in \eqref{eq:bvp1} is equivalent to 
\begin{align*}
	\left(\bold{P}_{1} \right) 
	\begin{cases}
		\dfrac{ d u_{1} }{ dt } = Lg \left( u_{1} \right), & t \in \left[0, t_{1} \right], 
	\end{cases} \quad
	\left(\bold{P}_{i} \right) 
	\begin{cases}
		\dfrac{ d u_{i} }{ dt } = Lg \left( u_{i} \right), & t \in \left[t_{i-1}, t_{i} \right], \\[2ex]
		 u_{i} \left( t_{i-1} \right) = u_{i-1} \left( t_{i-1} \right),
	\end{cases}, \quad	
\end{align*}
where $2 \leq i \leq m$. 
These boundary conditions are thus imposed on the internal nodes of the partition $\mathcal{P}_m$ only. The boundary conditions at the end points will be discussed in Section~\ref{s:connectingorbitmap}.
If $\left( \bold{P}_{i} \right)_{i=1}^{m}$ admits a solution, then each $u_{i}$ is real-analytic since $g$ is. 
Therefore, there exists weights $\nu_{i} >1$ and \emph{real} coefficients $a^{i} \in \ell^{1}_{\nu_{i},n}$ such that 
\begin{align*}
	u_{i} = a^{i}_{0} + 2 \sum_{k=1}^{\infty} a^{i}_{k} T^{i}_{k}
\end{align*}
in $C\left(\left[t_{i-1},t_{i}\right]\right)$. Here $\left( T^{i}_{k} \right)_{k \in \NN_{0}}$ are the \emph{shifted} Chebyshev-polynomials on $\left[t_{i-1}, t_{i} \right]$
defined by 
\begin{align}
	T^{i}_{k}(t) = T_{k} \left( \frac{ 2t - t_{i} - t_{i-1} }{ t_{i} - t_{i-1} } \right), \quad k \in \NN_{0}.  \label{e:Trescaled}
\end{align} 

Next, define the map $c : \ell^{1}_{\nu_{i},n} \rightarrow \ell^{1}_{\nu_{i},n}$ by  
\begin{align}
	\label{eq:c(a)}
	c(a) := \sum_{\alpha \in \mathcal{A} } g_{\alpha} a^{\alpha}, \quad a^{\alpha} := \prod_{j=1}^{n} \left[a\right]_{j}^{\alpha_{j}},
\end{align}
where the latter product is understood to be the symmetric discrete convolution $\ast$. Then 
\begin{align*}
	g \left( u_{i} \right) = c_{0} \left( a^{i} \right) + 2 \sum_{k=1}^{\infty} c_{k} \left( a^{i} \right) T^{i}_{k},
\end{align*}
since the Chebyshev coefficients of the product of two functions is given by the symmetric discrete convolution.
Formally, we should incorporate the index $i$ into the notation for $c$ to emphasize its dependence on the weight $\nu_{i}$.
However, since the domain of $c$ can usually be easily inferred from the context, 
we haven chosen to omit the index from the notation. 
\begin{remark}
		\label{remark:Dc_comp}
		Throughout this paper we will need to analyze $Dc \left( a^{i} \right)$, where
		$a^{i} \in \ell^{1}_{\nu_{i},n}$ and $1 \leq i \leq m$, on numerous occasions. 
		For this reason, we state here for future reference how this derivative can be computed in an efficient way. 
		Since $\left( \ell^{1}_{\nu_{i},n}, \ast \right)$ is a Banach algebra, we may use the ``usual'' rules of calculus
		to compute the derivative of $c$. In particular, direct differentiation of \eqref{eq:c(a)} with respect to
		$a^{i}$ shows that 
		\begin{align}
				\label{eq:Dch}
				D\left[c\right]_{j} \left( a^{i} \right) \tilde a^{i} = \sum_{l=1}^{n} \hat g^{ijl} \ast \left[ \tilde a^{i} \right]_{l}, \quad
				\tilde a^{i} \in \ell^{1}_{\nu_{i},n}, \ 1 \leq j \leq n, 
		\end{align}
		where $\left [ c \right]_{j}$ denotes the $j$-th component of $c$ and 
		$\hat g^{ijl} \in \ell^{1}_{\nu_{i}}$ are the Chebyshev coefficients of 
		\begin{align}
				\label{eq:gijl}
				\frac{\partial g_{j} }{ \partial x_{l}} \left( a^{i}_{0} + 2 \sum_{k=1}^{\infty} a^{i}_{k} T^{i}_{k} \right), 
				\quad 1 \leq l \leq n.
		\end{align}
\end{remark}

Substitution of the Chebsyshev expansions for $\left(u_{i}\right)_{i=1}^{m}$ into $\left( \bold{P}_{i} \right)_{i=1}^{m}$ 
yields an equivalent system of equations for the coefficients and gives rise to the following map: 
\begin{definition}[Chebyshev map for ODEs] 
	\label{def:ODEmap}
	Let $\left( \nu_{i} \right)_{i=1}^{m}$ and $\left( \tilde \nu_{i} \right)_{i=1}^{m}$ be collections 
	of weights such that $1 < \tilde \nu_{i} < \nu_{i}$ for all $1 \leq i \leq m$. The Chebyshev map for ODEs
	is the function
	$F_{u}: \bigoplus_{i=1}^{m} \ell^{1}_{\nu_{i},n} \rightarrow \ell^{1}_{\tilde \nu_{1},n} / \CC^{n} \oplus
	\bigoplus_{i=2}^{m} \ell^{1}_{\tilde \nu_{i},n }$ 
	defined by
	\begin{align*}
		F_{u} \left( a \right) := 
		\left( f_{1} \left(a^{1} \right), f_{2} \left(a^{1}, a^{2} \right), \ldots, f_{m} \left(a^{m-1}, a^{m} \right) \right), 
	\end{align*}
	where $a = \left(a^{1}, \ldots, a^{m} \right)$, $f_{1}: \ell^{1}_{\nu_{1}, n } \rightarrow 
	\ell^{1}_{\tilde \nu_{1},n}/ \CC^{n}$
	is given by 
	\begin{align*}
	\left( f_{1} \left( a^{1} \right) \right)_{k} &:= 
		k a^{1}_{k} - \dfrac{ L\left(t_{1}-t_{0}\right) }{4} \left( c_{k-1} \left( a^{1} \right) - c_{k+1}\left( a^{1} \right) \right), 
		\quad k \in \NN,
	\end{align*}
	and $f_{i}: \ell^{1}_{\nu_{i-1},n} \times \ell^{1}_{\nu_{i},n} \rightarrow \ell^{1}_{\tilde \nu_{i},n}		$ by 
	\begin{align*}
	f_{i} \left(a^{i-1}, a^{i} \right) &:= 
		\begin{cases}
			\displaystyle
			a^{i}_{0} - a^{i-1}_{0} + 2 \sum_{l=1}^{\infty} \left( \left(-1\right)^{l} a^{i}_{l} - a^{i-1}_{l} \right), & k=0, \\[3ex]
			k a^{i}_{k} - \dfrac{ L \left( t_{i}-t_{i-1} \right) }{4} \left( c_{k-1} \left( a^{i} \right) - c_{k+1} \left( a^{i} \right) \right), & k \in \NN,
		\end{cases}
	\end{align*}
	for $2 \leq i \leq m$. 
\end{definition}
\begin{remark}
	The map $F_{u}$ is well-defined, since $\left( k a^{i}_{k} \right)_{k \in \NN_{0}} \in \ell^{1}_{\tilde \nu_{i},n}$ 
	for any $a^{i} \in \ell^{1}_{\nu_{i},n}$ and $1<\tilde \nu_{i} < \nu_{i}$.  
\end{remark}

Let us stress the subtle difference between the latter map and the one constructed in \cite{domaindecomposition};
in the current setting we allow for \emph{complex} Chebyshev coefficients. The main reason for this is that the parameterization maps
$P$ and $Q$ (see the previous section) are in principle complex-valued and we wish to use them to proof that 
$u(0) \in W^{u}_{\text{loc}} \left( \tilde p_{0} \right)$ and $u(L) \in W^{s}_{\text{loc}} \left( \tilde q_{0} \right)$. We will conclude
a-posteriori that the Chebyshev coefficients are in fact real by invoking symmetry arguments. 
Indeed, for any $a = \left(a^{1}, \ldots, a^{m} \right)$, define 
$\overline{a} := \left( \overline{a^{1}}, \cdots, \overline{a^{m}} \right)$ by 
$\left( \overline{ a^{i} } \right)_{k} := \overline{a^{i}_{k}}$. We will conclude that an element $a$ is real,
i.e., $\overline{a} = a$, by using the following observation: 
\begin{lemma}
	\label{lemma:symmetry_F_u}
	The map $F_{u}$ is compatible with conjugation, i.e., $F_{u} \left( \overline{a} \right) = \overline{ F_{u}(a) }$. 
	\begin{proof}
		Let $a = \left( a^{1}, \ldots, a^{m} \right) \in \prod_{j=1}^{m} \ell^{1}_{\nu_{j},n}$ be arbitrary. 
		We will prove that $c \left( \overline{a^{j}} \right) = \overline{ c \left( a^{j} \right)}$ for $1 \leq j \leq m$. 
		The desired result follows directly from this observation. 
		Recall that a Chebyshev series is a Fourier series up to coordinate transformation. To be more precise, 
		let $1 \leq j \leq m$, set $a^{j}_{-k} := a^{j}_{k}$ for $k \in \NN$, and define
		\begin{align*}
			s_{j} \left( \theta \right) := \frac{ t_{j} - t_{j-1} }{2} \left( \cos \theta + 1 \right) + t_{j-1}, \quad \theta \in \left[0, \pi \right]. 
		\end{align*}
		Then 
		\begin{align*}
			a^{j}_{0} + 2 \sum_{k=1}^{\infty} a^{j}_{k} T^{j}_{k} 
			\left( s_{j} \left( \theta \right) \right) =
			\sum_{k \in \ZZ} a^{j}_{k} e^{ik \theta}, \quad \theta \in \left[0, \pi \right],
		\end{align*}
		by~\eqref{e:Trescaled} and the definition of the Chebyshev polynomials (see Definition \ref{def:chebpoly}). 
		Note that the latter series converges uniformly to an analytic $2 \pi$-periodic function,
		since $a^{j} \in \ell^{1}_{\nu_{j},n}$. Furthermore, since $a^{j}_{-k} = a^{j}_{k}$ (by definition), it follows that 
		\begin{align*}
			g \left( a^{j}_{0} + 2 \sum_{k=1}^{\infty} a^{j}_{k} T^{j}_{k}  \right) = \sum_{k \in \ZZ} c_{k} \left( a^{j} \right) e^{ik\theta}, 
		\end{align*}
		where $c$ is defined in \eqref{eq:c(a)} and we have set $c_{-k} \left( a^{j} \right) := c_{k} \left( a^{j} \right)$ for $k \in \NN$. 
		In particular, 
		\begin{align*}
			\overline{ g \left( a^{j}_{0} + 2 \sum_{k=1}^{\infty} a^{j}_{k} T^{j}_{k} \right) } = 
			\sum_{k \in \ZZ} \overline{c_{k} \left( a^{j} \right)} e^{ik \theta}.
		\end{align*}  
		On the other hand, since $g$ is real-analytic, similar reasoning shows that 
		\begin{align*}
			\overline{ g \left( a^{j}_{0} + 2 \sum_{k=1}^{\infty} a^{j}_{k} T^{j}_{k} \right) }  = 
			g \left( \overline{a^{j}_{0}} + 2 \sum_{k=1}^{\infty} \overline{a^{j}_{k}} T^{j}_{k} \right) = 
			\sum_{k \in \ZZ} c_{k} \left( \overline{a^{j}} \right) e^{ik \theta}.
		\end{align*}
		Therefore, since a (pointwise) convergent Fourier series is unique, we conclude that $c \left( \overline{ a^{j} } \right) = \overline{c \left( a^{j} \right)}$. 
	\end{proof}
\end{lemma}

Finally, observe that by construction we now have the following result: 
\begin{proposition} 
	\label{prop:orbit}
	Suppose $a \in \bigoplus_{i=1}^{m} \ell^{1}_{\nu_{i},n}$ is symmetric, i.e., $\overline{a} = a$, then 
	$F(a)=0$ if and only if the functions  
	$\left \{ u_{i}  =  a^{i}_{0} + 2 \sum_{k=1}^{\infty} a^{i}_{k} T^{i}_{k} : 1 \leq i \right.$ $\left. \leq m \right \}$
	constitute a solution of $\left( \bold{P}_{i} \right)_{i=1}^{m}$.
\end{proposition}

\subsection{The connecting orbit map}
\label{s:connectingorbitmap}
In this section we set up a zero finding problem for \eqref{eq:bvp1}. We have already set up appropriate 
zero finding mappings for the ODE and charts on the (un)stable manifolds. What remains is imposing appropriate 
phase conditions. 

\paragraph{Boundary conditions}
We can now replace the boundary conditions in 
\eqref{eq:bvp1} with explicit equations. Let $P$ and $Q$ denote the local parameterizations of 
the complex unstable and stable manifold as before, respectively. Then the conditions $u(0) \in W^{u}_{\text{loc}} \left( \tilde p_{0} \right)$ 
and $u(1) \in W^{s}_{\text{loc}} \left( \tilde q_{0} \right)$ are equivalent to the problem of finding 
coordinates $\theta \in \mathbb{B}^{\text{sym}}_{\nu_{u}}$ and $\phi \in \mathbb{B}^{\text{sym}}_{\nu_{s}}$ such that 
\begin{align*}
	u_{1}(0) - P \left( \theta \right) 
	= a^{1}_{0} + 2 \sum_{k=1}^{\infty} \left(-1\right)^{k} a^{1}_{k}
	- \sum_{k \in \NN_{0}^{n_{u}}} p_{k} \theta^{k} &= 0, \\[1ex] 
	u_{m}(1) - Q \left( \phi \right) 
	= a^{m}_{0} + 2 \sum_{k=1}^{\infty} a^{m}_{k}
	- \sum_{k \in \NN_{0}^{n_{s}}} q_{k} \phi^{k} &= 0.
\end{align*}
As mentioned before, we will verify a-posteriori that $\theta^{\star} = \theta$ and $\phi^{\star} = \phi$
so that $P \left( \theta \right)$ and $Q \left( \theta \right)$
are points on the \emph{real} (un)stable manifolds. 

\paragraph{Length of the eigenvectors}
Recall that the first order Taylor coefficients of $P$ and $Q$ are determined
up to a rescaling. To get rid of this extra degree of freedom we prescribe the
length and orientation of the eigenvectors of $Dg \left( \tilde p_{0} \right)$ and
$Dg \left( \tilde q_{0} \right)$. More precisely, we require that 
\begin{align*}
	\left \langle p_{k}, \hat p_{k} \right \rangle- \epsilon_{u,k} = 0, \quad \left \vert k \right \vert &= 1, \ k \in \NN_{0}^{n_{u}}, \\[2ex]
	\left \langle q_{k}, \hat q_{k} \right \rangle- \epsilon_{s,k} = 0, \quad \left \vert k \right \vert &=1,  \ k \in \NN_{0}^{n_{s}},
\end{align*} 
where $\hat p_{k}, \hat q_{k} \in \CC^{n}$ are prescribed vectors and $\epsilon_{u,k}, \epsilon_{s,k}>0$.  
In practice, $\hat p_{k}$ and  $\hat q_{k}$ are numerical approximations of the eigenvectors of 
$Dg \left( \tilde p_{0} \right)$ and $Dg \left( \tilde q_{0} \right)$, respectively, and $\epsilon_{u,k}, \epsilon_{s,k}$
are their respective squared lengths. In order to respect the symmetry 
$\star$, we impose the following ordering:
\begin{align}
		\label{eq:eigvec_u_order}
		\hat p_{k^{\star}} = \overline{ \hat p_{k} }, \quad 
		\epsilon_{u,k^{\star}} = \epsilon_{u,k}, \quad \left \vert k \right \vert = 1, \  k \in \NN_{0}^{n_{u}}, \\[2ex]
		\label{eq:eigvec_s_order}
		\hat q_{k^{\star}} = \overline{ \hat q_{k} }, \quad 
		\epsilon_{s,k^{\star}} = \epsilon_{s,k}, \quad \left \vert k \right \vert = 1, \  k \in \NN_{0}^{n_{s}}.
\end{align}
We recall again that the length of the eigenvectors 
determines the decay rate of the Taylor coefficients and hence the size of the domains of $P$ and $Q$. 

\paragraph{Translation invariance in time}
Finally, we introduce a phase condition to fix the time parameterization of the connecting orbit. 
In \cite{MR635945,MR2359336} a phase condition specifically tailored for continuation is presented. The idea is to fix
a reference function $\tilde u$ and to minimize the functional 
\begin{align*}
	s \mapsto \int_{-\infty}^{\infty} \left \Vert u(t-s) - \tilde u (t) \right \Vert_{2}^{2} \ \mbox{d}t
\end{align*}
on some appropriate functions space, where $u$ is the connecting orbit.
This phase condition is used in popular software
packages for continuation such as \textsc{AUTO} and \textsc{Matcont} and was first suggested in \cite{MR635945}.

The intuition is that this phase condition enforces the connecting orbit to remain as close as possible 
(in the $L^{2}$-sense) to the reference solution 
with respect to small shifts in time. In practice, $\tilde u$ is the solution computed at the previous 
continuation step (or just the numerical approximation $\hat u$ in case we are not performing
continuation). In particular, a necessary condition for the latter functional to have a minimum at $s=0$ is
\begin{align}
	\label{eq:phase_Doedel}
	\int_{-\infty}^{\infty} \left \langle u(t) - \tilde u(t), \tilde u'(t)  \right \rangle \mbox{d}t = 0. 
\end{align}

We shall use \eqref{eq:phase_Doedel} to construct an appropriate phase condition in terms of Chebyshev
coefficients by approximating the integral on a finite domain. First, write 
\begin{align*}
	u = \sum_{i=1}^{m} u_{i} \bold{1}_{ \left[ t_{i-1}, t_{i} \right] }, \quad 
	\tilde u = \sum_{i=1}^{m} \tilde u_{i} \bold{1}_{ \left[ t_{i-1}, t_{i} \right] }.
\end{align*}
If the time of flight $L>0$ is sufficiently large, then 
\begin{align*}
	\int_{-\infty}^{\infty} \left \langle u(t) - \tilde u(t), \tilde u'(t) \right \rangle \mbox{d}t &\approx 
	\int_{-1}^{1} \left \langle u(t) - \tilde u(t), \tilde u'(t) \right \rangle \mbox{d}t \\[2ex] &= 
	\sum_{i=1}^{m} \int_{t_{i-1}}^{t_{i}} \left \langle u_{i}(t) - \tilde u_{i}(t), \tilde u_{i}'(t) \right \rangle \mbox{d}t. 
\end{align*}
Next, write 
\begin{align*}
	u_{i} = a^{i}_{0} + 2 \sum_{k=1}^{\infty} a^{i}_{k} T^{i}_{k}, \quad 
	\tilde u_{i} = b^{i}_{0} + 2 \sum_{k=1}^{\infty} b^{i}_{k} T^{i}_{k}, \quad 1 \leq i \leq m.
\end{align*}
For notational convenience, let us omit the superscripts from the Chebyshev coefficients 
and assume (for the moment) that $u_{i}$ and $\tilde u_{i}$ are scalar functions. Then 
\begin{align}
	& \left \langle u_{i} - \tilde u_{i}, \tilde u_{i}' \right \rangle_{L^{2}} \nonumber \\[2ex] & \quad =  
	2 \left( a_{0} - b_{0} \right) \sum_{l=1}^{\infty} \overline{b_{l}} 
	\left \langle T^{i}_{0}, \frac{ dT^{i}_{l} }{ dt} \right \rangle_{L^{2}} +
	4 \sum_{k,l=1}^{\infty} \left( a_{k} - b_{k} \right) \overline{b_{l}}
	\left \langle T^{i}_{k}, \frac{ dT^{i}_{l} }{ dt} \right \rangle_{L^{2}},
	\label{eq:phase_condition1}
\end{align}
where $\left \langle \cdot, \cdot \right \rangle_{L^{2}}$ denotes the standard 
complex inner product on $L^{2} \left( \left[t_{i-1}, t_{i} \right] \right)$.

Now, rescale time back to $[-1,1]$, use the coordinate transformation $\theta = \arccos \left( t \right)$ and the definition
of the Chebyshev polynomials to see that 
\begin{align*}
	\left \langle T^{i}_{k}, \frac{ d T^{i}_{l} }{ dt} \right \rangle_{L^{2}} &= 
	l \int_{0}^{\pi} \sin \left( l \theta \right) \cos \left( k \theta \right) d \theta \\[2ex] &= 
	\begin{cases}
		\dfrac{2l^{2}}{ l^{2} - k^{2}}, & k+l \equiv 1 \mod 2, \\[2ex]
		0, & \mbox{otherwise},
	\end{cases}
\end{align*}
for any $k, l \in \NN_{0}$. Finally, substitution of the latter expression into \eqref{eq:phase_condition1} yields
\begin{align}
	\label{eq:phase_condition2}
	\left \langle u_{i} - \tilde u_{i}, \tilde u_{i}' \right \rangle_{L^{2}} = 
	4 \left( \left( a_{0} - b_{0} \right) \sum_{l=0}^{\infty} \overline{b_{2l+1}} + 
	2 \sum_{s=1}^{\infty} \quad \ \sum_{ \mathclap{ \substack{ k+l = 2s+1 \\ k,l \in \NN}} } \  
	\left( a_{k} - b_{k} \right) \overline{b_{l}} \frac{l^{2}}{ l^{2} - k^{2} } \right).
\end{align}
If $u_{i}$ and $\tilde u_{i}$ are vector-valued, then we need to carry out the above computations
component-wise and sum over the components. 

In practice, we choose the Chebyshev coefficients of $\tilde u$ to be real, since in the end we wish
to establish the existence of a real-valued connecting orbit. Altogether, this motivates the following definition:
\begin{definition}[Phase condition for translation invariance in time]
	\label{def:time_parameterization}
	Let $b = \left( b^{1}, \ldots, \right.$ $\left. b^{m} \right) \in \bigoplus_{i=1}^{m} \ell^{1}_{\nu_{i},n}$ 
	be given symmetric sequences, i.e., $\bar{b} = b$, such that $b^{i}_{k} = 0$ for $k \geq N_{i}$, for some $N \in \NN^{m}$.
	The phase condition for translation invariance in time is the map 
	$\eta :  \bigoplus_{i=1}^{m} \ell^{1}_{\nu_{i},n} \rightarrow \CC$ defined by the following truncated version of \eqref{eq:phase_condition2}:
	\begin{align}
		& \eta \left( a^{1}, \ldots, a^{m} \right) := \nonumber \\[2ex] & \quad \sum_{i=1}^{m} \sum_{j=1}^{n}  
		\left( \left[ a^{i}_{0} - b^{i}_{0} \right]_{j} \sum_{k=0}^{\left \lfloor \frac{N_{i}-1}{2} \right \rfloor -1} \left[ b^{i}_{2k+1} \right]_{j} + 
		2 \sum_{s=1}^{N_{i}-2} \quad \ \sum_{ \mathclap{ \substack{ k+l = 2s+1 \\ 1 \leq k,l \leq N_{i}-1}} } \ \left[ a^{i}_{k} - b^{i}_{k} \right]_{j}
		\left[ b^{i}_{l} \right]_{j} \frac{l^{2}}{ l^{2} - k^{2} } \right).
		\label{eq:phase_cond}
	\end{align}
\end{definition}
\begin{remark}
	The expression for $\eta$ might seem complicated at first sight. Note, however, that $\eta$ is really just an affine
	linear map depending on finitely many components $a^i_k$ only. 
\end{remark}

We are now ready to set up the connecting orbit map. To this end, let 
\begin{align*}
	\nu := \left( \nu_{u}, \nu_{s}, \nu_{1}, \ldots, \nu_{m} \right), \quad 
	\tilde \nu := \left( \tilde \nu_{u}, \tilde \nu_{s}, \tilde \nu_{1}, \ldots, \tilde \nu_{m} \right)
\end{align*} be given weights such that $\tilde \nu < \nu$, $\nu_{u}, \nu_{s} >0$, $\nu_{i} >1$ for $1 \leq i \leq m$, and set 
\begin{align*}
	\mathcal{X}_{\nu} :&=  \mathbb{B}_{\nu_{u}} \times \mathbb{B}_{\nu_{s}} \times \CC^{n_{u}} \times \CC^{n_{s}} \times 
		\prod_{i=1}^{m} \ell^{1}_{\nu_{i},n} \times W^{1}_{\nu_{u},n, n_{u}} \times W^{1}_{\nu_{s},n, n_{s}}, \\[2ex]  
	\mathcal{Y}_{\tilde \nu} :&=  \CC^{2n+1+n_{u}+n_{s}} \times 
		\ell^{1}_{\tilde \nu_{1},n} / \CC^{n} \times \prod_{i=2}^{m} \ell^{1}_{\tilde \nu_{i},n} \times W^{1}_{\tilde \nu_{u},n, n_{u}} 
		\times W^{1}_{\tilde \nu_{s},n, n_{s}}.
\end{align*}
\begin{definition}[Chebyshev-Taylor map for connecting orbits]
	\label{def:connection_map}
	The Chebyshev-Taylor map $F: \mathcal{X}_{\nu} \rightarrow \mathcal{Y}_{\tilde \nu}$ for
	connecting orbits is defined by 
	\begin{align*}
		F(x) := 
		\begin{bmatrix}
			\displaystyle a^{1}_{0} + 2 \sum_{k=1}^{\infty} \left(-1\right)^{k} a^{1}_{k}
			- \sum_{k \in \NN_{0}^{n_{u}}} p_{k} \theta^{k}  \\[3ex]
			\displaystyle a^{m}_{0} + 2 \sum_{k=1}^{\infty} a^{m}_{k}
			- \sum_{k \in \NN_{0}^{n_{s}}} q_{k} \phi^{k} \\[4ex]
			\left[ \left \langle p_{e_{k}}, \hat p_{e_{k}} \right \rangle- \epsilon_{u,k} \right]_{k=1}^{n_{u}} \\[2ex]
			\left[ \left \langle q_{e_{k}}, \hat q_{e_{k}} \right \rangle- \epsilon_{s,k} \right]_{k=1}^{n_{s}} \\[2ex]
			\eta \left( a^{1}, \ldots, a^{m} \right), \\[2ex]
			F_{u} \left( a^{1}, \ldots, a^{m} \right) \\[2ex]
			F_{P} \left( \lambda^{u}, p \right) \\[2ex]
			F_{Q} \left( \lambda^{s}, q \right)
		\end{bmatrix},
	\end{align*}
	where $x = \left( \theta, \phi, \lambda^{u}, \lambda^{s}, a^{1}, \ldots, a^{m}, p, q \right)$. 
\end{definition}
\begin{remark}
	\label{remark:range_notation}
	We shall frequently denote elements in $\mathcal{Y}_{\tilde \nu}$ by 
	\begin{align*}
		y= \left( y_{t_{0}}, y_{t_{m}}, y_{\hat p_{1}}, y_{\hat q_{1}}, y_{\eta}, y_{a}, y_{p}, y_{q}\right),
	\end{align*}
	where 
	\begin{itemize}
	\item $y_{t_{0}}, y_{t_{m}} \in \CC^{n}$ correspond to the equations for the boundary conditions at 
		$t=t_{0}$ and $t=t_{m}$, respectively, 
	\item $y_{\hat p_{1}} \in \CC^{n_{u}}$, $y_{\hat q_{1}} \in \CC^{n_{s}}$ correspond to the equations for
		fixing the length and orientation of the eigenvectors, 
	\item $y_{\eta} \in \CC$ corresponds to the phase condition for fixing the time parameterization of the orbit,
	\item $y_{a},y_{p}$ and $y_{q}$ correspond to the Chebyshev and Taylor coefficients, respectively, as before. 
\end{itemize}
\end{remark}

The only reason for introducing the weights $\tilde \nu$ is to specify the codomain of $F$. 
These weights are irrelevant though, since we will establish the existence of a connecting orbit
by analyzing a fixed point map from $\mathcal{X}_{\nu}$ into itself, see Section \ref{sec:fixedpoint}. 
For this reason, we only specify a norm on $\mathcal{X}_{\nu}$.
Namely, we set 
\begin{align*}
	\left \Vert x \right \Vert_{\mathcal{X}_{\nu}} :=  
	 \max \left \{ \max_{1 \leq i \leq n_{u} }\left \vert \theta_{i} \right \vert, \right.& \left.
		  \max_{1 \leq i \leq n_{s} }\left \vert \phi_{i} \right \vert, \ 
		  \max_{1 \leq i \leq n_{u} }\left \vert \lambda^{u}_{i} \right \vert, \
		  \max_{1 \leq i \leq n_{s} }\left \vert \lambda^{s}_{i} \right \vert, \right. 
		  \\[2ex] &  \left.
		  \max_{1 \leq i \leq m} \left \Vert a^{i} \right \Vert_{\nu_{i},n},  \
	\left \Vert p \right \Vert_{\nu_{u},n}, \ \left \Vert q \right \Vert_{\nu_{s},n}  \right \}, 
\end{align*}
where $x = \left( \theta, \phi, \lambda^{u}, \lambda^{s}, a^{1}, \ldots, a^{m}, p, q \right)$.

\paragraph{Symmetry revisited}
Next, we examine the compatibility of $F$ with respect to the symmetries
introduced in the previous sections. In particular, the involution operations on the space of Taylor and Chebyshev
coefficients yield an symmetry operation 
$\boldsymbol{\star}$ on $\mathcal{X}_{\nu}$ defined by 
\begin{align*}
	x^{\boldsymbol{\star}} := \left( \theta^{\star}, \phi^{\star}, \left( \lambda^{u} \right)^{\star}, 
			\left( \lambda^{s} \right)^{\star}, \bar{a}, p^{\star}, q^{\star} \right).
\end{align*}
Similarly, we define an involution on the range $\mathcal{Y}_{\tilde \nu}$, also denoted by
$\boldsymbol{\star}$, via
\begin{align*}
	y^{\boldsymbol{\star}} := 
	\left( \overline{y_{t_{0}}},  \overline{y_{t_{m}}}, y^{\star}_{ \hat p_{1}}, 
	y^{\star}_{ \hat q_{1}}, 
	 \overline{y_{\eta}}, 
		\overline{y_{a}}, y_{P}^{\star}, y_{Q}^{\star}  \right).
\end{align*}
\begin{lemma}
	\label{lemma:symmetry}
	The map $F$ is compatible with $\boldsymbol{\star}$, i.e., 
	$F \left( x^{\boldsymbol{\star}} \right)  = F(x)^{\boldsymbol{\star}}$
	for any $x \in \mathcal{X}_{\nu}$. 
	\begin{proof}
			We start by considering the phase conditions associated to the unstable
			manifold. First, observe that 
			\begin{align*}
				\overline{
					a^{1}_{0} + 2 \sum_{k=1}^{\infty} \left(-1\right)^{k} a^{1}_{k} 
					 - \sum_{k \in \NN_{0}^{n_{u}}} p_{k} \theta^{k} } &= 
					\overline{a^{1}_{0}} + 2 \sum_{k=1}^{\infty} \left(-1\right)^{k} \overline{a^{1}_{k}} 
					- \sum_{k \in \NN_{0}^{n_{u}}} p^{\star}_{k} \left( \theta^{\star} \right)^{k}	
			\end{align*}			
			by definition of $\star$ and reordering of the series
			associated to the unstable manifold. 
			Furthermore, note that 
			\begin{align*}
					\left \langle p^{\star}_{k}, \hat p_{k} \right \rangle = 
					\overline{ \left \langle p_{k^{\star}}, \hat p_{k^{\star}} \right \rangle }, \quad
					\left \vert k \right \vert = 1, \quad k \in \NN_{0}^{n_{u}},
			\end{align*}
			since $\hat p_{k}$ was ordered in a symmetric way, see \eqref{eq:eigvec_u_order}. 
			The computations for the stable manifold are analogous. 				
			Next, observe that $\eta \left( \overline{a} \right) = \overline{\eta \left( a \right)}$, since
			the Chebyshev coefficients $b$ of the reference orbit $\tilde u$ are real. Finally, recall that 
			$F_{u} \left( \overline{a} \right) = \overline{F_{u} \left( a \right)}$,
			$F_{P} \left( \left(\lambda^{u}\right)^{\star}, p^{\star} \right) = F_{P} ( \lambda^{u}, p)^{\star}$ and 
			$F_{Q} \left(  \left(\lambda^{s}\right)^{\star}, q^{\star} \right) = F_{Q} ( \lambda^{s}, q)^{\star}$
			by Lemmas \ref{lemma:symmetry_F_u} and \ref{lemma:symmetry_F_Q}, respectively. 
			Altogether, this proves the result. 
	\end{proof}
\end{lemma}

We are now ready to formulate an appropriate characterization of a connecting orbit:
\begin{proposition}
	\label{prop:symmetry_connection}
	Suppose $F$ has a unique zero $x$ in some open neighborhood $U \subset \mathcal{X}_{\nu}$ 
	and assume that $x^{\boldsymbol{\star}} \in U$. Then $p_{0}, q_{0} \in \RR^{n}$ are equilibria of $g$, 
	$P \left( \theta \right) \in W^{u}_{\text{loc}} \left( p_{0} \right)$, 
	$Q \left( \phi \right) \in W^{s}_{\text{loc}} \left( q_{0} \right)$, and the map $u$ defined by the Chebyshev coefficients is an
	isolated connecting orbit from $p_{0}$ to $q_{0}$. 
	\begin{proof}
		Suppose $F(x) =0$, then the previous lemma implies that $x^{\boldsymbol{\star}} = x$, since
		$x$ is the only zero in $U$ and $x^{\boldsymbol{\star}} \in U$. Consequently, the Chebyshev
		coefficients $a$ are real and $P(\theta)$ and $Q( \phi)$ are points on the
		\emph{real} (un)stable manifolds by Proposition \ref{prop:real_manifold}. 
		Therefore, $u$~is a connecting orbit from $p_{0}$ to $q_{0}$ by Proposition \ref{prop:orbit}. Moreover, the connecting
		orbit is isolated, since $x$ is. 
	\end{proof}
\end{proposition}
\begin{remark}
		\label{remark:symmetry}
		In practice, we seek a zero of $F$ in a closed ball $B_{r} \left( \hat x \right)$ of radius $r>0$ centered at
		an approximate zero $\hat x$ obtained through numerical simulation.
	   	 The numerical computations yield an approximate zero which is almost symmetric (up to machine precision).
	    	We enforce that $ \left( \hat x \right)^{\star} = \hat x$ by going through ``all'' the 
		elements of $\hat x$ and imposing the exact symmetry conditions. For example, for the Taylor coefficients $\hat q$, 
		we determine all the multi-indices $k \in \mathcal{K}^{s}$ such that $k^{\star} \in \mathcal{K}^{s}$ and then 
		redefine $\hat q_{k^{\star}}$, for each $k^{\star} \not = k$, by setting it equal to $\overline{\hat q_{k}}$  (if $k^{\star} = k$ we set it equal to 
		$\normalfont \text{Re} \left( \hat q_{k} \right)$).  The symmetry implies that  
		$\left \Vert x - \hat x \right \Vert_{\mathcal{X}_{\nu}} = 
		\left \Vert x^{\star} - \hat x \right \Vert_{\mathcal{X}_{\nu}}$ for all $x \in \mathcal{X}_{\nu}$. Hence
		$B_{r} \left( \hat x \right)^{\star} = B_{r} \left( \hat x \right)$, which motivates the assumption that 
		$x^{\star} \in U$. 
\end{remark}

\paragraph{Transversality}
We end this section with a sufficient condition for proving that a connecting orbit is transverse.
The key observation is summarized in the following lemma:
\begin{lemma}
	\label{lemma:ker(DF(x))}
	Suppose $a, \tilde a \in \bigoplus_{i=1}^{m} \ell^{1}_{\nu_{i},n}$ are real.
	Let $u,w:[0,1] \rightarrow \RR^{n}$ denote the maps associated to $a$ and $\tilde a$, respectively, i.e., 
	\begin{align*}
			u :&= \sum_{i=1}^{m} \bold{1}_{ \left[ t_{i-1}, t_{i} \right] } u_{i}, \quad
			u_{i} :=  a^{i}_{0} + 2 \sum_{k=1}^{\infty} a^{i}_{k} T^{i}_{k}, \\[2ex]
			w :&= \sum_{i=1}^{m} \bold{1}_{ \left[ t_{i-1}, t_{i} \right] } w_{i}, \quad
			w_{i} :=  \tilde a^{i}_{0} + 2 \sum_{k=1}^{\infty} \tilde a^{i}_{k} T^{i}_{k}.
	\end{align*}
	Then $\tilde a \in \normalfont \text{ker} \left( DF_{u} \left( a \right) \right)$ if and only if $w'(t) = L Dg \left( u(t) \right) w(t)$ on $[0,1]$. 
	\begin{proof}
			A straightforward computations shows that $\tilde a \in \normalfont \text{ker} \left( DF_{u} \left( a \right) \right)$ if and only if 
			\begin{align*}
					\begin{cases}
							k \tilde a^{i}_{k} - \dfrac{L \left( t_{i}-t_{i-1}\right)}{4}
							\left( Dc_{k-1} \left( a^{i} \right) \tilde a^{i} - Dc_{k+1} \left( a^{i} \right) \tilde a^{i} \right)= 0, 
							& 1 \leq i \leq m, \ k \in \NN, \\[2ex]			
							\displaystyle
							\tilde a^{i}_{0} - \tilde a^{i-1}_{0} + \displaystyle 2 \sum_{l=1}^{\infty}
							\left( \left(-1\right)^{l} \tilde a^{i}_{l} - \tilde a^{i-1}_{l} \right) = 0, &
							2 \leq i \leq m.									
							\displaystyle 
					\end{cases}
			\end{align*}
			Furthermore, substitution of the expression in \eqref{eq:Dch} for $Dc \left( a^{i} \right)\tilde a^{i}$ shows that 
			the above system of equations is equivalent to 
			\begin{align*}
					\begin{cases}
						\dfrac{dw_{i}}{dt}(t) = L Dg \left( u_{i}(t) \right)w_{i}(t), & t \in \left[ t_{i-1}, t_{i} \right], \
						1 \leq i \leq m, \\[2ex]
						w_{i-1} \left( t_{i-1} \right) = w_{i} \left( t_{i-1} \right), & 2 \leq i \leq m,
					\end{cases}
			\end{align*}
			which proves the statement (see Proposition \ref{prop:orbit}).
	\end{proof}	
\end{lemma}

We are now ready to formulate a sufficient criterium for establishing the transversality of a connecting orbit. 
\begin{proposition}
		\label{prop:transversality}
		Suppose $x \in \mathcal{X}_{\nu}$ is symmetric and $F(x)=0$. If $DF(x)$ is injective,  then $x$ corresponds 
		to a transverse connecting orbit.
		\begin{proof}
				It is shown in Proposition \ref{prop:symmetry_connection} that $x$ corresponds to a connecting orbit $u$
				from $p_{0}$ to $q_{0}$ with the property that $u(0) = P \left( \theta \right) \in W^{u}_{\text{loc}} \left( p_{0} \right)$ 
				and $u(1) = Q \left( \phi \right) \in W^{s}_{\text{loc}} \left( q_{0} \right)$. 
				To show that $u$ is transverse, first observe that the mappings 
				\begin{align*}
					P \circ \iota_{u} : \mathbb{B}^{\text{sym},\text{re}}_{\nu_{u}} \subset \RR^{n_{u}} \rightarrow \RR^{n},
					\quad 
					Q \circ \iota_{s} : \mathbb{B}^{\text{sym},\text{re}}_{\nu_{s}} \subset \RR^{n_{s}} \rightarrow \RR^{n},
				\end{align*}
				are parameterizations of $W_{\text{loc}}^{u}\left( p_{0} \right)$ and
				$W_{\text{loc}}^{s}\left( q_{0} \right)$, respectively, by Proposition \ref{prop:real_manifold} and Remark~\ref{remark:real_manifold}.				
				Hence teh amp $\tilde \theta \mapsto \varphi \left( t, P \circ \iota_{u} \left( \tilde \theta \right) \right)$, where 
				$\varphi$ denotes the flow generated by~$Lg$, is a diffeomorphism from 
				$\mathbb{B}^{\text{sym},\text{re}}_{\nu_{u}}$ into
				$W^{u} \left( p_{0} \right)$ for any $t \in \RR$. Therefore, its derivative 
				$D_{x} \varphi \left(t, P \left( \theta \right) \right) D_{\theta} P \left(  \theta  \right) \iota_{u}$ (evaluated at $\tilde \theta = \iota_{u}^{-1} \left( \theta \right)$)
				is an isomorphism from $\RR^{n_{u}}$ onto $T_{ \varphi \left( t, P \left( \theta \right) \right)} W^{u} \left( p_{0} \right)$.
				Similarly, $D_{\phi} Q \left( \phi \right) \iota_{s}$ is an isomorphism from $\RR^{n_{s}}$ onto 
				$T_{ Q \left( \phi \right)} W^{s} \left( q_{0} \right)$. Consequently, since $\varphi \left(1, P \left( \theta \right) \right) = Q \left( \phi \right) = u(1)$, 
				the linear map 
				\begin{align*}
						\Phi_{1} := 
						\begin{bmatrix}
								D_{x} \varphi \left(1, P \left( \theta \right) \right) D_{\theta} P \left(  \theta  \right) \iota_{u} &
								- D_{\phi} Q \left( \phi \right) \iota_{s}
						\end{bmatrix}
				\end{align*}
				is a surjection from $\RR^{n_{u}} \times \RR^{n_{s}} = \RR^{n+1}$ onto 
				$\left( T_{ u(1) } W^{u} \left( p_{0} \right) + 
				T_{u(1)} W^{s} \left( q_{0} \right) \right) \subset \RR^{n}$. 
				
				Now, suppose $DF(x)$ is injective but $u$ is not transverse.
				Then the intersection of $W^{u} \left( p_{0} \right)$
				and $W^{s} \left ( q_{0} \right)$ is (in particular) not transverse
				at $u(1)$, since $u$ is transverse if and only if it is transverse at a point. Hence the map $\Phi_{1} : \RR^{n+1} \rightarrow \RR^{n}$
				cannot be surjective. Therefore, $\dim \left( \text{ker} \left( \Phi_{1} \right) \right) \geq 2$. Consequently, 
				there exist two linearly independent vectors 
				$
				\begin{bmatrix}
						\tilde \theta_{1} \\ \tilde \phi_{1}
				\end{bmatrix}, 
				\begin{bmatrix}
						\tilde \theta_{2} \\ \tilde \phi_{2}
				\end{bmatrix} \in \text{ker} \left( \Phi_{1} \right) \subset \RR^{n_{u}} \times \RR^{n_{s}}
				$.					
				We will show that this leads to a contradiction by constructing a nontrivial element in the kernel of $DF(x)$.
								
				Define $\xi_{1}, \xi_{2}: [0,1] \rightarrow \RR^{n}$ by 
				\begin{align*}
						\xi_{j}(t) := D_{x} \varphi \left(t, P \left( \theta \right) \right) D_{\theta} P \left( \theta \right) \iota_{u} \left( \tilde \theta_{j} \right),
						\quad j \in \{1,2\},
				\end{align*}
				then a straightforward computation shows that 
				\begin{align}
						\label{eq:trans_bvp}
						\begin{cases}
								\dfrac{d\xi_{j}}{dt}(t) = L Dg \left( u(t) \right) \xi_{j}(t), & t \in [0,1], \\[2ex]
								\xi_{j}(0) = D_{\theta} P \left( \theta \right) \iota_{u} \left( \tilde \theta_{j} \right),  \\[2ex]
								\xi_{j}(1) = D_{\phi} Q \left( \phi \right) \iota_{s} \left( \tilde \phi_{j} \right),
						\end{cases}
						\qquad j \in \{1,2\},
				\end{align}
				where the boundary condition at $t=1$ follows from the fact that 
				$
				\begin{bmatrix}
						\tilde \theta_{1} \\ \tilde \phi_{1}
				\end{bmatrix}, 
				\begin{bmatrix}
						\tilde \theta_{2} \\ \tilde \phi_{2}
				\end{bmatrix} \in \text{ker} \left( \Phi_{1} \right)
				$.					
				Further note that $\xi_{1}(t)$ and $\xi_{2}(t)$ are linearly independent for each $t \in [0,1]$, since the vectors 
				$
				\begin{bmatrix}
						\tilde \theta_{1} \\ \tilde \phi_{1}
				\end{bmatrix}$
				and
				$
				\begin{bmatrix}
						\tilde \theta_{2} \\ \tilde \phi_{2}
				\end{bmatrix} 
				$ 
				are, and the operators $D_{\theta} P \left( \theta \right) \iota_{u}$ and $D_{\phi} Q \left( \phi \right) \iota_{s}$ are injective. 	
				Consequently, since 
				\begin{align*}
					u = \sum_{i=1}^{m} \bold{1}_{ \left[ t_{i-1}, t_{i} \right]} \left( a^{i}_{0} + 2\sum_{k=1}^{\infty} a^{i}_{k} T^{i}_{k} \right), 
					\quad a^{i} \in \ell^{1}_{\nu_{i},n},
				\end{align*}
				where $a^{1}, \ldots, a^{m}$ are real, there exist (unique) real Chebyshev coefficients 
				$b_{1}, b_{2} \in \bigoplus_{i=1}^{m} \ell^{1}_{\nu_{i},n}$ such that
				\begin{align*}
						\xi_{j}  = \sum_{i=1}^{m}  \bold{1}_{ \left[ t_{i-1}, t_{i} \right]}
										\left( \left(b^{i}_{j} \right)_{0} + 2\sum_{k=1}^{\infty} \left(b^{i}_{j} \right)_{k} T^{i}_{k} \right), 
										\quad j \in \{1,2\}. 
				\end{align*}
				In particular, note that any linear combination of $b_{1}$ and $b_{2}$ corresponds to a solution of \eqref{eq:trans_bvp} and is thus
				an element in $\text{ker} \left( DF_{u} \left( a \right) \right)$ by Lemma \ref{lemma:ker(DF(x))}. 
				
				Now, set $\tilde \lambda^{u} := \bold{0}_{n_{u}}$,  $\tilde \lambda^{s} := \bold{0}_{n_{s}}$,
				$\tilde p := 0$, $\tilde q :=0$ and
				\begin{align*}
						h_{j} := \left( \iota_{u} \left( \tilde \theta_{j} \right), \iota_{s} \left( \tilde \phi_{j} \right), \tilde \lambda^{u}, \tilde \lambda^{s}, b_{j}, \tilde p, \tilde q \right)
						\in \mathcal{X}_{\nu},
						\quad j \in \{1,2\}.
				\end{align*} 
				If $D\eta(a) b_{j} =0$ for some $j \in \{1,2\}$, where $\eta$ is the phase condition defined
				in \eqref{eq:phase_cond}, then a straightforward computation shows that 
				$h_{j} \in \text{ker} \left( DF(x) \right)$.
				Otherwise, without loss of generality, we may assume that $D \eta (a) b_{1} \not=0$ and set 
				\begin{align*}
						h :=  h_{2} - \frac{ D\eta(a)b_{2} }{ D\eta(a)b_{1}} h_{1}. 
				\end{align*}
				A straightforward computation then shows that $0 \neq h \in \text{ker} \left( DF(x) \right)$. 
 Therefore, we have reached a contradiction, since $DF(x)$ 
				is assumed to be injective. 
				Hence $u$ must be transverse.  
				\end{proof}
\end{proposition}
\begin{remark}
	In practice, the injectivity of $DF(x)$ follows directly from our
	computer-assisted proof (a contraction argument) and is thus
	obtained for ``free'', see Remark \ref{remark:transversality}.
\end{remark}

\section{Functional analytic setup}
\label{sec:setup}

In this section we set up a functional analytic framework for establishing the existence of an isolated
zero of $F$. We start by introducing some notation and a finite dimensional reduction of $F$. 
We then combine numerical simulation and analysis on paper to set up a Newton-like operator $T$ whose fixed
points correspond to zeros of $F$. Finally, we derive a finite number of inequalities to establish that
$T$ is a contraction in a neighborhood of an approximate zero. 

\subsection{Projections}
In this section we define projections on both the range and domain. These projections
will help structure the calculations in the following sections. 

\paragraph{Projections in $\mathcal{X}_{\nu}$}
Write $x = \left( \theta, \phi, \lambda^{u}, \lambda^{s}, a^{1}, \ldots, a^{m}, p, q \right) \in \mathcal{X}_{\nu}$. 
Let $1 \leq i \leq m$, $1 \leq j \leq n$ and $\boldsymbol{k} \subset \NN_{0}$. Define projections 
$\Pi_{a}^{i} : \mathcal{X}_{\nu} \rightarrow \ell^{1}_{\nu_{i},n}$, \ 
$\Pi_{a}^{ij}, \ \Pi^{ij\boldsymbol{k}}_{a} : \mathcal{X}_{\nu} \rightarrow \ell^{1}_{\nu_{i}}$
onto the Chebyshev coefficients by
\begin{align*}
	\Pi^{i}_{a} \left( x \right) := a^{i}, \quad 
	\Pi^{ij}_{a} \left( x \right) := \left[ a^{i} \right]_{j}, \quad
	\left( \Pi^{ij\boldsymbol{k}}_{a} \left( x \right) \right)_{l} := 
	\begin{cases}
		\left[ a^{i}_{l} \right]_{j} & l \in \boldsymbol{k}, \\
		0 & \mbox{otherwise}.
	\end{cases}
\end{align*}
In particular, if $\boldsymbol{k}$ is a singleton, we identify $\Pi^{ij\boldsymbol{k}}_{a}(x) \simeq \left[ a^{i}_{\boldsymbol{k}} \right]_{j}$.
Similarly, let $\mathcal{I} \subset \NN_{0}^{n_{u}}$, $\mathcal{J} \subset \NN_{0}^{n_{s}}$ and define projections
$\Pi_{P} : \mathcal{X}_{\nu} \rightarrow W^{1}_{\nu_{u,n}}$, $\Pi^{j}_{P}, \ \Pi^{j \mathcal{I}}_{P} : \mathcal{X}_{\nu} \rightarrow W^{1}_{\nu_{u}}$  
and $\Pi_{Q} : \mathcal{X}_{\nu} \rightarrow W^{1}_{\nu_{s,n}}$, 
$\Pi^{j}_{Q},  \ \Pi^{j \mathcal{J}}_{Q} : \mathcal{X}_{\nu} \rightarrow W^{1}_{\nu_{s}}$
onto the Taylor coefficients of the (un)stable manifolds by 
\begin{align*}
	\Pi_{P} \left( x \right) :&= p, \quad \Pi^{j}_{P} \left( x \right) := \left [p\right]_{j}, \quad 
	\left( \Pi^{j \mathcal{I}}_{P} \left( x \right) \right)_{k} := 
	\begin{cases}
		\left[ p_{k} \right]_{j}, & k \in \mathcal{I}, \\
		0, & k \not \in \mathcal{I},
	\end{cases} \\[2ex]	
	\Pi_{Q} \left( x \right) :&= q, \quad \Pi^{j}_{Q} \left( x \right) := \left [q \right]_{j},  \quad
	\left( \Pi^{j \mathcal{J}}_{Q} \left( x \right) \right)_{k} := 
	\begin{cases}
		\left[ q_{k} \right]_{j}, & k \in \mathcal{J}, \\
		0, & k \not \in \mathcal{J}.
	\end{cases}	
\end{align*}
As before, if $\mathcal{I}$ and $\mathcal{J}$ are singletons, we identify 
$\left( \Pi^{j \mathcal{I}}_{P} \left( x \right) \right) \simeq \left[ p_{\mathcal{I}} \right]_{j}$,  
$\left( \Pi^{j \mathcal{J}}_{Q} \left( x \right) \right) \simeq \left[ q_{\mathcal{J}} \right]_{j}$. 
Finally, we define projections $\Pi^{j}_{\theta}, \Pi^{j}_{\phi}, \Pi^{j}_{\lambda^{u}}, \Pi^{j}_{\lambda^{s}} :
\mathcal{X}_{\nu} \rightarrow \CC$ by 
\begin{align*}
	\Pi^{j}_{\theta} \left( x \right) :&= \theta_{j}, \quad
	\Pi^{j}_{\lambda^{u}} \left( x \right) := \lambda^{u}_{j}, \quad 1 \leq j \leq n_{u}, \\[2ex]
	\Pi^{j}_{\phi} \left( x \right) :&= \phi_{j}, \quad 
	\Pi^{j}_{\lambda^{s}} \left( x \right) := \lambda^{s}_{j}, 
	\quad 1 \leq j \leq n_{s}.
\end{align*}
\begin{remark}
		In order to keep the notation and number of symbols to a minimum, we have used 
		the symbol $k$ as a ``dummy'' index which, depending on the context, 
		can be either an element in $\NN_{0}, \NN_{0}^{n_{u}}$ or $\NN_{0}^{n_{s}}$.
\end{remark}

We shall denote the collection of projections onto the components of $\mathcal{X}_{\nu}$ by $\mathbb{P}$, i.e., 
\begin{align*}
	\mathbb{P} := & \left \{ \Pi^{j}_{\theta} : 1 \leq j \leq n_{u} \right \} \cup \left \{ \Pi^{j}_{\phi} : 1 \leq j \leq n_{s} \right \}
			      \cup \left \{ \Pi^{j}_{\lambda^{u}} : 1 \leq j \leq n_{u} \right \} \\[2ex]
			      & \cup \left \{ \Pi^{j}_{\lambda^{s}} : 1 \leq j \leq n_{s} \right \}  
			      \cup \left \{ \Pi^{ij}_{a} : 1 \leq i \leq m, \ 1 \leq j \leq n\right \} \\[2ex]
			      & \cup \left \{ \Pi^{j}_{P} : 1 \leq j \leq n \right \} \cup \left \{ \Pi^{j}_{Q} : 1 \leq j \leq n \right \}.
\end{align*}
Observe that 
\begin{align*}
	\left \Vert x \right \Vert_{\mathcal{X}_{\nu}} = \max_{ \Pi \in \mathbb{P}} \left \Vert \Pi(x) \right \Vert_{\Pi \left( \mathcal{X}_{\nu} \right)},
\end{align*}
where $\left \Vert \cdot \right \Vert_{\Pi \left( \mathcal{X}_{\nu} \right)}$ denotes the norm on 
$\Pi \left( \mathcal{X}_{\nu} \right) \in \left \{ \CC, \ell^{1}_{\nu_{i}}, W^{1}_{\nu_{u}}, W^{1}_{\nu_{s}} \right \}$.
In the following we shall be a bit more sloppy in describing subsets of $\mathbb{P}$ by omitting the ranges
for the components of the projections. For instance, whenever we write $\left \{ \Pi_{\theta}^{j}, \Pi_{\phi}^{j}, \Pi^{j}_{\lambda^{u}}, \Pi^{j}_{\lambda_{s}} \right \}$, 
we mean to say that this set contains \emph{all} components associated to these projections, i.e., it contains
$\Pi^{j}_{\theta}$ for $1 \leq j \leq n_{u}$, $\Pi^{j}_{\phi}$ for $1 \leq j \leq n_{s}$, etc. We will adopt the same convention
for the projections into the range which are introduced at the end of this section. 

\paragraph{Galerkin projection into $\mathcal{X}_{\nu}$}
Let $N \in \NN^{m}$ be a given truncation parameter and define operators 
$\Pi^{i}_{N} :  \ell^{1}_{\nu_{i},n} \rightarrow \ell^{1}_{\nu_{i},n}$ by  
\begin{align*}
	\left( \Pi^{i}_{N} \left( a^{i} \right) \right)_{k}:= 
	\begin{cases}
		a^{i}_{k}, & 0 \leq k \leq N_{i} -1, \\[2ex]
		\bold{0}_{n}, & k \geq N_{i},
	\end{cases}, \qquad \ 1 \leq i \leq m. 
\end{align*}
Similarly, let $K = \left( K^{u}, K^{s} \right) \in \NN^{n_{u}} \times \NN^{n_{s}}$, set
\begin{align*}
	\mathcal{K}^{u} :&= \left \{ k \in \NN_{0}^{n_{u}}: k_{i} \leq K^{u}_{i}, \ 1 \leq i \leq n_{u} \right \}, \\[2ex]
	\mathcal{K}^{s} :&= \left \{ k \in \NN_{0}^{n_{s}}: k_{i} \leq K^{s}_{i}, \ 1 \leq i \leq n_{s} \right \},
\end{align*}
and define 
$\Pi_{K^{u}} : W^{1}_{\nu_{u},n} \rightarrow W^{1}_{\nu_{u},n}$, 
$\Pi_{K^{s}} : W^{1}_{\nu_{s},n} \rightarrow W^{1}_{\nu_{s},n}$ 
by 
\begin{align*}
	\left( \Pi_{K_{u}} \left( p \right) \right)_{k} :=
	\begin{cases}
		p_{k}, & k \in \mathcal{K}^{u},\\[2ex]
		\bold{0}_{n}, & k \not \in \mathcal{K}^{u},
	\end{cases} \qquad
	\left( \Pi_{K^{s}} \left( q \right) \right)_{k} :=
	\begin{cases}
		q_{k}, & k \in \mathcal{K}^{s}, \\[2ex]
		\bold{0}_{n}, & k \not \in \mathcal{K}^{s}, 
	\end{cases}
\end{align*}
respectively. Finally, define the \emph{Galerkin-projection} 
$\Pi_{\text{dom}}^{NK} : \mathcal{X}_{\nu} \rightarrow \mathcal{X}_{\nu}$ into the domain by 
\begin{align*}
	\Pi^{NK}_{\text{dom}}:= I_{2 \left( n_{u} + n_{s} \right)} \ \oplus \ \bigoplus_{i=1}^{m} \Pi^{i}_{N} \ \oplus 
			\Pi_{K^{u}} \ \oplus \ \Pi_{K^{s}},
\end{align*}
where $I_{2 \left( n_{u} + n_{s} \right)}$ is the identity on $\CC^{2 \left( n_{u} + n_{s} \right)}$, and set
$\mathcal{X}^{NK} := \Pi^{NK}_{\text{dom}} \left( \mathcal{X}_{\nu} \right)$.
\begin{remark}
	In practice, we choose an ordering on the set of multi-indices and identify 
	$\left[ \Pi_{K_{u}} \left( p \right) \right]_{j}$ with a column vector in 
	$\CC^{ \prod_{i=1}^{n_{u}} \left( K_{i}^{u} + 1 \right)}$ (and we do the same for the
	Taylor coefficients associated to the stable manifold). The Chebyshev coefficients
	$\Pi^{i}_{N} \left( a \right)$ are identified with a column vector in $\CC^{nN_{i}}$ as well. 
	Altogether, this yields an identification of $\normalfont \Pi^{NK}_{\text{dom}}(x)$ with a vector in $\CC^{\kappa}$,
	where
	\begin{align*}
		\kappa &  = n \left( \sum_{i=1}^{m} N_{i}  + \prod_{i=1}^{n_{u}} \left( K^{u}_{i} + 1 \right) + 
				\prod_{i=1}^{n_{s}} \left( K^{s}_{i} + 1 \right) \right)
				+ 2 \left( n_{u} + n_{s} \right) 
		\\[2ex] &= n \left( 2 + \sum_{i=1}^{m} N_{i} +  \prod_{i=1}^{n_{u}} \left( K^{u}_{i} + 1 \right) + 
				\prod_{i=1}^{n_{s}} \left( K^{s}_{i} + 1 \right) + 
				 \right)
				+ 2,
	\end{align*}
	since $n_{u} + n_{s} = n +1$. 
	In particular, $\mathcal{X}^{NK} \simeq \CC^{\kappa}$. 
\end{remark}

\paragraph{Projections in $\mathcal{Y}_{\tilde \nu}$} 
Recall that the range consists of elements of the form 
$y = \left( y_{t_{0}}, y_{t_{m}}, \right.$ 
$\left. y_{\hat p_{1}}, y_{\hat q_{1}}, y_{\eta},y_{a}, y_{p}, y_{q}\right)$,
see Remark \ref{remark:range_notation}. 
We shall abuse notation and denote the projections onto the Chebyshev and Taylor coefficients 
in the range in the same way as for the domain. Furthermore, we define projections 
$\Pi^{j}_{t_{0}}, \Pi^{j}_{t_{m}}, \Pi^{j}_{\hat p_{1}}, \Pi^{j}_{\hat q_{1}}, \Pi_{\eta}: \mathcal{Y}_{\tilde \nu} \rightarrow \CC$ by 
\begin{align*}
	\Pi^{j}_{t_{0}} (y) :&= y_{t_{0}}, \quad \Pi^{j}_{t_{m}} (y) := y_{t_{m}}, \quad 1 \leq j \leq n, \\[2ex]
	\Pi^{j}_{\hat p_{1}} \left( y \right) :&= \left[ y_{\hat p_{1}} \right]_{j}, \quad 1 \leq j \leq n_{u}, \\[2ex]
	\Pi^{j}_{\hat q_{1}} \left( y \right) :&= \left[ y_{\hat q_{1}} \right]_{j}, \quad 1 \leq j \leq n_{s},  \\[2ex]
	\Pi_{\eta} \left( y \right) :&= y_{\eta}.
\end{align*}

The truncation operators on the range are defined in the same way as for the domain. For this reason 
we shall use the same notation to denote them. 
There is one slight modification in the projection onto the Chebyshev coefficients associated to the first domain, 
however, namely we set
\begin{align*}
	\left( \tilde \Pi^{1}_{N} \left( a^{1} \right) \right)_{k} := 
	\begin{cases}
		a^{1}_{k}, & 1 \leq k \leq N_{1}-1, \\[1ex]
		\bold{0}, & k \geq N_{1}.
	\end{cases}
\end{align*}
Finally, we define the Galerkin-projection $\Pi^{NK}_{\text{ran}}$ into the range by 
\begin{align*}
	\Pi^{NK}_{\text{ran}} := I_{\CC^{3n+2}} \oplus \tilde \Pi^{1}_{N} \oplus \bigoplus_{i=2}^{m} \Pi^{i}_{N} \oplus \Pi_{K^{u}} \oplus \Pi_{K^{s}}
\end{align*}	
and set $\mathcal{Y}^{NK} := \Pi^{NK}_{\text{ran}} \left( \mathcal{Y}_{\tilde \nu} \right)$. 
\begin{remark}
	Observe that $\mathcal{Y}^{NK} \simeq \CC^{\kappa}$, 
	since $n_{u} + n_{s} = n+1$. Hence $\mathcal{Y}^{NK}$ and $\mathcal{X}^{NK}$
	have the same dimension. 
\end{remark}

\subsection{An equivalent fixed-point problem}
\label{sec:fixedpoint}
In this section we construct a Newton-like operator and set up an equivalent fixed point problem. 
We start by introducing a finite dimensional reduction of the zero finding problem amenable to 
numerical computations: 
\begin{definition}[Finite dimensional reduction]
	The finite dimensional reduction $F_{NK} : \mathcal{X}^{NK} \rightarrow \mathcal{Y}^{NK}$
	of the connecting orbit map $F$ is defined by 
	\begin{align*}
		F_{NK} := \Pi_{\text{ran}}^{NK} \circ F \vert_{\mathcal{X}^{NK}}. 
	\end{align*}
\end{definition}

Next, we construct an approximation of $DF \left( \hat x \right)$ and its inverse by combining numerical computations
and analysis on paper. To this end, assume that we have computed 
\begin{itemize}
	\item[(\bf{A1})] an approximate zero 
			$\hat x = \left( \hat \theta, \hat \phi, \hat \lambda^{u}, \hat \lambda^{s}, \hat a, \hat p, \hat q \right) \in \mathcal{X}^{NK}$ 
			of $F_{NK}$ such that $\left( \hat x \right)^{\star} = \hat x$,
	\item[(\bf{A2})] an approximate \emph{injective} inverse $A_{NK}$ of $DF_{NK} \left( \hat x \right)$. 
\end{itemize}
If the truncation parameters are sufficiently large and the grid is sufficiently fine, we expect the linear part
of the mappings $F_{u}$, $F_{P}$ and $F_{Q}$ to be dominant in a small neighborhood of the approximate zero.
This motivates the following definitions: 
\begin{definition}[Approximate derivative]
	The approximate derivative $\widehat{DF}: \mathcal{X}_{\nu} \rightarrow \mathcal{Y}_{\tilde \nu}$ at 
	$\hat x$ is defined by 
	\begin{align*}
		\Pi \widehat{DF}h :&= \Pi (y),
		\quad \text{for} \ \Pi \in \left \{ \Pi^{j}_{t_{0}}, \Pi^{j}_{t_{m}}, \Pi^{j}_{\hat p_{1}}, \Pi^{j}_{\hat q_{1}}, \Pi_{\eta} \right \}, \\[2ex]
		\left[\Pi^{i}_{a} \widehat{DF}h \right]_{k} :&= 
		\begin{cases}
			\left( \Pi_{a}^{i} (y) \right)_{k}, & \delta_{i1} \leq k \leq N_{i}-1,
			\\[2ex]
			k \left( \Pi_{a}^{i} (h) \right)_{k}, & k \geq N_{i}, 
		\end{cases}, \qquad 1 \leq i \leq m, \\[2ex]
		\left[\Pi_{P} \widehat{DF}h \right]_{k} :&= 
		\begin{cases}
			\left( \Pi_{P} (y) \right)_{k}, & 
			k \in \mathcal{K}^{u}, \\[2ex]
			\left \langle \hat \lambda^{u},k \right \rangle \left( \Pi_{P} (h) \right)_{k}, & 
			k \not \in \mathcal{K}^{u},
		\end{cases} \\[2ex]
		\left[\Pi_{Q} \widehat{DF} h \right]_{k} :&= 
		\begin{cases}
			\left( \Pi_{Q} (y) \right)_{k}, & 
			k \in \mathcal{K}^{s}, \\[2ex]
			\left \langle \hat \lambda^{s},k \right \rangle \left( \Pi_{Q} (h) \right)_{k}, & 
			k \not \in \mathcal{K}^{s},
		\end{cases}
	\end{align*}
	where $y = DF_{NK} \left( \hat x \right) \Pi^{NK}_{\text{dom}}(h)$. 
\end{definition}
\begin{definition}[Approximate inverse]
	\label{def:approximate_inverse}
	The approximate inverse $A: \mathcal{Y}_{\nu} \rightarrow \mathcal{X}_{\nu}$ of $DF \left( \hat x \right)$ is defined by 
	\begin{align*}
		\Pi A h :&= \Pi (x), \quad \text{for} \
		\Pi \in \left \{ \Pi^{j}_{\theta}, \Pi^{j}_{\phi}, \Pi^{j}_{\lambda^{u}}, \Pi^{j}_{\lambda^{s}} \right \}, \\[2ex]
		\left( \Pi^{i}_{a}Ah \right)_{k} :&= 
		\begin{cases}
			\left( \Pi^{i}_{a} (x) \right)_{k}, & 0 \leq k \leq N_{i}-1, \\[2ex]
			k^{-1} \left( \Pi^{i}_{a}(h) \right)_{k}, & k \geq N_{i}, 
		\end{cases} \\[2ex]
		\left( \Pi_{P}Ah \right)_{k} :&= 
		\begin{cases}
			\left( \Pi_{P} (x) \right)_{k}, & k \in \mathcal{K}^{u}, \\[2ex]
			\left \langle \hat \lambda^{u},k \right \rangle^{-1} \left( \Pi_{P}(h) \right)_{k}, &
			k \not \in \mathcal{K}^{u},
		\end{cases} \\[2ex]
		\left( \Pi_{Q}Ah \right)_{k} :&= 
		\begin{cases}
			\left( \Pi_{Q} (x) \right)_{k}, & k \in \mathcal{K}^{s}, \\[2ex]
			\left \langle \hat \lambda^{s},k \right \rangle^{-1} \left( \Pi_{Q}(h) \right)_{k}, &
			k \not \in \mathcal{K}^{s},
		\end{cases}
	\end{align*}
	where $x = A_{NK}\Pi_{\text{ran}}^{NK}(h)$. 
\end{definition}
\begin{remark}
	Note that $A$ is injective, since $A_{NK}$ is (by assumption \textup{(\textbf{A2})}). 
\end{remark}

We are now ready to construct a Newton-like operator for $F$ based at the approximate zero: 
\begin{definition}[Newton-like operator]
	The Newton-like operator $T$ for $F$ based at $\hat x$ is defined by $T:=I-AF$. 
\end{definition}
A straightforward computation shows that $T$ maps $\mathcal{X}_{\nu}$ into itself by construction of the approximate inverse $A$. 
The weights $\tilde \nu$ are therefore irrelevant. Furthermore, observe that $T(x)=x$ if and only if $F(x)=0$, 
since $A$ is injective. We conclude this section with a theorem which can be used to prove that $T$ is a contraction 
in a neighborhood of $\hat x$ by checking a \emph{finite} number of inequalities. 
The theorem is based on a parameterized Newton-Kantorovich method and is  
often referred to as the radii-polynomial approach (see \cite{MR3454370} for instance).
\begin{theorem}[Contraction mapping principle with variable radius]
	\label{thm:contraction}
	Suppose for each $\Pi \in \mathbb{P}$ there exist bounds $Y_{\Pi}, Z_{\Pi}(r)>0$ such that 
		\begin{align}
		\label{eq:Y}
		\left \Vert \Pi \left( T \left( \hat x \right) - \hat x \right) \right \Vert_{\Pi \left( \mathcal{X}_{\nu} \right)} &\leq Y_{\Pi}, \\[2ex]
		\label{eq:Z}
		\sup_{v,h \in B_{1} \left( 0 \right)} 
		\left \Vert \Pi DT \left( \hat x + rv \right)h \right \Vert_{\Pi \left( \mathcal{X}_{\nu} \right)} &\leq Z_{\Pi}(r),
	\end{align}	
	where $\left \Vert \cdot \right \Vert_{\Pi \left( \mathcal{X}_{\nu} \right)}$ denotes the norm on $\Pi \left( \mathcal{X}_{\nu} \right)$.
	If there exists a radius $\hat r>0$ such that 
	\begin{align}
		\label{eq:radiipoly}
		Z_{\Pi} \left( \hat r \right) \hat r + Y_{\Pi} < \hat r
	\end{align}
	for all $\Pi \in \mathbb{P}$, then $T : B_{\hat r} \left( \hat x \right) \rightarrow B_{\hat r} \left( \hat x \right)$ is a contraction. 
	\begin{proof}
			A proof can be found in \cite{Yamamoto}. 
	\end{proof}
\end{theorem}
\begin{remark}
	\label{remark:transversality}
	If $T$ is a contraction on $B_{\hat r} \left( \hat x \right)$, then there exists a unique zero 
	$\tilde x \in B_{\hat r} \left( \hat x \right)$ of $F$. In particular, $\left( \tilde x \right)^{ \boldsymbol{\star} } \in B_{\hat r} \left( \hat x \right) $, 
	since  $\left( \hat x \right)^{\boldsymbol{\star}} = \hat x$ and $\boldsymbol{\star}$ is norm-preserving. 
	Therefore, since $\tilde x$ is unique in $B_{\hat r} \left( \hat x \right)$, we haven proven the existence of 
	a real connecting orbit by Proposition \ref{prop:symmetry_connection}. Moreover, $DF \left( \tilde x \right)$ is injective, since
	\begin{align*}
		\left \Vert I - ADF \left( \tilde x \right) \right \Vert_{\mathcal{B} \left( \mathcal{X}_{\nu}, \mathcal{X}_{\nu} \right)}  = 
		\left \Vert DT \left( \tilde x \right) \right \Vert_{\mathcal{B} \left( \mathcal{X}_{\nu}, \mathcal{X}_{\nu} \right)} <1
	\end{align*}
	by \eqref{eq:Z} and \eqref{eq:radiipoly}. 
	Hence the connecting orbit is transverse by Proposition 
\ref{prop:transversality}. 
\end{remark}
\begin{remark}
The radius ${\hat r}$ in Theorem~\ref{thm:contraction} serves as an error bound
on the solution. Indeed, since for any $\nu \geq 1$ the $\ell^1_\nu$ norm
controls the $C^0$ norm, it follows that along the part of the heteroclinic
orbit between the local invariant manifolds, which is described by the (domain
decomposed) Chebyshev series, the distance in phase space between the numerical
approximation and the solution is bounded by ${\hat r}$. An analogous bound holds
for the distance between the parts of the orbit lying near the (un)stable
manifolds: the distance between the orbit and the numerical approximation of
the local manifold is bounded by ${\hat r}$. More precisely locating the position of
the orbit within the local stable manifold (and analogously for the
unstable one) requires solving the corresponding linear flow with initial data ${\hat \phi}$,
see~\eqref{eq:linear_flow}, and considering its image under the mapping
$\hat{Q}: \phi \to \sum_{k \in \mathcal{K}^s} \hat{q}_k \phi^k$.
\end{remark}

\section{Bounds for proving contraction}
\label{sec:bounds}
In this section we compute the bounds as stated in Theorem \ref{thm:contraction} to prove that
$T$ is a contraction in a neighborhood of $\hat x$. To compute these bounds, we need to project 
and perform analysis on the various subspaces of $\mathcal{X}_{\nu}$. Since the analysis  
for the unstable and stable manifold is the same, we will only write down the arguments 
in detail for the unstable manifold and simply state the analogous result for the stable manifold. 
We have aimed to compute the sharpest bounds whenever possible, but there are occasions in 
which we have chosen to use slightly less optimal bounds when the reduction in computational complexity
outweighed the potential loss in accuracy.

\subsection{$Y$-bounds}
\label{sec:Ybounds}
In this section we compute bounds for the residual 
\begin{align*}
	T\left( \hat x \right) - \hat x = -AF\left( \hat x \right)
\end{align*}
as stated in Theorem \ref{thm:contraction}. 
To this end, observe that  
\begin{align*}
	\Pi AF \left( \hat x \right) = \Pi A_{NK} F_{NK} \left( \hat x \right), \quad 
	\Pi \in \left \{ \Pi_{\theta}^{j}, \Pi_{\phi}^{j}, \Pi^{j}_{\lambda^{u}}, \Pi^{j}_{\lambda_{s}} \right \}.
\end{align*}
Furthermore, for $k \in \NN_{0}$, $1 \leq i \leq m$, $1 \leq j \leq n$, we have that 
\begin{align*}
	& \Pi_{a}^{ijk} A F \left( \hat x \right)
	\\[2ex] & \quad = 
		\begin{cases}
			\Pi^{ijk}_{a} A_{NK} F_{NK} \left( \hat x \right), & 0 \leq k \leq N_{i}-1, \\[2ex]
			-\dfrac{L \left( t_{i} - t_{i-1} \right)}{4k} \left[ c_{k-1} \left( \hat a^{i} \right) - c_{k+1} \left( \hat a^{i} \right) \right]_{j}, &
			N_{i} \leq k \leq N_{g_{j}} \left( N_{i}-1 \right) +1, \\[2ex]
			0, & \text{otherwise},
		\end{cases}
\end{align*}
since $\left[ c_{k} \left( \hat a \right) \right]_{j} = 0$ for all $k \geq N_{g_{j}} \left(N_{i}-1\right)+1$, see
the definition of $c$ in \eqref{eq:c(a)}. Similarly, for $k \in \NN_{0}^{n_{u}}$, 
\begin{align*}
	& \Pi_{P}^{jk} A F \left( \hat x \right) = 
	\begin{cases}
		\Pi_{P}^{jk} A_{NK} F_{NK} \left( \hat x \right), & k \in \mathcal{K}^{u}, \\[2ex]
		-\left \langle k, \hat \lambda^{u} \right \rangle^{-1}
		\left[ C_{k} \left( \hat p \right) \right]_{j}, & k \in \mathcal{J}^{j}_{u}, \\[2ex]
		0, & \mbox{otherwise},
	\end{cases}
\end{align*}
where
\begin{align*}
	\mathcal{J}^{j}_{u} :&= 
	\left \{ k \in \NN_{0}^{n_{u}} :  k_{i} \leq N_{g_{j}} K_{i}^{u} \ , 1 \leq i \leq n_{u} \right \} \cap \left( \mathcal{K}^{u} \right)^{c},
\end{align*}
since $\left[ C_{k} \left( \hat p \right) \right]_{j} =0$ for $k \not \in \mathcal{K}^{u}$, see \eqref{eq:C(a)}. 

The above computations show that there are only a \emph{finite} number of non-vanishing terms in $AF \left( \hat x \right)$.
Therefore, $AF \left( \hat x \right)$ can be computed with the aid of a computer. 
It is now a straightforward task to compute the $Y$-bounds by taking the appropriate norms of the above expressions: 
\begin{proposition}[$Y$-bounds]
	Let $1 \leq i \leq m$ and $1 \leq j \leq n$. The bounds 
	\begin{align*}
		Y_{\Pi} :& = \left \vert \Pi A_{NK} F_{NK} \left( \hat x \right) \right \vert, \quad \Pi \in \left \{ \Pi_{\theta}^{j}, \Pi_{\phi}^{j}, \Pi^{j}_{\lambda^{u}}, \Pi^{j}_{\lambda_{s}} \right \}, \\[2ex]
		Y^{ij}_{a} :&= \left \Vert \Pi^{ij}_{a} A_{NK} F_{NK} \left( \hat x \right) \right \Vert_{\nu_{i}} \\[1ex] & \quad + 
		\frac{L \left( t_{i} - t_{i-1} \right) }{2} \sum_{k=N_{i}}^{ N_{g_{j}} \left( N_{i} - 1 \right) + 1 }
		\left \vert \left[ c_{k-1} \left( \hat a^{i} \right) - c_{k+1} \left( \hat a^{i} \right) \right]_{j} \right \vert 
		\frac{\nu_{i}^{k}}{k}, \\[2ex]
		Y^{j}_{P} :&= \left \Vert \Pi^{j}_{P} A_{NK} F_{NK} \left( \hat x \right) \right \Vert_{\nu_{u}} 
		+ \sum_{k \in \mathcal{J}^{j}_{u}} \left \vert \left[ C_{k} \left( \hat p \right) \right]_{j} \right \vert 
		\left \vert \left \langle k, \hat \lambda^{u} \right \rangle \right \vert^{-1} \nu_{u}^{ \left \vert k \right \vert } , \\[2ex]
		Y^{j}_{Q} :&= \left \Vert \Pi^{j}_{Q} A_{NK} F_{NK} \left( \hat x \right) \right \Vert_{\nu_{s}} 
		+ \sum_{k \in \mathcal{J}^{j}_{s}} \left \vert \left[ C_{k} \left( \hat q \right) \right]_{j} \right \vert 
		\left \vert \left \langle k, \hat \lambda^{s} \right \rangle \right \vert^{-1} \nu_{s}^{ \left \vert k \right \vert }
	\end{align*}
	satisfy the estimate in \eqref{eq:Y}. 
\end{proposition}

\subsection{Z-bounds}
\label{sec:Z}
In this section we compute bounds for $DT$ as stated in Theorem \ref{thm:contraction}. 
To this end, let $r>0$, $v,h \in B_{1}(0)$ be arbitrary and observe that 
\begin{align*}
	DT \left( \hat x + rv \right)h = 
	\left( I - A \widehat{DF} \right)h - 
	A \left( DF \left( \hat x + rv \right) - DF \left( \hat x \right) + DF \left( \hat x \right) - \widehat{DF} \right)h. 
\end{align*}
We shall use this decomposition to compute \emph{quadratic} polynomials $Z_{\Pi}(r)$
which satisfy the condition in \eqref{eq:Z}.
Furthermore, throughout this section we shall write 
\begin{align*}
	h = \left( \tilde \theta, \tilde \phi, \tilde \lambda^{u}, \tilde \lambda^{s}, \tilde a, \tilde p, \tilde q\right).
\end{align*}

Let us start with analyzing the easiest term which measures the quality of the approximate derivative and inverse: 
\begin{lemma}
	\label{lemma:H}
	Let $\Pi \in \mathbb{P}$, then 
	\begin{align*}
		\left \Vert \Pi \left( I - A \widehat{DF} \right) \right \Vert_{
		\mathcal{B} \left( \mathcal{X}_{\nu}, \Pi \left( \mathcal{X}_{\nu} \right) \right)} \leq 
		\left \Vert \Pi \left( I_{NK} - A_{NK} DF_{NK} \left( \hat x \right) \right) \right \Vert_{
		\mathcal{B} \left( \mathcal{X}^{NK}_{\nu}, \Pi \left( \mathcal{X}^{NK}_{\nu} \right) \right)}.
	\end{align*}
	\begin{proof}
		It suffices to observe that
		\begin{align*}
			 \Bigl( I - A \widehat{DF} \Bigr) h = 
			\Bigl( I_{NK} - A_{NK} DF_{NK} ( \hat x )\Bigr)\Pi_{NK}h , 
		\end{align*}
		where $I_{NK}$ is the identity on $\mathcal{X}^{NK}_{\nu}$. The latter equality holds because 
		$A$ and $\widehat{DF}$  are exact inverses of each other
		on the subspaces associated to the ``tails'' in $W^{1}_{\nu_{u},n}$, $W^{1}_{\nu_{s},n}$ and 
		$\bigoplus_{i=1}^{m} \ell^{1}_{\nu_{i},n}$. 
	\end{proof}
	\begin{remark}
			Note that the computation of the stated bound is finite for each $\Pi \in \mathbb{P}$, since 
			$I_{NK} - A_{NK} DF_{NK} \left( \hat x \right)$ is a finite dimensional matrix. 
	\end{remark}
\end{lemma}

\subsubsection{Chebyshev series: convolution terms}
\label{sec:conv_terms}
In this section we develop tools for analyzing the terms 
\begin{align*}
	\Pi_{a}^{ij\NN} \left( DF \left( \hat x \right) - \widehat{DF} \right), \quad 1 \leq i \leq m,  \ 1 \leq j \leq n,
\end{align*}
which will be used extensively in Section \ref{sec:firstorder} to compute the $Z$-bounds. 
We start with the observation that 
\begin{align}
	\label{eq:Pi_cheb0}
	& \left( \Pi_{a}^{i} \left( DF \left( \hat x \right) - \widehat{DF} \right)h \right)_{k} 
	\nonumber \\[2ex] & \quad = 
	- \dfrac{L \left( t_{i} - t_{i-1} \right) }{4} 
	\begin{cases}
		Dc_{k-1} \left( \hat a^{i} \right) \tilde a^{i}_{\infty} - Dc_{k+1} \left( \hat a^{i} \right) \tilde a^{i}_{\infty}, & 1 \leq k \leq N_{i}-1, \\[2ex]
		Dc_{k-1} \left( \hat a^{i} \right) \tilde a^{i}- Dc_{k+1} \left( \hat a^{i} \right) \tilde a^{i}, & k \geq N_{i}, 
	\end{cases}
\end{align}
for $1 \leq i \leq m$, where 
\begin{align*}
	\left( \tilde a^{i}_{\infty} \right)_{k} := 
	\begin{cases}
		\bold{0}_{n}, & 0 \leq k \leq N_{i}-1, \\[2ex]
		\tilde a^{i}_{k}, & k \geq N_{i}. 
	\end{cases}
\end{align*}
The goal is to construct workable matrix representations for both the finite (truncated) and tail part of \eqref{eq:Pi_cheb0}.

To construct suitable matrix representations for \eqref{eq:Pi_cheb0}, recall that 
\begin{align*}
		D\left[c\right]_{j} \left( \hat a^{i} \right) \tilde a^{i} = \sum_{l=1}^{n} \hat g^{ijl} \ast \left[ \tilde a^{i} \right]_{l}, \quad
		1 \leq i \leq m, \ 1 \leq j \leq n,
\end{align*}
where the coefficients $\hat g^{ijl}$ are defined in \eqref{eq:gijl}, see Remark \ref{remark:Dc_comp}. 
In particular, note that $\hat g^{ijl}_{k} =0 $ for $k \geq M_{jl} \left( N_{i} -1 \right) +1$, where 
$M_{jl} := \text{order} \left( \dfrac{ \partial g_{j}} {\partial x_{l} } \right)$, since
$\hat a^{i}_{k} = 0$ for $k \geq N_{i}$.
Therefore, motivated by the above observations,
we consider a sequence $a \in \ell^{1}_{\nu}$ such that $a_{k} =0$
for $k \geq \tilde M:= M \left( N-1 \right) +1$, where $M \in \NN$, $N \in \NN$, $\nu >1$, and 
construct \emph{explicit} matrix representations for the mappings
$B(a), \Gamma(a) : \ell^{1}_{\nu} \rightarrow \ell^{1}_{\nu}$ defined by 
\begin{align*}
		\left[ B(a) \tilde a \right]_{k} :&= 
		\begin{cases}
				\left( a \ast \tilde a_{\infty} \right)_{k-1} - \left( a \ast \tilde a_{\infty} \right)_{k+1}, & 1 \leq k \leq N-1, \\[2ex]
				0, & k = 0 \ \text{or} \ k \geq N, 
		\end{cases} \\[2ex]
		\left[ \Gamma(a)\tilde a \right]_{k} :&= 	
		\begin{cases}
				0, & 0 \leq k \leq N-1, \\[2ex]
				\left( a \ast \tilde a \right)_{k-1} - \left( a \ast \tilde a \right)_{k+1}, & k \geq N,
		\end{cases}			
\end{align*}
where 
\begin{align}
		\left( \tilde a_{\infty} \right)_{k} := 
		\begin{cases}
				0, & 0 \leq k \leq N-1, \\[2ex]
				\tilde a_{k}, & k \geq N. 
		\end{cases}
		\label{eq:a_inf}
\end{align}
\begin{remark}
		The parameters $N$ and $\nu$ in this section are not to be confused with the vector valued
		ones used throughout this paper. In practice, we set $a = \hat g^{ijl}$, $N = N_{i}$ and
		$M = M_{jl}$. In particular, observe that 
		\begin{align} 
			\label{eq:DFh_N}
			\left[ \Pi_{a}^{ij} \left( DF \left( \hat x \right) - \widehat{DF} \right)h \right]_{k=1}^{N_{i}-1}  &= 
			- \dfrac{L \left( t_{i} - t_{i-1} \right) }{4} \left[ \sum_{l=1}^{n} B \left( \hat g^{ijl} \right) \Pi^{il}_{a}(h) \right]_{k=1}^{N_{i}-1}, \\[2ex]
			\label{eq:DFh_inf}
			\left[ \Pi_{a}^{ij} \left( DF \left( \hat x \right) - \widehat{DF} \right)h \right]_{k=N_{i}}^{\infty} &= 
			- \dfrac{L \left( t_{i} - t_{i-1} \right) }{4} 
			 \displaystyle \left[ \sum_{l=1}^{n} \Gamma \left( \hat g^{ijl} \right)\Pi^{il}_{a} \left( h \right) \right]_{k = N_{i}}^{\infty},
		\end{align}
		by \eqref{eq:Pi_cheb0}.
\end{remark}

We begin by extending $a$ and $\tilde a$ to ``bi-infinite'' sequences, by setting 
$a_{-k} := a_{k}$ and $\tilde a_{-k} := \tilde a_{k}$ for $k \in \NN$, and constructing a bi-infinite
matrix representation for the map 
\begin{align*}
\tilde a \mapsto \left[ \left( a \ast \tilde a \right)_{k-1} - \left( a \ast \tilde a \right)_{k+1} \right]_{k \in \ZZ}.
\end{align*} 
We then convert this bi-infinite matrix representation to an ``one-sided'' matrix representation 
by using appropriate reflections, which in turn will be used to construct the desired matrix representations for $B(a)$ and 
$\Gamma(a)$. 
To be more precise, first observe that 
\newcommand{\WidestEntry}{$a_{1-\tilde M}$}
\newcommand{\SetToWidest}[1]{\makebox[\widthof{\WidestEntry}]{$#1$}}
\begin{gather}
\begin{align*}
	& \left( a \ast \tilde a  \right)_{k \in \ZZ} \\[2ex] & \quad =  
	\left[	\begin{array}{*{12}{c}}	
				\SetToWidest{\ddots} 	& 					 		& 	\SetToWidest{\ddots}	  		 & 	
				\SetToWidest{\ddots}  	& 	\SetToWidest{\ddots}  		& \tikzmarkin[ver=style gray]{col 0} 	 & \SetToWidest{\ddots} &  & & &  
				\\[2ex]
					  				&	a_{\tilde M-1}				& 		  					&	\SetToWidest{\ddots}		& 	\SetToWidest{\ddots}		&
				\SetToWidest{a_{-1}}		&				& a_{1- \tilde M}  		  &						&					         &								
				\\[2ex] \tikzmarkin[hor=style gray_hor]{row 0}
									&				& a_{\tilde M-1}   		  &				    	&	\SetToWidest{\ddots}		&
				\SetToWidest{a_{0}}		&	\SetToWidest{\ddots}		&			  & a_{1- \tilde M} &	   		&									 
			 	\\[2ex]
									&						&			 		  & a_{\tilde M-1}   			& 			&
				\SetToWidest{a_{1}}		&	\SetToWidest{\ddots}		&  \SetToWidest{\ddots}	  & 			    	& a_{1 - \tilde M}   					&			 
				\\[2ex]
									&						& 			 		  & 		           	& \SetToWidest{\ddots}  				& \tikzmarkend{col 0}   
									& \SetToWidest{\ddots}   		& \SetToWidest{\ddots}  	  & \SetToWidest{\ddots}  		   & 				& \SetToWidest{\ddots} 
		\end{array} 
	\right]
	\left[
		\begin{array}{*{14}{c}}
				\\[2ex]
				\vdots
				\\[2ex]
				\tilde a_{-1}
				\\[2ex]	
				\tilde a_{0}  \tikzmarkend{row 0}	
				\\[2ex]	
				\tilde a_{1}
				\\[2ex]	
				\vdots
				\\[2ex]	
		\end{array} 
	\right].
\end{align*}
\end{gather}
Here we have identified elements in $\ell^{1}_{\nu}$ with bi-infinite column vectors (with respect
to the ordering as depicted above). The bandwidth of this bi-infinite matrix is $\tilde M-1$, since
$a_{k}$ vanishes for $\left \vert k \right \vert \geq \tilde M$. 
The shaded regions in grey indicate the position of the ``zeroth'' row and column.
Set $b_{k} := a_{k-1} - a_{k+1}$ for $- \tilde M \leq k \leq \tilde M$, then it follows from the
above expression that 
\newcommand{\WidestEntryTwo}{$b_{-\tilde M}$}
\newcommand{\SetToWidestTwo}[1]{\makebox[\widthof{\WidestEntryTwo}]{$#1$}}
\begin{gather}
\begin{align*}
	& \left[ \left( a \ast \tilde a  \right)_{k-1} - \left( a \ast \tilde a  \right)_{k+1} \right]_{k \in \ZZ} \\[2ex] & \quad =  
	\left[	\begin{array}{*{11}{c}}	
				\SetToWidestTwo{\ddots} 	& 					 		& 	\SetToWidestTwo{\ddots}	  		 & 	
				\SetToWidestTwo{\ddots}  	& 	\SetToWidestTwo{\ddots}  		& \tikzmarkin[ver=style gray]{col 1} 	 & \SetToWidestTwo{\ddots} &  & & &  
				\\[2ex]
					  				&	\SetToWidestTwo{b_{\tilde M}}				& 		  					&	\SetToWidestTwo{\ddots}		& 	\SetToWidestTwo{\ddots}		&
				\SetToWidestTwo{b_{-1}}		&				& b_{- \tilde M}  		  &						&					         &												
				\\[2ex] \tikzmarkin[hor=style gray_hor]{row 1}
									&				& \SetToWidestTwo{b_{\tilde M}}   		  &				    	&	\SetToWidestTwo{\ddots}		&
				\SetToWidestTwo{b_{0}}		&	\SetToWidestTwo{\ddots}		&			  & b_{- \tilde M} &	   		&									 
			 	\\[2ex]
									&						&			 		  & \SetToWidestTwo{b_{\tilde M}}   			& 			&
				\SetToWidestTwo{b_{1}}		&	\SetToWidestTwo{\ddots}		&  \SetToWidestTwo{\ddots}	  & 			    	& b_{- \tilde M}   					&			 
				\\[2ex]
									&						& 			 		  & 		           	& \SetToWidest{\ddots}  				& \tikzmarkend{col 1}   
									& \SetToWidestTwo{\ddots}   		& \SetToWidestTwo{\ddots}  	  & \SetToWidestTwo{\ddots}  		   & 				& \SetToWidestTwo{\ddots} 
		\end{array}
	\right]
	\left[
		\begin{array}{*{14}{c}}
				\\[2ex]
				\vdots
				\\[2ex]
				\tilde a_{-1}
				\\[2ex]	
				\tilde a_{0}  \tikzmarkend{row 1}	
				\\[2ex]	
				\tilde a_{1}
				\\[2ex]	
				\vdots
				\\[2ex]	
		\end{array}
	\right].
\end{align*}
\end{gather}
In particular, this bi-infinite matrix has bandwidth $\tilde M$. 

Next, we convert the latter matrix representation into a one-sided representation on $\NN_{0}$ 
by ``reflecting'' all elements on the left hand-side of the zeroth column to the right and ignoring the rows
with negative indices. This yields  
\newcommand{\WidestEntryThree}{$b_{-\tilde M}$}
\newcommand{\SetToWidestThree}[1]{\makebox[\widthof{\WidestEntryThree}]{$#1$}}
\begin{align}
	 \left[ \left(a \ast \tilde a \right)_{k-1}  -  \left(a \ast \right. \right. & \left. \left. \tilde a \right)_{k+1}  \right]_{k \in \NN_{0}} \nonumber \\[2ex] & =  
	 \left[
	 	\begin{array}{*{7}{c}}
	 			\SetToWidestThree{b_{0}}	&		\SetToWidestThree{\ldots}		& 			\SetToWidestThree{b_{- \tilde M}}		&&&&
				\\[2ex]		
				\SetToWidestThree{\vdots}			&		\SetToWidestThree{\ddots}		& \SetToWidestThree{\ddots}
							&		\SetToWidestThree{\ddots}		&&& 		
				\\[2ex]
				\SetToWidestThree{b_{\tilde M}}		& \SetToWidestThree{\ddots}	 	& 		\SetToWidestThree{b_{0}} &
				\SetToWidestThree{\ddots}  &  \SetToWidestThree{b_{- \tilde M}} &&
				\\[2ex] 
				& \SetToWidestThree{b_{\tilde M}}		& \SetToWidestThree{\ddots}		& 		\SetToWidestThree{b_{0}} & 
				\SetToWidestThree{\ddots} &  \SetToWidestThree{b_{- \tilde M}}   &			
				\\[2ex]
				&		&		\SetToWidestThree{\ddots} 	& 		\SetToWidestThree{\ddots} 	& 	
				\SetToWidestThree{\ddots} & 	\SetToWidestThree{\ddots}  & 	\SetToWidestThree{\ddots}  
		 \end{array}
	 \right] 
	 \left[
	 \begin{array}{*{11}{c}}
	 		\tilde a_{0} \\[2ex]
			\tilde a_{1} \\[2ex]
			\vdots \\[2ex]
			 \vdots \\[2ex]
			 \vdots \\[2ex]
	 \end{array}
	 \right] 
	 \nonumber \\[2ex] & \qquad \qquad + 
	 \begin{bmatrix}
	 		b_{1} & \ldots & b_{\tilde M} \\[1ex]
			\vdots & 	\iddots	& \\[1ex]
			b_{\tilde M}
	 \end{bmatrix}
	 \left[
	  \begin{array}{*{11}{c}}
	  \tilde a_{1} \\[1ex]
	  \vdots \\[1ex]
	  \tilde a_{\tilde M}
	  \end{array}
	  \right]. 
	  \label{eq:S(a)}
\end{align}
Altogether, the sum of the above two matrices, which we will denote by $S(a)$, 
constitutes an infinite dimensional matrix representation
of the map 
\begin{align*}
	\tilde a \mapsto \left[ \left(a \ast \tilde a \right)_{k-1} - \left(a \ast \tilde a \right)_{k+1}  \right]_{k \in \NN_{0}}.
\end{align*} 

Let $\bold{\tilde B}(a) \in \CC^{(N-1) \times \left( \tilde M + N \right) }$ denote the \emph{finite dimensional} submatrix 
of $S(a)$ defined by 
\begin{align*}
		\bold{\tilde B}(a) := \left \{ S(a)_{ij} : 1 \leq i \leq N-1, \ 0 \leq j \leq \tilde M + N -1 \right \}. 
\end{align*}
Note that we are
using the convention that the indexing of the rows and columns start at zero rather than at one.
In view of \eqref{eq:a_inf}, set the elements in the columns of $\bold{\tilde B}(a)$ with index $0 \leq j \leq N-1$ to zero and
let $\bold{B}(a)$ denote the resulting matrix. Finally, let $\Gamma_{\infty}(a)$ denote the infinite
dimensional matrix which consists of the rows of $S(a)$ with index $N$ and higher. Then 
\begin{align}
	\label{eq:BGamma}
	\left[ B(a) \tilde a \right]_{k=1}^{N-1} = \bold{B}(a) \left[ \tilde a_{k} \right]_{k=0}^{\tilde M + N -1}, \quad 
	\left[ \Gamma(a) \tilde a \right]_{ k \geq N} = \Gamma_{\infty}(a) \tilde a
\end{align}
by construction. 

In preparation for the analysis in Section \ref{sec:firstorder}, we show 
that the operator norm of $\Gamma(a)$ can be computed by considering a sufficiently large 
finite dimensional submatrix of $\Gamma_{\infty}(a)$. 
\begin{lemma}[Operator norm of $\Gamma(a)$]
	\label{lemma:Gamma_norm}
		Let $\boldsymbol{\Gamma}_{N}(a) \in \CC^{\left( 3 \tilde M - N +1 \right) \times \left( 2 \tilde M + 1 \right)}$ denote the submatrix 
		$
			\left \{ S(a)_{ij} : N \leq i \leq 3\tilde M, \ 0 \leq j \leq 2 \tilde M \right \}
		$
		of $S(a)$. Then 
		\begin{align*}
				\left \Vert \Gamma (a) \right \Vert_{\mathcal{B} \left( \ell^{1}_{\nu}, \ell^{1}_{\nu} \right)} = 
				\left \Vert 
				\begin{bmatrix}
						\bold{0}_{N \times \left( 2 \tilde M + 1\right)} \\[2ex]
						\boldsymbol{\Gamma}_{N}(a)
				\end{bmatrix}
				\right \Vert_{\mathcal{B} \left( \ell^{1}_{\nu}, \ell^{1}_{\nu} \right)}.
		\end{align*}
		\begin{proof}
				It follows directly from the expression in \eqref{eq:S(a)} that 
				\begin{align*}
						\left \Vert \Gamma (a) \xi_{k'} \right \Vert_{\nu} = 
						\sum_{k = - \tilde M}^{\tilde M} \left \vert b_{k} \right \vert \nu^{k}, 
						\quad \forall k' \geq 2 \tilde M,
				\end{align*}
				where $\left( \xi_{k'} \right)_{k' \in \NN_{0}}$ are the corner points introduced in Definition \ref{def:cornerpoints}. 
This implies that $\left \Vert \Gamma (a) \xi_{k'} \right \Vert_{\nu} = 
\left \Vert \Gamma (a) \xi_{2\tilde{M}} \right \Vert_{\nu}$
for all $k' \geq 2 \tilde{M}$, hence
				\begin{align*}
						\left \Vert \Gamma (a) \right \Vert_{\mathcal{B} \left( \ell^{1}_{\nu}, \ell^{1}_{\nu} \right)} = 
						\max_{0 \leq k' \leq 2 \tilde M } \left \Vert \Gamma (a) \xi_{k'} \right \Vert_{\nu} = 
						\left \Vert 
							\begin{bmatrix}
								\bold{0}_{N \times \left( 2 \tilde M + 1\right)} \\[2ex]
								\boldsymbol{\Gamma}_{N}(a)
						\end{bmatrix}
						\right \Vert_{\mathcal{B} \left( \ell^{1}_{\nu}, \ell^{1}_{\nu} \right)}					,
				\end{align*}
				by Proposition \ref{prop:l1_operator_norm} and the definition of $\boldsymbol{\Gamma}_{N}(a)$. 
		\end{proof}
\end{lemma}
\begin{remark}
	The latter results shows that the operator norm of $\Gamma(a)$
	is determined by its first $2 \tilde M +1$ columns. 
\end{remark}

\subsubsection{First order bounds}
\label{sec:firstorder}
In this section we compute bounds for 
\begin{align}
	\label{eq:Z1_decomp}
	A \left( DF \left( \hat x \right) - \widehat{DF} \right)h
\end{align}
by projecting it onto all the relevant subspaces of $\mathcal{X}_{\nu}$. For notational convenience, we shall 
write $y = \left(DF \left( \hat x \right) - \widehat{DF} \right)h$ throughout this section. 
We start by computing the difference between the exact and approximate derivative. 

A straightforward computation shows  that 
\begin{align}
	\label{eq:Pi_t0}
	\Pi_{t_{0}} (y)  = 
	2 \sum_{k=N_{1}}^{\infty} \left(-1\right)^{k} \tilde a^{1}_{k} - \sum_{k \not \in \mathcal{K}^{u}}
	\tilde p_{k} \hat \theta^{k}, \quad
	\Pi_{t_{m}}  (y)  = 
	2 \sum_{k=N_{m}}^{\infty} \tilde a^{m}_{k} - \sum_{k \not \in \mathcal{K}^{s}}
	\tilde q_{k} \hat \phi^{k}.
\end{align}
Furthermore, $\Pi(y) = 0$ for $\Pi \in \left \{ \Pi_{\hat p_{1}}, \Pi_{\hat q_{1}}, \Pi_{\eta} \right \}$, 
since the equations associated to $\Pi F$ are linear and only depend on elements in the finite dimensional
subspace $\mathcal{X}^{NK}$. Next, set 
\begin{align*}
	\boldsymbol{k}_{ijl} := \left[0, \left(M_{jl}+1\right)\left(N_{i}-1\right)+1 \right] \cap \NN_{0}, \quad
	1 \leq i \leq m, \ 1 \leq j, l \leq n, 
\end{align*}
where $M_{jl} = \text{order} \left( \dfrac{ \partial g_{j} }{ \partial x_{l} } \right)$, then it follows from
\eqref{eq:DFh_N}, \eqref{eq:DFh_inf} and \eqref{eq:BGamma} that  
\begin{align}
	\label{eq:Pi_cheb}
	 \left( \Pi_{a}^{ij} (y) \right)_{k} =
	\begin{cases}
		\displaystyle \sum_{l=N_{i}}^{\infty} \left( -1\right)^{l} \left[ \tilde a^{i}_{l} \right]_{j} - \sum_{l=N_{i-1}}^{\infty} \left[\tilde a^{i-1}_{l} \right]_{j}, & k=0, \\[4ex]
		- \dfrac{L \left( t_{i} - t_{i-1} \right) }{4} 
		\left[ \displaystyle \sum_{l=1}^{n} \bold{B} \left( \hat g^{ijl} \right) \Pi_{a}^{il \boldsymbol{k}_{ijl} } \left( h \right) \right]_{k}, & 1 \leq k \leq N_{i}-1, \\[4ex]
		- \dfrac{L \left( t_{i} - t_{i-1} \right) }{4} 
		\left[ \displaystyle \sum_{l=1}^{n} \Gamma_{\infty} \left( \hat g^{ijl} \right)\Pi^{il}_{a} \left( h \right) \right]_{k} & k \geq N_{i}, 
	\end{cases}
\end{align}
for $2 \leq i \leq m$, see Section \ref{sec:conv_terms}. Here the matrix-vector product $\bold{B} \left( \hat g^{ijl} \right) \Pi_{a}^{il \boldsymbol{k}_{ijl} } \left( h \right)$ 
is interpreted by using the identification 
\begin{align*}
	\Pi_{a}^{il \boldsymbol{k}_{ijl} } \left( h \right) \simeq 
	\begin{bmatrix} 
		\left[ \tilde a^{i}_{0} \right]_{l} & \ldots & \left[ \tilde a^{i}_{ \max {\boldsymbol{k}_{ijl}} } \right]_{l}
	\end{bmatrix}^{T}.
\end{align*}
The same formula holds for $i=1$ and $k \in \NN$. In particular, 
if $i=1$, then there is no component to consider for $k=0$.
Finally, we compute that 
\begin{align}
		\left( \Pi_{P} (y) \right)_{k} =
		\begin{cases}
			\bold{0}_{n}, & k \in \mathcal{K}^{u}, \\[2ex]
			-DC_{k} \left( \hat p \right) \tilde p, & k \not \in \mathcal{K}^{u}, 
		\end{cases} \qquad
		\left( \Pi_{Q} (y) \right)_{k} =
		\begin{cases}
			\bold{0}_{n}, & k \in \mathcal{K}^{s}, \\[2ex]
			-DC_{k} \left( \hat q \right) \tilde q, & k \not \in \mathcal{K}^{s}.
		\end{cases} 
	\label{eq:Pi_parm}
\end{align}

Altogether, the above formulae give rise to the decomposition 
\begin{align}
	y = \sum_{  \tilde \jmath =1 }^{n} &\left( 
		\Pi^{ \tilde \jmath }_{t_{0}}(y) + \Pi^{ \tilde \jmath }_{t_{m}}(y) + 
		\sum_{\tilde \imath=2}^{m} \Pi_{a}^{ \tilde \imath  \tilde \jmath  0}(y) + 
		\sum_{\tilde \imath=1}^{m} \left[ \Pi_{a}^{ \tilde \imath  \tilde \jmath  \left[ 1 : N_{\tilde \imath}-1 \right]}(y) + 
		\Pi_{a}^{ \tilde \imath \tilde \jmath  \left[N_{\tilde \imath} : \infty \right)}(y) \right] 
		\right. \nonumber \\[2ex] & \qquad \left.
		+ \Pi^{ \tilde \jmath  \mathcal{K}_{c}^{u}}_{P}(y) + \Pi^{ \tilde \jmath  \mathcal{K}_{c}^{s}}_{Q}(y)
		\right),
	\label{eq:decomp}
\end{align}
where we have set 
\begin{align*}
	\mathcal{K}^{u}_{c} :&= \NN_{0}^{n_{u}} \setminus \mathcal{K}^{u}, \quad 
	\mathcal{K}^{s}_{c} := \NN_{0}^{n_{u}} \setminus \mathcal{K}^{s}, \quad  \\[2ex]
	\left[1:N_{\tilde \imath}-1\right] :&= \left[1,N_{\tilde \imath}-1\right] \cap \NN, \quad
	\left[N_{\tilde \imath} : \infty \right) := \left[N_{\tilde \imath}, \infty \right) \cap \NN.
\end{align*}
The strategy is to compute
bounds for \eqref{eq:Z1_decomp} by individually composing each term in the above decomposition with $A$
and analyzing the associated projections into the domain. 

\begin{remark}
	Observe that $A\Pi^{ \tilde \imath  \tilde \jmath  \left[N_{\tilde \imath} : \infty \right)}_{a}$, 
	$A\Pi^{ \tilde \jmath  \mathcal{K}_{c}^{u}}_{P}$ and $A\Pi^{ \tilde \jmath  \mathcal{K}_{c}^{s}}_{Q}$, are 
	``diagonal'' and ``uncoupled'' in the sense that 
	the only nonzero projections into the domain are 
	$\Pi_{a}^{ \tilde \imath  \tilde \jmath } A\Pi^{ \tilde \imath  \tilde \jmath  \left[N_{\tilde i} : \infty \right)}$, 
	$\Pi^{ \tilde \jmath }_{P} A\Pi^{ \tilde \jmath  \mathcal{K}_{c}^{u}}_{P}$ and 
	$\Pi^{ \tilde \jmath }_{Q} A\Pi^{ \tilde \jmath  \mathcal{K}_{c}^{s}}_{Q}$.
\end{remark}

\paragraph{Boundary conditions}
We start by considering the terms associated to $A\Pi^{ \tilde \jmath }_{t_{0}}(y)$, $A\Pi^{ \tilde \jmath }_{t_{m}}(y)$ and
$A \Pi^{ \tilde \imath  \tilde \jmath  0}(y)$, which are related to the boundary conditions. 

\begin{lemma}
	\label{lemma:boundary}
	Let $\Pi \in \mathbb{P}$, $2 \leq  \tilde \imath \leq m$ and $1 \leq  \tilde \jmath  \leq n$, then 
	\begin{align}
		\label{eq:b_t0}
		\left \Vert \Pi A \Pi^{ \tilde \jmath }_{t_{0}}(y) \right \Vert_{\Pi \left( \mathcal{X}_{\nu} \right)}
		&\leq \left \Vert \Pi A_{NK} \Pi^{ \tilde \jmath }_{t_{0}} \right \Vert_{\Pi \left( \mathcal{X}_{\nu} \right)}
		\left( \nu_{1}^{-N_{1}} + \max_{1 \leq l \leq n_{u} } \left \vert \dfrac{ \hat \theta_{l} }{ \nu_{u}} \right \vert^{ K^{u}_{l} +1 } \right), \\[2ex]
		\label{eq:b_tm}
		\left \Vert \Pi A \Pi^{ \tilde \jmath }_{t_{m}}(y) \right \Vert_{\Pi \left( \mathcal{X}_{\nu} \right)}
		&\leq \left \Vert \Pi A_{NK} \Pi^{ \tilde \jmath }_{t_{m}} \right \Vert_{\Pi \left( \mathcal{X}_{\nu} \right)}
		\left( \nu_{m}^{-N_{m}} + \max_{1 \leq l \leq n_{s} } \left \vert \dfrac{ \hat \phi_{l} }{ \nu_{s}} \right \vert^{ K^{s}_{l} +1 } \right), \\[2ex]
		\label{eq:b_ij0}
		\left \Vert \Pi A \Pi^{ \tilde \imath  \tilde \jmath  0}_{a}(y) \right \Vert_{\Pi \left( \mathcal{X}_{\nu} \right)}
		&\leq \left \Vert \Pi A_{NK} \Pi^{ \tilde \imath  \tilde \jmath  0}_{a} \right \Vert_{\Pi \left( \mathcal{X}_{\nu} \right)}
		\left( \nu_{ \tilde \imath}^{-N_{\tilde \imath} } + \nu_{\tilde \imath-1}^{-N_{ \tilde \imath - 1}} \right),
	\end{align}
	where $\left \Vert \cdot \right \Vert_{\Pi \left( \mathcal{X}_{\nu} \right)}$ 
	denotes the corresponding norm on $\Pi \left( \mathcal{X}_{\nu} \right)$. 
	\begin{proof}
		First observe that 
		\begin{align*}
			\left \Vert a^{i} \mapsto 2 \sum_{k=N_{i}}^{\infty} a^{i}_{k} \right \Vert
			_{\mathcal{B} \left( \ell^{1}_{\nu_{i}}, \CC \right)} = 
			\left \Vert a^{i} \mapsto 2 \sum_{k=N_{i}}^{\infty} (-1)^{k} a^{i}_{k} \right \Vert 
			_{\mathcal{B} \left( \ell^{1}_{\nu_{i}}, \CC \right)} = 
			\nu_{i}^{-N_{i}}, \quad 1 \leq i \leq m,
		\end{align*}
		by Proposition \ref{prop:l1_operator_norm}. Similarly, 
		\begin{align*}
			\left \Vert p \mapsto \sum_{k \not \in \mathcal{K}^{u}} p_{k} \hat \theta^{k} \right \Vert
			_{\mathcal{B} \left( W^{1}_{\nu_{u}}, \CC \right)} &= 
			\sup \left \{ \left \vert \frac{ \hat \theta }{ \nu_{u} } \right \vert^{k} :  \ k \in \NN_{0}^{n_{u}}, \ \exists 1 \leq l \leq n_{u} \ \text{such that} \ 				k_{l} \geq K_{l}^{u} + 1 \right \} \\[2ex]
			& \leq \max_{1 \leq l \leq n_{u} } \left \vert \frac{ \hat \theta_{l} }{ \nu_{u} } \right \vert^{K_{l}^{u}+1} ,
		\end{align*}
		by Proposition \ref{prop:operator_norm_mult}, where in the last line we used that $\hat \theta \in \mathbb{B}_{\nu_{u}}$. Finally, the above bounds and
		the expressions in \eqref{eq:Pi_t0} and \eqref{eq:Pi_cheb} show that 
		\begin{align*} 
			\left \vert \Pi^{ \tilde \jmath }_{t_{0}} \left( y \right) \right \vert &\leq \nu_{1}^{-N_{1}} 
				+ \max_{1 \leq l \leq n_{u} } \left \vert \dfrac{ \hat \theta_{l} }{ \nu_{u}} \right \vert^{ K^{u}_{l} +1 }, \quad
			\left \vert \Pi^{ \tilde \jmath }_{t_{m}} \left( y \right) \right \vert \leq \nu_{1}^{-N_{m}} 
				+ \max_{1 \leq l \leq n_{s}}  \left \vert \dfrac{ \hat \phi_{l} }{ \nu_{s}} \right \vert^{ K^{s}_{l} +1 }, \\[3ex]
			\left \vert \Pi^{ \tilde \imath  \tilde \jmath 0} \left( y \right) \right \vert &\leq 
			\nu_{\tilde \imath}^{-N_{\tilde \imath}} +  \nu_{\tilde \imath-1}^{-N_{\tilde \imath-1}},
		\end{align*}
		which proves the statement. 
	\end{proof}
\end{lemma}
\begin{remark}
		The computation of the stated bounds is finite for each $\Pi \in \mathbb{P}$, 
		since $A_{NK}$ is a finite dimensional matrix.  
\end{remark}

\paragraph{Chebyshev coefficients}
Next, we consider the terms associated to $A\Pi^{ \tilde \imath  \tilde \jmath  \NN}_{a}(y)$.
We start with the observation that 
\begin{align}
	\label{eq:AB(g)}
	A\Pi^{ \tilde \imath  \tilde \jmath  \left[ 1 : N_{\tilde \imath}-1 \right]}_{a} (y) = 
	- \frac{L \left( t_{\tilde \imath} - t_{\tilde \imath-1} \right)}{4} 
	\sum_{l=1}^{n} A_{NK} \Pi_{a}^{ \tilde \imath  \tilde \jmath  \left[ 1 : N_{\tilde \imath}-1 \right]} 
	\bold{B} \left( \hat g^{  \tilde \imath   \tilde \jmath  l } \right) 
	\Pi_{a}^{  \tilde \imath  l \boldsymbol{k}_{ \tilde \imath  \tilde \jmath  l} } \left( h \right)
\end{align}
by \eqref{eq:Pi_cheb}. Note that 
$A_{NK} \Pi^{ \tilde \imath  \tilde \jmath  \left[ 1 : N_{\tilde i}-1 \right]}_{a} 
\bold{B} \left( \hat g^{  \tilde \imath  \tilde \jmath  l } \right)$ is a
\emph{finite dimensional} matrix which can be explicitly computed on a computer. 
In particular, 
\begin{align*}
	\Pi A_{NK} \Pi_{a}^{ \tilde \imath  \tilde \jmath  \left[ 1 : N_{\tilde \imath}-1 \right]} \bold{B} \left( \hat g^{  \tilde \imath  \tilde \jmath  l } \right),
	\quad \Pi \in \mathbb{P}, 
\end{align*} 
corresponds to a finite dimensional matrix representation
of a linear operator on $\ell^{1}_{\nu_{\tilde \imath} }$. Hence the computation of its operator norm is finite. 
\begin{lemma}[Scalar and Taylor projections]
	\label{lemma:Z_ij}
	Let $1 \leq  \tilde \imath \leq m$, $1 \leq  \tilde \jmath  \leq n$ and  
	 $\Pi \in \left \{ \Pi^{j}_{\theta}, \Pi^{j}_{\phi}, \Pi_{\lambda^{u}}^{j}, \Pi^{j}_{\lambda^{s}} 
	    \Pi^{j}_{P}, \Pi^{j}_{Q} \right \}$, then 
	\begin{align*}
		\left \Vert \Pi A \Pi_{a}^{ \tilde \imath   \tilde \jmath  \NN}(y) \right \Vert_{\Pi \left( \mathcal{X}_{\nu} \right)} \leq 
		\frac{L \left( t_{\tilde \imath} - t_{\tilde \imath-1} \right)}{4} 
		\sum_{l=1}^{n} \left \Vert 
		\Pi A_{NK} \Pi_{a}^{ \tilde \imath   \tilde \jmath  \left[ 1 : N_{\tilde i}-1 \right]} 
		\bold{B} \left( \hat g^{  \tilde \imath   \tilde \jmath  l } \right) \right \Vert_{
		\mathcal{B} \left( \ell^{1}_{\nu_{\tilde \imath}}, \Pi \left( \mathcal{X}_{\nu} \right) \right)}.
	\end{align*}
	\begin{proof}
			It suffices to observe that  
			\begin{align*}
					\Pi A \Pi_{a}^{ \tilde \imath   \tilde \jmath  \NN}(y) = 
					\Pi A_{NK} \Pi_{a}^{ \tilde \imath   \tilde \jmath \left[1: N_{\tilde \imath} -1 \right]}  (y),
					\quad \Pi \in \left \{ \Pi^{j}_{\theta}, \Pi^{j}_{\phi}, \Pi_{\lambda^{u}}^{j}, \Pi^{j}_{\lambda^{s}} 
	   			     \Pi^{j}_{P}, \Pi^{j}_{Q} \right \},
	   	   \end{align*} 
		   by construction of the approximate inverse $A$. Hence the result follows directly from 
		   \eqref{eq:AB(g)}.
	\end{proof}
\end{lemma}

To analyze the terms $\Pi_{a}^{ij} A \Pi_{a}^{ \tilde \imath \tilde \jmath  \NN}(y)$ for $1 \leq i \leq m$, $1 \leq j \leq n$, 
we first derive a more explicit expression for the tail $A\Pi_{a}^{ \tilde \imath  \tilde \jmath \left[ N_{\tilde \imath} : \infty \right)}(y)$. 
For this purpose, define a (infinite dimensional) diagonal matrix $D^{i}_{\infty}$ by 
\begin{align*}
	D^{i}_{\infty} := 
	\begin{bmatrix}
		\dfrac{1}{N_{i}} 	& 				& \\
					& \dfrac{1}{N_{i}+1} 		& \\
					& 				& \ddots
	\end{bmatrix}, \quad 1 \leq i \leq m. 
\end{align*}
Then it follows from \eqref{eq:Pi_cheb} and the definition of the approximate inverse that
\begin{align*}
	\Pi^{ij \left[N_{i}: \infty \right)}_{a} A\Pi_{a}^{ \tilde \imath  \tilde \jmath \left[ N_{\tilde \imath} : \infty \right)}(y) = 
	\begin{cases}
		0, & \left(i,j\right) \not = \left( \tilde \imath, \tilde \jmath \right), \\[2ex]
		- \displaystyle \frac{L \left(t_{i} - t_{i-1} \right) }{4} \sum_{l=1}^{n} D^{i}_{\infty} \Gamma_{\infty} \left( \hat g^{ijl} \right) \Pi_{a}^{il}(h), 
		& \left(i,j\right) = \left( \tilde \imath, \tilde \jmath \right).
	\end{cases}
\end{align*}
Altogether, by combining the latter result with \eqref{eq:AB(g)}, we conclude that 
if $\left(i,j\right) \not =  \left( \tilde \imath, \tilde \jmath \right)$, then 
\begin{align}
	\label{eq:Pi_cheb_case_1}
	\Pi_{a}^{ij} A \Pi_{a}^{ \tilde \imath \tilde \jmath \NN}(y) = 
	- \frac{L \left( t_{\tilde \imath} - t_{\tilde \imath-1} \right)}{4} 
	\sum_{l=1}^{n} \Pi_{a}^{ij \left[0:N_{i}-1\right]} A_{NK} \Pi_{a}^{ \tilde \imath  \tilde \jmath  \left[ 1 : N_{\tilde \imath}-1 \right]} 
	\bold{B} \left( \hat g^{  \tilde \imath   \tilde \jmath  l } \right) 
	\Pi_{a}^{  \tilde \imath  l \boldsymbol{k}_{ \tilde \imath  \tilde \jmath  l} } \left( h \right).
\end{align}
Otherwise, if $\left(i,j\right) =  \left( \tilde \imath, \tilde \jmath \right)$, then 
\begin{align}
	& \Pi_{a}^{ij} A \Pi_{a}^{ ij \NN}(y)  \nonumber \\[2ex] & \quad  
	= -  \frac{L \left(t_{i} - t_{i-1} \right) }{4} 
	 \displaystyle \sum_{l=1}^{n}
		\begin{bmatrix}
			\begin{array}{cc}
				\Pi_{a}^{ij \left[0:N_{i}-1\right]} A_{NK} \Pi_{a}^{ ij  \left[ 1 : N_{i}-1 \right]} 
				\bold{B} \left( \hat g^{ ij  l } \right) & \bold{0}_{N \times \infty}
			\end{array} \\[2ex]
			D^{i}_{\infty} \Gamma_{\infty} \left( \hat g^{ijl} \right) 
		\end{bmatrix}  \Pi_{a}^{il}(h).
	\label{eq:Pi_cheb_case_2}
\end{align}

We are now ready to compute the desired bounds.   
Before we proceed, observe that the operator norms of the infinite dimensional matrices in \eqref{eq:Pi_cheb_case_2} 
can be computed by considering sufficiently large \emph{finite dimensional} submatrices 
by the same reasoning as in Lemma \ref{lemma:Gamma_norm}. The details are given in the lemma below. 
\begin{lemma}[Chebyshev projections]
	\label{lemma:Z_ijij}
	Let $1 \leq \tilde \imath, i \leq m$, $1 \leq  \tilde \jmath , j \leq n$. 	
	If $\left(  \tilde \imath ,  \tilde \jmath  \right) \not = \left(i,j\right)$, then 
	\begin{align*}
	 \left \Vert \Pi_{a}^{ij}A \Pi_{a}^{ \tilde \imath  \tilde \jmath  \NN} (y) \right \Vert_{\nu_{i}} 
		\leq
		\frac{L \left( t_{ \tilde \imath } - t_{\tilde i-1} \right) }{4}  \sum_{l=1}^{n} \left \Vert
			\Pi_{a}^{ij \left[0:N_{i}-1\right]} A_{NK} \Pi_{a}^{ \tilde \imath  \tilde \jmath  \left[ 1 : N_{\tilde \imath}-1 \right]} 
			\bold{B} \left( \hat g^{  \tilde \imath  \tilde \jmath  l } \right)  
		\right \Vert_{\mathcal{B} \left( \ell^{1}_{\nu_{\tilde \imath}}, \ell^{1}_{\nu_{i}} \right)}.
	\end{align*}
	Otherwise, if $\left(  \tilde \imath ,  \tilde \jmath  \right)  = \left(i,j\right)$, then 
	\begin{align*}
	 	& \left \Vert \Pi_{a}^{ij}A \Pi_{a}^{ \tilde \imath  \tilde \jmath  \NN} (y) \right \Vert_{\nu_{i}} \\[2ex] & \quad \leq
		\frac{L \left(t_{i} - t_{i-1} \right) }{4} 
		 \displaystyle \sum_{l=1}^{n}
			\left \Vert \begin{bmatrix}
				\begin{array}{cc}
					\Pi_{a}^{ij \left[0:N_{i}-1\right]} A_{NK} \Pi_{a}^{ ij  \left[ 1 : N_{i}-1 \right]} 
					\bold{B} \left( \hat g^{ ij  l } \right) & \bold{0}_{N \times \left( \tilde M_{ijl} - N_{i} +1 \right)}
				\end{array} \\[2ex]
				\bold{D}^{ijl}_{N} \ \bold{\Gamma}_{N} \left( \hat g^{ijl} \right) 
			\end{bmatrix} 
			\right \Vert_{\mathcal{B} \left( \ell^{1}_{\nu_{i}}, \ell^{1}_{\nu_{i}} \right)},
	\end{align*}	
	where 
	\begin{align*}
		\bold{D}^{ijl}_{N} := 
			\begin{bmatrix}
				\dfrac{1}{N_{i}} 	& 				& \\
							& \ddots 			& \\
							& 				& \dfrac{1}{ 3 \tilde M_{ijl}}
			\end{bmatrix}, 	
			\quad \tilde M_{ijl} :=M_{jl} \left( N_{i} - 1 \right) +1, 
	\end{align*}
	and $\bold{\Gamma}_{N} \left( \hat g^{ijl} \right)$ is defined in Lemma \ref{lemma:Gamma_norm}. 
	\begin{proof}
		The statement for $\left(  \tilde \imath ,  \tilde \jmath  \right) \not = \left(i,j\right)$ follows directly from 
		\eqref{eq:Pi_cheb_case_1} and the triangle inequality. To prove the result for
		$\left(  \tilde \imath ,  \tilde \jmath  \right) = \left(i,j\right)$, consider the linear operator 
		\begin{align}
			\left[ \tilde a^{i} \right]_{l} \mapsto 
			\begin{bmatrix}
				\begin{array}{cc}
					\Pi_{a}^{ij \left[0:N_{i}-1\right]} A_{NK} \Pi_{a}^{ ij  \left[ 1 : N_{i}-1 \right]} 
					\bold{B} \left( \hat g^{ ij  l } \right) & \bold{0}_{N \times \infty}
				\end{array} \\[2ex]
				D^{i}_{\infty} \Gamma_{\infty} \left( \hat g^{ijl} \right) 
			\end{bmatrix} \left[ \tilde a^{i} \right]_{l}. 
			\label{eq:PiAPi}
		\end{align}
		A similar computation as in Lemma \ref{lemma:Gamma_norm} shows that 
		\begin{align*}
			& \left \Vert \begin{bmatrix}
				\begin{array}{cc}
					\Pi_{a}^{ij \left[0:N_{i}-1\right]} A_{NK} \Pi_{a}^{ ij  \left[ 1 : N_{i}-1 \right]} 
					\bold{B} \left( \hat g^{ ij  l } \right) & \bold{0}_{N \times \infty}
				\end{array} \\[2ex]
				D^{i}_{\infty} \Gamma_{\infty} \left( \hat g^{ijl} \right) 
			\end{bmatrix} \xi_{k'} \right \Vert_{\nu_{i}}	\\[2ex] 
			& \qquad = 	
			\sum_{k= - \tilde M_{ijl}}^{\tilde M_{ijl}} 
			\left \vert \hat g^{ijl}_{k-1} - \hat g^{ijl}_{k+1} \right \vert 
			\dfrac{ \nu_{i}^{k}}{  k  + k' }, \quad
			\forall k' \geq 2 \tilde M_{ijl},
		\end{align*}
		where $\left( \xi_{k'} \right)_{k' \in \NN_{0}}$ are the corner points introduced in Definition \ref{def:cornerpoints}.
		Note that the latter quantity is decreasing for $k' \geq 2 \tilde M_{ijl}$. Hence, by Proposition \ref{prop:l1_operator_norm}, the operator norm of 
		\eqref{eq:PiAPi} is completely determined by its columns with index $0 \leq k' \leq 2 \tilde M_{ijl}$. The corresponding 
		submatrix is given by 
		\begin{align*}
			 \begin{bmatrix}
				\begin{array}{cc}
					\Pi_{a}^{ij \left[0:N_{i}-1\right]} A_{NK} \Pi_{a}^{ ij  \left[ 1 : N_{i}-1 \right]} 
					\bold{B} \left( \hat g^{ ij  l } \right) & \bold{0}_{N \times \left( \tilde M_{ijl}  - N_{i} +1 \right)}
				\end{array} \\[2ex]
				\bold{D}^{ijl}_{N} \ \bold{\Gamma}_{N} \left( \hat g^{ijl} \right) 
			\end{bmatrix}.
		\end{align*}
		Therefore, the result now follows from \eqref{eq:Pi_cheb_case_2} and the triangle inequality. 
	\end{proof}
\end{lemma}

\paragraph{Taylor coefficients}
Finally, we consider the terms $A \Pi^{ \tilde \jmath  \mathcal{K}_{c}^{u}}_{P}(y)$ and $A \Pi^{ \tilde \jmath  \mathcal{K}_{c}^{s}}_{Q}(y)$.
In particular, recall that the only nonzero projections in this case are $\Pi^{ \tilde \jmath }_{P}A \Pi^{ \tilde \jmath  \mathcal{K}_{c}^{u}}_{P}(y)$
and $\Pi^{ \tilde \jmath }_{Q}A \Pi^{ \tilde \jmath  \mathcal{K}_{c}^{u}}_{Q}(y)$.
To study these terms, we will use the following result: 
\begin{lemma}
	\label{lemma:Lambda}
	Let $M \in \NN$ and set $\mathcal{K}^{u}_{M}:= \left \{ k \in \NN^{n_{u}}_{0}: k_{i} \leq MK^{u}_{i}, \ 1 \leq i \leq n_{u} \right \}$.
	Suppose $p \in W^{1}_{\nu_{u}}$ satisfies $p_{k} = 0$ whenever $k \not \in \mathcal{K}^{u}_{M}$ and define an operator
	$\Lambda_{u} (p) : W^{1}_{\nu_{u}} \rightarrow W^{1}_{\nu_{u}}$ by 
	\begin{align*}
		\left(  \Lambda_{u}(p)w \right)_{k} := 
		\begin{cases}
			0, & k \in \mathcal{K}^{u}, \\[1ex]
			\left \langle k, \hat \lambda^{u} \right \rangle^{-1} \left( p \boldsymbol{\ast}w \right)_{k}, & k \not \in \mathcal{K}^{u}.
		\end{cases}
	\end{align*}
	Then $\Lambda_{u}(p)$ is bounded and 
	\begin{align*}
		& \left \Vert \Lambda_{u}(p) \right \Vert_{ \mathcal{B} \left( W^{1}_{\nu_{u}}, W^{1}_{\nu_{u}} \right) }
		\\[2ex] & \quad  \leq 
		\max_{ \substack{ 1 \leq j \leq n_{u} } } 
		\left \{ 
			\sum_{ \substack{ k \in \mathcal{K}^{u}_{M} \\  k_{j} \geq K^{u}_{j} + 1 - l_{j} } }
			\frac{\left \vert p_{k} \right \vert  \nu_{u}^{\left \vert k \right \vert} }{
			\left( k_{j} + l_{j} \right) \left \vert \normalfont \text{Re} \left( \hat \lambda^{u}_{j} \right) \right \vert + 
			\sum_{ \substack{1 \leq i \leq n_{u} \\ i \not = j} } k_{i} \left \vert \normalfont \text{Re} \left( \hat \lambda^{u}_{i} \right) 
			\right \vert } :  0 \leq l_{j} \leq K^{u}_{j} +1\right \}.
	\end{align*}
 
	\begin{proof}
		It follows directly from the Banach algebra estimate that $\Lambda^{u}(p)$ is bounded. To obtain the stated bound
		for operator norm, first note that 
		\begin{align*}
			\left( p \boldsymbol{\ast}w \right)_{k} &= 
			\sum_{ \substack{ \alpha + \beta = k \\ \alpha \in \mathcal{K}^{u}_{M} \\ \beta \in \NN_{0}^{n_{u}} } } p_{\alpha} w_{\beta} =
			\sum_{ \max \left \{0, k_{i} - MK^{u}_{i} \right \} \leq \beta_{i} \leq k_{i}	} p_{k - \beta} w_{\beta}, \quad		 		
		\end{align*}
		for any $k \in \NN_{0}^{n_{u}}$.
		In particular, 
		\begin{align*}
			\left( p \boldsymbol{\ast} \boldsymbol{\xi}_{l} \right)_{k} = 
			\begin{cases}
				p_{k - l} \nu_{u}^{-\left \vert l \right \vert} , & \max \left \{0, k_{i} - MK^{u}_{i} \right \} \leq l_{i} \leq k_{i}, \ 1 \leq i \leq 					n_{u}, \\[2ex]
				0, & \mbox{otherwise},
			\end{cases}
		\end{align*}
		for any $l \in \NN_{0}^{n_{u}}$, 
		where $\left( \boldsymbol{\xi}_{l} \right)_{l \in \NN_{0}^{n_{u}}}$ are the corner points 
		introduced in Definition \ref{def:cornerpoints_array}. 
		Consequently, 
		\begin{align}
			\left \Vert \Lambda_{u}(p) \boldsymbol{\xi}_{l}  \right \Vert_{\nu_{u}} &= 
			\sum_{ \substack{ k \in \mathcal{K}_{c}^{u} \cap \left( l + \mathcal{K}^{u}_{M} \right)  } }
			\left \vert p_{k-l} \right \vert
			\left \vert \left \langle k, \hat \lambda^{u} \right \rangle \right \vert^{-1}
			 \nu_{u}^{ \left \vert k \right \vert - \left \vert l \right \vert} \nonumber \\[2ex] \quad &=  
			 \sum_{ \substack{ k \in \mathcal{K}_{M}^{u} \\ \exists j : k_{j} \geq K^{u}_{j} +1 - l_{j} } }
			 \left \vert p_{k} \right \vert \left \vert \left \langle k+l, \hat \lambda^{u} \right \rangle \right \vert^{-1}
			 \nu_{u}^{\left \vert k \right \vert}. 
			 \label{eq:Taylor_tail}
		\end{align}		
		
		Next, observe that for any $l \in \NN_{0}^{n_{u}}$ and $1 \leq j \leq n_{u}$, 
		\begin{align*}
			\left \vert \left \langle k + l, \hat \lambda^{u} \right \rangle \right \vert & = 
			\left \vert \sum_{i=1}^{n_{u}} \left( k_{i} + l_{i} \right) \hat \lambda^{u}_{i} \right \vert 
			\\[2ex] & \geq 
			\left \vert \sum_{i=1}^{n_{u}} \left( k_{i} + l_{i} \right) \text{Re} \left( \hat \lambda^{u}_{i} \right) \right \vert
			\\[2ex] &\geq \left( k_{j} + l_{j} \right)  \left \vert \text{Re} \left( \hat \lambda^{u}_{j} \right) \right \vert
			+ \sum_{ \substack{1 \leq i \leq n_{u} \\ i \not = j} } 
			k_{i} \left \vert \normalfont \text{Re} \left( \hat \lambda^{u}_{i} \right) \right \vert, 
		\end{align*}
		where in the last line we used the fact that all the $\hat \lambda^{u}_{i}$ have the same sign.	
		Therefore, the term in \eqref{eq:Taylor_tail} is bounded by 
		\begin{align*}
			\max_{1 \leq j \leq n_{u}}
			\sum_{ \substack{ k \in \mathcal{K}^{u}_{M} \\  k_{j} \geq K^{u}_{j} + 1 - l_{j} } }
			\frac{\left \vert p_{k} \right \vert  \nu_{u}^{\left \vert k \right \vert} }{
			\left( k_{j} + l_{j} \right) \left \vert \normalfont \text{Re} \left( \hat \lambda^{u}_{j} \right) \right \vert + 
			\sum_{ \substack{1 \leq i \leq n_{u} \\ i \not = j} } k_{i} \left \vert \normalfont \text{Re} \left( \hat \lambda^{u}_{i} \right) \right 				\vert }.
		\end{align*}			
		Finally, note that for any fixed $1 \leq j \leq n_{u}$, the above sum is strictly decreasing for 
		$l_{j} \geq K_{j}^{u}+1$. Hence the desired result now follows from Proposition \ref{prop:operator_norm_mult}. 
	\end{proof}
\end{lemma} 
\begin{remark}
	A similar statement holds for the map associated to the stable manifold. The corresponding operator 
	(for $q \in W^{1}_{\nu_{s}}$) is denoted by $\Lambda_{s} \left(q \right) : W^{1}_{\nu_{s}} \rightarrow  W^{1}_{\nu_{s}}$.
\end{remark}

We are now ready to compute the required bounds. To this end, observe that 
\begin{align}
	\label{eq:Dc_taylor}
	D\left[C\right]_{j} \left( \hat p \right) \tilde p = \sum_{l=1}^{n} \hat G_{u}^{jl} \boldsymbol{\ast} \left[ \tilde p \right]_{l}, \quad 1 \leq j \leq n, 
\end{align}
by the reasoning in Remark \ref{remark:Dc_comp}, 
where $\left [ C \right]_{j}$ denotes the $j$-th component of $C$ and $\hat G_{u}^{jl}$ are the Taylor coefficients of 
\begin{align*}
	\frac{ \partial g_{j} }{ \partial x_{l} } \left( x \mapsto \sum_{k \in \mathcal{K}^{u}} \hat p_{k} x^{k} \right), 
	\quad 1 \leq l \leq n.
\end{align*}
The coefficients $\hat G_{s}^{jl}$ associated to the stable manifold are defined similarly. 
\begin{lemma}[Projection onto the Taylor coefficients]
	\label{lemma:tail_taylor}
	Let $1 \leq j \leq n$, then 
	\begin{align}
		\label{eq:Lambda_u}
		 \left \Vert \Pi_{P}^{j} A \Pi_{P}^{j \mathcal{K}^{u}_{c}}(y)
		 \right \Vert_{W^{1}_{\nu_{u}} } &\leq
		\sum_{l=1}^{n} \left \Vert \Lambda_{u} \left( \hat G_{u}^{jl} \right) \right \Vert_{
		\mathcal{B} \left( W^{1}_{\nu_{u}}, W^{1}_{\nu_{u}} \right)}, \\[2ex]
		\label{eq:Lambda_s}
		 \left \Vert \Pi_{Q}^{j} A \Pi_{Q}^{j \mathcal{K}^{s}_{c}}(y)
		 \right \Vert_{W^{1}_{\nu_{s}} } &\leq
		\sum_{l=1}^{n} \left \Vert \Lambda_{s} \left( \hat G_{s}^{jl} \right) \right \Vert_{
		\mathcal{B} \left( W^{1}_{\nu_{s}}, W^{1}_{\nu_{s}} \right)}.
	\end{align}
	\begin{proof}
		It suffices to observe that 
		\begin{align*}
			\Pi^{j}_{P} A  \Pi^{ j  \mathcal{K}^{u}_{c}}_{P}(y) =
			\sum_{l=1}^{n} \Lambda_{u} \left( \hat G_{p}^{jl} \right) \left[ \tilde p \right]_{l}
		\end{align*}
		by \eqref{eq:Pi_parm}, \eqref{eq:Dc_taylor} and the definition of the approximate inverse. 
	\end{proof}
\end{lemma}
The norms in the right-hand sides of~\eqref{eq:Lambda_u}  and~\eqref{eq:Lambda_s} are controlled by using Lemma~\ref{lemma:Lambda}.

\paragraph{First order coefficients of $Z_{\Pi}(r)$}
We are now ready to construct the first order
terms of the (quadratic) polynomials $Z_{\Pi}(r)$ for $\Pi \in \mathbb{P}$. For this purpose, we first
introduce some additional notation. We will denote the bounds 
in Lemma \ref{lemma:H}, which measure the quality of the approximate derivative and inverse, by $H_{\Pi}$. 
The bounds in \eqref{eq:b_t0}, \eqref{eq:b_tm} and \eqref{eq:b_ij0}, which are related to the boundary conditions,
will be denoted by $Z^{1,\tilde \jmath t_{0}}_{\Pi}$, $Z^{1,\tilde \jmath t_{m}}_{\Pi}$ and $Z^{1,\tilde \imath \tilde \jmath 0 }_{\Pi}$,
respectively. The bounds in Lemmas \ref{lemma:Z_ij} and
\ref{lemma:Z_ijij}, which are related to the differential equation,
will be denoted by $Z^{1,\tilde \imath \tilde \jmath \left[1:N_{\tilde \imath -1}\right]}_{\Pi}$
and $Z^{1,\tilde \imath \tilde \jmath \NN}_{\Pi^{ij}_{a}}$, respectively. Finally, the bounds
in \eqref{eq:Lambda_u} and \eqref{eq:Lambda_s}, which are related to the invariance
equation for the charts on the (un)stable manifolds, will be denoted by $Z^{1,j}_{P}$ and $Z^{1,j}_{Q}$, respectively. 

With the above notation in place, we define 
\begin{align*}
		Z^{1}_{\Pi} :&= H_{\Pi} + \sum_{\tilde \jmath =1}^{n}
		\left( Z^{1,\tilde \jmath t_{0}}_{\Pi} + Z^{1,\tilde \jmath t_{m}}_{\Pi} + 
		\sum_{\tilde \imath=2}^{m} Z^{1,\tilde \imath \tilde \jmath 0}_{\Pi} + 
		\sum_{\tilde \imath=1}^{m} Z^{1,\tilde \imath \tilde \jmath \left[1:N_{\tilde \imath -1}\right]}_{\Pi} 
		\right), 
\end{align*}
for $\Pi \in \left \{ \Pi_{\theta}^{j}, \Pi_{\phi}^{j}, \Pi^{j}_{\lambda^{u}}, \Pi^{j}_{\lambda_{s}} \right \}$, and 
\begin{align*}
		Z^{1}_{\Pi^{ij}_{a}} :&= H_{\Pi^{ij}_{a}} + \sum_{\tilde \jmath =1}^{n}
		\left( Z^{1,\tilde \jmath t_{0}}_{\Pi^{ij}_{a}} + Z^{1,\tilde \jmath t_{m}}_{\Pi^{ij}_{a}} + 
		\sum_{\tilde \imath=2}^{m} Z^{1,\tilde \imath \tilde \jmath 0}_{\Pi^{ij}_{a}} + 
		\sum_{\tilde \imath=1}^{m} Z^{1,\tilde \imath \tilde \jmath \NN}_{\Pi^{ij}_{a}} 
		\right), \quad 1 \leq i \leq m, \ 1 \leq j \leq n, \\[2ex]
		Z^{1}_{\Pi^{j}_{P}} :&= H_{\Pi^{j}_{P}} + 
		\sum_{\tilde \jmath=1}^{n} \left( Z^{1,\tilde \jmath t_{0}}_{\Pi^{j}_{P}} + Z^{1,\tilde \jmath t_{m}}_{\Pi^{j}_{P}} + 
		\sum_{\tilde \imath=2}^{m} Z^{1,\tilde \imath \tilde \jmath 0}_{\Pi^{j}_{P}} + 
		\sum_{\tilde \imath=1}^{m} Z^{1,\tilde \imath \tilde \jmath \left[1:N_{\tilde \imath -1}\right]}_{\Pi^{j}_{P}} \right)
		+ Z^{1,j}_{P}, \quad 1 \leq j \leq n, \\[2ex]
		Z^{1}_{\Pi^{j}_{Q}} :&= H_{\Pi^{j}_{Q}} + 
		\sum_{\tilde \jmath=1}^{n} \left( Z^{1,\tilde \jmath t_{0}}_{\Pi^{j}_{Q}} + Z^{1,\tilde \jmath t_{m}}_{\Pi^{j}_{Q}} + 
		\sum_{\tilde \imath=2}^{m} Z^{1,\tilde \imath \tilde \jmath 0}_{\Pi^{j}_{Q}} + 
		\sum_{\tilde \imath=1}^{m} Z^{1,\tilde \imath \tilde \jmath \left[1:N_{\tilde \imath -1}\right]}_{\Pi^{j}_{Q}} \right)
		+ Z^{1,j}_{Q}, \quad 1 \leq j \leq n. 
\end{align*}

\subsubsection{Second order bounds}
\label{sec:Z2}
In this section we compute bounds for 
\begin{align}
		\label{eq:DF_perturb}
		A \left( DF \left( \hat x + rv \right) - DF \left( \hat x \right) \right)h.
\end{align}
We will compute the desired bounds by projecting \eqref{eq:DF_perturb} onto 
the relevant subspaces of $\mathcal{X}_{\nu}$ as in the previous section.
We start with the observation that 
\begin{align*}
		A \left( DF \left( \hat x + rv \right) - DF \left( \hat x \right) \right)h = 
		\int_{0}^{1} A D^{2}F \left( \hat x + \tau r v \right) \left[v,h\right] \mbox{d} \tau \ r
\end{align*}
by the (generalized) Mean Value Theorem. For notational convenience, we shall write
$y \left( \tau \right) = D^{2}F \left( \hat x + \tau r v \right) \left[v,h\right]$ and $v = \left( \dot{\theta}, \dot{\phi}, \dot{\lambda}^{u}, \dot{\lambda}^{s}, \dot{a}, \dot{p}, \dot{q} \right)$.
Furthermore, we will denote the max-norm on both $\CC^{n_{u}}$ and $\CC^{n_{s}}$ by $\left \Vert \cdot \right \Vert_{\infty}$.

Observe that $\Pi \left (y \left( \tau \right) \right) = 0$ for $\Pi \in \left \{ \Pi_{\hat p_{1}}, \Pi_{\hat q_{1}} , \Pi_{\eta}  \right \}$,
since the equations associated to these projections are linear. Furthermore, a straightforward computation shows that 
\begin{align} 
		\label{eq:Pi_t0second}
		& \Pi^{\tilde \jmath}_{t_{0}}  \left( y \left( \tau \right) \right) = 
		- \sum_{l=1}^{n_{u}} \sum_{k \in \NN_{0}^{n_{u}}}
		k_{l} \left( \hat \theta + \tau r \dot{\theta} \right) ^{k-e_{l}} 
		\left[ \tilde p_{k} \dot{\theta}_{l} +  \dot{p}_{k} \tilde \theta_{l}  \right]_{\tilde \jmath} \nonumber \\[2ex] & \quad - 
		\sum_{i,l=1}^{n_{u}}  \sum_{k \in \NN_{0}^{n_{u}}}
		k_{l} \left( k_{i} - \delta_{il} \right) \left[ \hat p_{k} + \tau r \dot{p}_{k} \right]_{\tilde \jmath} 
		\left( \hat \theta + \tau r \dot{\theta} \right) ^{k-e_{l}-e_{i}} \tilde \theta_{l} \dot{\theta}_{i}
\end{align}
for $1 \leq \tilde \jmath \leq n$. 
An analogous formula holds for $\Pi^{\tilde \jmath}_{t_{m}} \left (y \left( \tau \right) \right)$.
Next, let $2 \leq i \leq m$, then 
\begin{align}
		& \left( \Pi_{a}^{i}  \left (y \left( \tau \right) \right) \right)_{k} \\[2ex] & \quad = 
		\begin{cases}
				\bold{0}_{n}, & k = 0, \\[2ex]
				- \dfrac{L \left( t_{i} - t_{i-1} \right) }{4} 
				\left( D^{2}c_{k-1} \left( \hat a^{i} + \tau r \dot{a}^{i} \right) \left[ \dot{a}^{i}, \tilde a^{i} \right] 
				-  D^{2}c_{k+1} \left( \hat a^{i} + \tau r \dot{a}^{i} \right) \left[ \dot{a}^{i}, \tilde a^{i} \right] \right), & k \in \NN.
		\end{cases} 
		\label{eq:D2c}
\end{align}
This formula is also valid for $i=1$ and $k \in \NN$, but in this case there is no component to consider for $k=0$. Finally, 
observe that  
\begin{align}
		& \left( \Pi_{P}  \left (y \left( \tau \right) \right) \right)_{k} = \nonumber \\[2ex] & \qquad 
		\begin{cases}
				D^{2}g \left( \hat p_{0} + \tau r \dot{p}_{0} \right) \left[ \dot{p}_{0}, \tilde p_{0} \right], & k = 0, \\[2ex]
				D^{3}g \left( \hat p_{0} + \tau r \dot{p}_{0} \right) \left[ \dot{p}_{0}, \tilde p_{0}, \hat p_{e_{i}} + \tau r \dot{p}_{e_{i}} \right] + 
				D^{2} g \left( \hat p_{0} + \tau r \dot{p}_{0} \right)  \left[\tilde p_{0}, \dot{p}_{e_{i}} \right] \\[2ex] \quad + \
				D^{2} g \left( \hat p_{0} + \tau r \dot{p}_{0} \right)  \left[ \dot{p}_{0}, \tilde p_{e_{i}} \right] - 
				\left( \tilde \lambda^{u}_{i} \dot{p}_{e_{i}} + \dot{\lambda}^{u}_{i} \tilde p_{e_{i}}  \right), & k = e_{i}, \ 1 \leq i \leq n_{u}, \\[2ex]
				\left \langle \tilde \lambda^{u}, k \right \rangle \dot{p}_{k} + \left \langle \dot{\lambda}^{u}, k \right \rangle \tilde p_{k}
				-D^{2}C_{k} \left( \hat p + \tau r \dot{p} \right) \left[ \dot{p}, \tilde p \right], & \left \vert k \right \vert \geq 2.
		\end{cases}
		\label{eq:Z2_Pi_Taylor}
\end{align}
The formula for $\Pi_{Q} \left (y \left( \tau \right) \right)$ is analogous. 

Altogether, the above formulae give rise to the decomposition
\begin{align*}
		 y \left( \tau \right) = \sum_{\tilde \jmath = 1}^{n} & \left( 
		\Pi^{\tilde \jmath}_{t_{0}} \left (y \left( \tau \right) \right) + 
		\Pi_{t_{m}}^{\tilde \jmath} \left (y \left( \tau \right) \right)+  
		\sum_{\tilde \imath = 1}^{m} \Pi_{a}^{\tilde \imath \tilde \jmath \NN} \left (y \left( \tau \right) \right)
		+ \ \Pi^{\tilde \jmath}_{P} \left (y \left( \tau \right) \right) + \Pi^{\tilde \jmath}_{Q} \left (y \left( \tau \right) \right)
		\right).
\end{align*}
We will now follow the same strategy as in the previous section and compute bounds for \eqref{eq:DF_perturb}
by individually composing each term in the above decomposition with $A$ and analyzing the associated projections into the domain.

\paragraph{Boundary conditions}
We start by considering the terms associated to 
\begin{align}
		\label{eq:Z2_boundary}
		\int_{0}^{1} A \Pi  \left (y \left( \tau \right) \right) \mbox{d} \tau, 
		\quad \Pi \in \left \{ \Pi_{t_{0}}^{\tilde \jmath}, \Pi_{t_{m}}^{\tilde \jmath} \right \},
\end{align}
which are related to the boundary conditions. To compute the desired bounds, we first analyze 
the two series in \eqref{eq:Pi_t0second}.
\begin{lemma}
		\label{lemma:series1}
		Suppose $\theta \in \normalfont \mbox{int} \ \mathbb{B}_{\nu_{u}} \setminus \{0\}$ and let
		$\zeta \in \CC^{n_{u}}$ be such that $\left \Vert \zeta \right \Vert_{\infty} \leq 1$. 
		Define a linear map $\varphi_{u,1}: W^{1}_{\nu_{u}} \rightarrow \CC$ by 
		\begin{align*}
				\varphi_{u,1}(p) := \sum_{l=1}^{n_{u}} \sum_{k \in \NN_{0}^{n_{u}}} 
				k_{l} p_{k} \theta^{k-e_{l}} \zeta_{l}.
		\end{align*}
		Then $\varphi_{u,1} \in \left( W^{1}_{\nu_{u}} \right)^{\ast}$ and 
		\begin{align*}
				\left \Vert \varphi_{u,1} \right \Vert_{\mathcal{B} \left( W^{1}_{\nu_{u}}, \CC \right)} \leq \Phi_{u,1} \left( \theta \right) :=
				\begin{cases}
				\nu_{u}^{-1}, & \log\left( \dfrac{\left \Vert \theta \right \Vert_{\infty}}{ \nu_{u}} \right) \leq -1, \\[2ex]
				-\left( e \left \Vert \theta \right \Vert_{\infty} \log \left( \dfrac{\left \Vert \theta \right \Vert_{\infty}}{ \nu_{u}} \right) \right)^{-1},
				& \normalfont \text{otherwise}. 
				\end{cases}
		\end{align*}
		\begin{proof}
				It is clear that $\varphi_{u,1}$ is bounded for any $\theta \in \normalfont \mbox{int} \ \mathbb{B}_{\nu_{u}}$, since 
				\begin{align}
					\label{eq:O(ptheta)}
					\left \vert p_{k} \theta^{k-e_{l}} \right \vert \leq  
					\frac{ \left \Vert \theta \right \Vert_{\infty}^{ \left \vert k \right \vert - 1} }{ \nu_{u}^{ \left \vert k \right \vert } }
					\left \Vert p \right \Vert_{\nu_{u}}, 
					\quad k \in \NN_{0}^{n_{u}}.
				\end{align}
				We will use this observation to compute a bound for the operator norm.
				Namely, let $j \in \NN_{0}^{n_{u}}$ be arbitrary, then 
				\begin{align*}
						\left \vert \varphi_{u,1} \left( \boldsymbol{\xi}_{j} \right) \right \vert \leq  
						\sum_{l=1}^{n_{u}} j_{l} 
						\left \vert \theta \right \vert^{j - e_{l}} \nu_{u}^{-\left \vert j \right \vert} \leq 
						 \left \Vert \theta \right \Vert_{\infty}^{-1}
						\left \vert j \right \vert 
						\left( \frac{ \left \Vert \theta \right \Vert_{\infty} }{ \nu_{u}} \right)^{ \left \vert j \right \vert },
				\end{align*}
				where $\left( \boldsymbol{\xi}_{l} \right)_{l \in \NN_{0}^{n_{u}}}$ are the corner points introduced
				in Definition \ref{def:cornerpoints_array}. In particular, observe that 
				$\left \Vert \varphi_{u,1} \right \Vert_{\mathcal{B} \left( W^{1}_{\nu_{u}}, \CC \right) } = 
				  \sup_{j \in \NN_{0}^{n_{u}} \setminus {\{0\}} } \left \vert \varphi_{u,1} \left( \boldsymbol{\xi}_{j} \right)\right \vert$
				by Proposition \ref{prop:operator_norm_mult}, since $\varphi_{u,1} \left( \boldsymbol{\xi}_{0} \right) = 0$. 
				Next, note that the map 
				$x \mapsto x \rho^{x}$, where $\rho \in (0,1)$, is strictly increasing on 
				$\left[0, \dfrac{-1}{\log \rho} \right]$, strictly decreasing on 
				$\left[\dfrac{-1}{\log \rho}, \infty \right)$ and has a global maximum on $\left[0,\infty\right)$ at
				$x = \dfrac{-1}{ \log{\rho} }$. Therefore, since 
				$\left \vert j \right \vert \geq 1$ for $j \in \NN_{0}^{n_{u}} \setminus {\{0\}}$
				and $\left \Vert \theta \right \Vert_{\infty} < \nu_{u}$, it follows that 
				$
						\sup_{j \in \NN_{0}^{n_{u}} \setminus {\{0\}} } \left \vert \varphi_{u,1} \left( \boldsymbol{\xi}_{j} \right) \right \vert \leq  
						\Phi_{u,1} \left( \theta \right),
				$
				which proves the result. 
		\end{proof}
\end{lemma}
\begin{remark}
	The analogues of $\varphi_{u,1}$ and $\Phi_{u,1}$ in the context of the stable manifold are 
	defined similarly and are
	denoted by $\varphi_{s,1}$ and $\Phi_{s,1} \left( \phi \right)$, respectively. 
\end{remark}
\begin{remark}
	In practice, $\theta$ is an interval enclosure of $\hat \theta$ (see Lemma \ref{lemma:Z2_t0tm}). 
	Recall that $\hat \theta$ corresponds to 
	a (numerically obtained) coordinate on the chart of the unstable manifold at which the connecting orbit $u$ starts. 
	This is why we may assume that $\theta \not = 0$. 
\end{remark}

The latter result provides a way to bound the first term in \eqref{eq:Pi_t0second}. To bound the 
second term \eqref{eq:Pi_t0second}, we perform a similar analysis. 
\begin{lemma}
		\label{lemma:series2}
		Suppose $\theta \in \normalfont \mbox{int} \ \mathbb{B}_{\nu_{u}} \setminus \{0\}$ and define 
		a linear map $\varphi_{u,2} : W^{1}_{\nu_{u}} \rightarrow \CC$ by 
		\begin{align*}
				\varphi_{u,2}(p) := 
				\sum_{i,l=1}^{n_{u}} \sum_{ k \in \mathcal{K}_{c}^{u} }
				k_{l} \left( k_{i} - \delta_{il} \right) p_{k} \theta^{k-e_{l}-e_{i}} 
				\tilde \theta_{l}\dot{\theta}_{i}.
		\end{align*}
		Then $\varphi_{u,2} \in \left( W^{1}_{\nu_{u}} \right)^{\ast}$ and 
		$
				\left \Vert \varphi_{u,2} \right \Vert_{\mathcal{B} \left( W^{1}_{\nu_{u}}, \CC \right)} 
				\leq \Phi_{u,2} \left( \theta \right),
		$
		where
		\begin{align*}
				\Phi_{u,2} \left( \theta \right) := 
				\left \Vert \theta \right \Vert_{\infty}^{-2}
				\begin{cases}
						\normalfont
						\left(  K^{u}_{\text{min}} +1  \right)^{2} \left( \dfrac{ \left \Vert \theta \right \Vert_{\infty} }{ \nu_{u}} \right)^{K^{u}_{\text{min}}+1 },
						& \normalfont K^{u}_{\text{min}} +1 \geq -2 \log \left( \dfrac{ \left \Vert \theta \right \Vert_{\infty} }{ \nu_{u}} \right)^{-1}, \\[2ex]
						4 \left( e \log \left( \dfrac{ \left \Vert \theta \right \Vert_{\infty} }{ \nu_{u}} \right) \right)^{-2}, & \normalfont \mbox{otherwise},
				\end{cases}				
		\end{align*}
		and $\normalfont K^{u}_{\text{min}} := \min_{1 \leq i \leq n_{u}} K^{u}_{i}$.
		\begin{proof}
				The boundedness of $\varphi_{u,2}$ follows directly from the observation in \eqref{eq:O(ptheta)} and the assumption that 
				$\theta \in \normalfont \mbox{int} \ \mathbb{B}_{\nu_{u}}$.
				To compute a bound for the operator norm, first observe that  $\varphi_{u,2} \left( \boldsymbol{\xi}_{j} \right) = 0$ for $j \in \mathcal{K}^{u}$.
				Hence 
				$
					\left \Vert \varphi_{u,2} \right \Vert_{\mathcal{B} \left( W^{1}_{\nu_{u}}, \CC \right)}  = 
					\sup_{j \in \mathcal{K}^{u}_{c}} \left \vert \varphi_{u,2} \left( \boldsymbol{\xi}_{j} \right) \right \vert 
				$
				by Proposition \ref{prop:operator_norm_mult}. Furthermore, a straightforward computation shows that 
				\begin{align*}
						\left \vert \varphi_{u,2} \left( \boldsymbol{\xi}_{j} \right) \right \vert &\leq 
						\sum_{i,l=1}^{n_{u}} j_{l} \left( j_{i} - \delta_{il} \right) 
						\left \vert \theta \right \vert ^{j-e_{l}-e_{i}}  \nu_{u}^{- \left \vert j \right \vert} \\[2ex] &\leq 
						\left \Vert \theta \right \Vert_{\infty}^{-2}
						\left( \frac{ \left \Vert \theta \right \Vert_{\infty} }{ \nu_{u} } \right)^{\left \vert j \right \vert}
						\sum_{i,l=1}^{n_{u}} j_{l} \left( j_{i} - \delta_{il} \right) \\[2ex] &\leq 
						\left \Vert \theta \right \Vert_{\infty}^{-2} \left \vert j \right \vert^{2}
						\left( \frac{ \left \Vert \theta \right \Vert_{\infty} }{ \nu_{u} } \right)^{\left \vert j \right \vert},
				\end{align*}	
				for any $j \in \mathcal{K}^{u}_{c}$, where we used that $\left \Vert \tilde \theta \right \Vert_{\infty}, 
				\left \Vert \dot{\theta} \right \Vert_{\infty} \leq 1$. 				
				Now, note that the map $x \mapsto x^{2} \rho^{x}$, where $\rho \in (0,1)$, is strictly increasing on 
				$\left[0, \dfrac{-2}{\log{\rho}} \right]$, strictly decreasing on $\left [\dfrac{-2}{\log{\rho}}, \infty \right)$
				and has a global maximum on $\left[0, \infty \right)$ at $x = \dfrac{-2}{\log{\rho}}$. Therefore, since
				$\left \vert j \right \vert \geq K^{u}_{\text{min}}+1$ for $j \in \mathcal{K}^{u}_{c}$, 
				it follows that 
				$
					\sup_{j \in \mathcal{K}^{u}_{c}} \left \vert \varphi_{u,2} \left( \boldsymbol{\xi}_{j} \right) \right \vert
					\leq \Phi_{u,2} \left( \theta \right),
				$
				which proves the result. 
		\end{proof}
\end{lemma}
\begin{remark}
	As before, the analogues of $\varphi_{u,2}$ and $\Phi_{u,2}$ in the context of
	the stable manifold are denoted by $\varphi_{s,2}$ and $\Phi_{s,2} \left( \phi \right)$, respectively. 
\end{remark}

We are now ready to compute bounds for \eqref{eq:Z2_boundary}. The computation of these bounds consists of a mixture of interval 
analysis on the computer and ordinary estimates derived with ``pen and paper''. More precisely, as mentioned at the beginning of this 
paper, in order to manage the rounding errors on the computer, all the bounds in this paper are computed with \emph{interval arithmetic}. 
Roughly speaking, this means that all the elementary operations on floating point numbers are replaced by operations on intervals with endpoints representable on a computer. 
In this way, one can compute rigorous bounds (and hence verify inequalities)  with the aid of a computer. 

Now, in order to compute bounds for \eqref{eq:Z2_boundary},  we first define interval enclosures for 
$\hat \theta$, $\hat p$, $\hat \phi$ and $\hat q$.
Let $r^{\ast}>0$ be an upper bound for the radius $r$ and set 
\begin{align*}
	\boldsymbol{\hat \theta} :&= \prod_{j=1}^{n_{u}} \left[ \hat \theta_{j} - r^{\ast}, \hat \theta_{j} + r^{\ast} \right], \quad
	\boldsymbol{ \hat p }_{k} := \prod_{j=1}^{n} \left[ \left[\hat p_{k} \right]_{j} 
	- \frac{r^{\ast}}{\nu_{u}^{\left \vert k \right \vert}}, \ \left[ \hat p_{k} \right]_{j} 
	+  \frac{r^{\ast}}{\nu_{u}^{\left \vert k \right \vert}} \right], \quad k \in \mathcal{K}^{u}, \\[2ex]
	\boldsymbol{\hat \phi} :&= \prod_{j=1}^{n_{s}} \left[ \hat \phi_{j} - r^{\ast}, \hat \phi_{j} + r^{\ast} \right], \quad
	\boldsymbol{ \hat q }_{k} := \prod_{j=1}^{n} \left[ \left[\hat q_{k} \right]_{j} 
	- \frac{r^{\ast}}{\nu_{s}^{\left \vert k \right \vert}}, \ \left[ \hat q_{k} \right]_{j} 
	+  \frac{r^{\ast}}{\nu_{s}^{\left \vert k \right \vert}} \right], \quad k \in \mathcal{K}^{s}. 
\end{align*}
We require that $r^{\ast}$ is sufficiently small and that $\nu_{u}, \nu_{s}$ are sufficiently
large so that $\boldsymbol{\hat \theta} \subset \mbox{int} \ \mathbb{B}_{\nu_{u}}$ and $\boldsymbol{\hat \phi} \subset \mbox{int} \ \mathbb{B}_{\nu_{s}}$.
\begin{remark}
	Strictly speaking, the endpoints of the above intervals should be floating point numbers 
	so that we can perform rigorous computations on a computer. 
	In practice, this amounts to computing slightly larger interval enclosures for $\hat \theta$, $\hat p$, $\hat \phi$ 
	and $\hat q$ (compared to the ones above). To avoid clutter in the notation,
	however, we have chosen to ignore this rather technical (but easily solved) issue.
\end{remark}
\begin{lemma}
		\label{lemma:Z2_t0tm}
		Let $1 \leq \tilde \jmath \leq n$, $0<r \leq r^{\ast}$ and $\Pi \in \mathbb{P}$, then
		\begin{align}
				\label{eq:Z2_t0}
				\sup_{\tau \in [0,1]} \left \Vert \Pi A \Pi^{\tilde \jmath}_{t_{0}} \left( y \left( \tau \right) \right) \right \Vert_{\Pi \left( \mathcal{X}_{\nu} 					\right)} &\leq \
				\left \Vert \Pi A_{NK} \Pi^{\tilde \jmath}_{t_{0}} \right \Vert_{\Pi \left( \mathcal{X}_{\nu} \right)} \sup \boldsymbol{\beta}_{u}^{\tilde 					\jmath}, \\[2ex]
				\label{eq:Z2_tm}
				\sup_{\tau \in [0,1]}  \left \Vert \Pi A \Pi^{\tilde \jmath}_{t_{m}} \left( y \left( \tau \right) \right) \right \Vert_{\Pi \left( \mathcal{X}_{\nu} 				\right)} &\leq \
				\left \Vert \Pi A_{NK} \Pi^{\tilde \jmath}_{t_{m}} \right \Vert_{\Pi \left( \mathcal{X}_{\nu} \right)} \sup \boldsymbol{\beta}_{s}^{\tilde 					\jmath},
		\end{align}
		where 
		\begin{align*}
			\boldsymbol{\beta}_{u}^{\tilde \jmath} :&= 2 \Phi_{u,1} \left( \boldsymbol{ \hat \theta } \right) 
			+  r^{\ast} \Phi_{u,2} \left( \boldsymbol{ \hat \theta } \right) 
			+ \sum_{i,l=1}^{n_{u}}  \sum_{ k \in \mathcal{K}^{u} }
			k_{l} \left( k_{i} - \delta_{il} \right) \left \vert \left[ \boldsymbol{ \hat p }_{k} \right]_{\tilde \jmath} \right \vert
			\left \vert \boldsymbol{\hat \theta} \right \vert^{k-e_{l}-e_{i}} 
			, \\[2ex]
			\boldsymbol{\beta}_{s}^{\tilde \jmath} :&= 2 \Phi_{s,1} \left( \boldsymbol{ \hat \phi } \right) 
			+ r^{\ast} \Phi_{s,2} \left( \boldsymbol{ \hat \phi } \right)
			+ \sum_{i,l=1}^{n_{s}} \sum_{ k \in \mathcal{K}^{s} }
			k_{l} \left( k_{i} - \delta_{il} \right) \left \vert \left[ \boldsymbol{ \hat q }_{k} \right]_{\tilde \jmath} \right \vert
			\left \vert \boldsymbol{\hat \phi} \right \vert^{k-e_{l}-e_{i}}.
		\end{align*}
		\begin{proof}
			Let $\tau \in [0,1]$ and $1 \leq \tilde \jmath \leq n$ be arbitrary and use 
			Lemma \ref{lemma:series1} to see 
			that the first term in \eqref{eq:Pi_t0second} is bounded by $2 \Phi_{u,1} \left( \hat \theta + \tau r \dot{\theta} \right)$.
			To bound the second term in \eqref{eq:Pi_t0second}, we first split it into two series; one
			over $\mathcal{K}^{u}$ and one over $\mathcal{K}^{u}_{c}$. We then use
			Lemma \ref{lemma:series2} and the fact that $\hat p_{k} = 0$ for $k \in \mathcal{K}^{u}_{c}$ to 
			estimate 
			\begin{align*}
				& \left \vert \sum_{i,l=1}^{n_{u}}  \sum_{k \in \NN_{0}^{n_{u}}}
				k_{l} \left( k_{i} - \delta_{il} \right) \left[ \hat p_{k} + \tau r \dot{p}_{k} \right]_{j} 
				\left( \hat \theta + \tau r \dot{\theta} \right) ^{k-e_{l}-e_{i}} \tilde \theta_{l} \dot{ \theta}_{i} \right \vert \\[2ex] & \quad \leq 				
				\left \vert \sum_{i,l=1}^{n_{u}}  \sum_{k \in \mathcal{K}^{u}}
				k_{l} \left( k_{i} - \delta_{il} \right) \left[ \hat p_{k} + \tau r \dot{p}_{k} \right]_{\tilde \jmath}
				\left( \hat \theta + \tau r \dot{\theta} \right) ^{k-e_{l}-e_{i}} \tilde \theta_{l} \dot{ \theta}_{i} \right \vert + 
				r \Phi_{u,2} \left( \hat \theta + \tau r \dot{\theta} \right).
			\end{align*}
			Altogether, we conclude that 
			\begin{align}
				& \left \vert \Pi_{t_{0}}^{\tilde \jmath} \left( y \left( \tau \right) \right) \right \vert \leq 
				2 \Phi_{u,1} \left( \hat \theta + \tau r \dot{\theta} \right) + r^{\ast} 
				\Phi_{u,2} \left( \hat \theta + \tau r \dot{\theta} \right) \nonumber \\[2ex] & 				\quad 
				+ \sum_{i,l=1}^{n_{u}}  \sum_{k \in \mathcal{K}^{u}}
				k_{l} \left( k_{i} - \delta_{il} \right) \left \vert \left[ \hat p_{k} + \tau r \dot{p}_{k} \right]_{\tilde \jmath} \right \vert 
				\left \vert \hat \theta + \tau r \dot{\theta} \right \vert ^{k-e_{l}-e_{i}},
				\label{eq:beta_u}
			\end{align}
			where we used the assumption that $r \leq r^{\ast}$. 
			Finally, observe that $\hat \theta + \tau r \dot{\theta} \in \boldsymbol{\hat \theta}$ and
			$\hat p_{k} + \tau r \dot{p}_{k} \in \boldsymbol{\hat p}_{k}$, since 
			$\Vert \dot{\theta}  \Vert_{\infty}, \ \left \Vert \dot{p} \right \Vert_{ W^{1}_{\nu_{u},n}} \leq 1$.
			Hence \eqref{eq:beta_u} is contained in $\boldsymbol{\beta}^{\tilde \jmath}_{u}$ for all $\tau \in [0,1]$, which proves the result. 
		\end{proof}
\end{lemma}
\begin{remark}
	Observe that the computation of the bound in this lemma is finite, since $\mathcal{K}^{u}$ is a finite set of multi-indices. 
\end{remark}

\paragraph{Chebyshev coefficients}
Next, we consider the terms associated to 
\begin{align}
		\label{eq:Z2_cheb}
		\int_{0}^{1} A \Pi^{\tilde \imath \tilde \jmath \NN}_{a} \left (y \left( \tau \right) \right) \mbox{d} \tau, 
		\quad 1 \leq \tilde \imath \leq m, \ 1 \leq \tilde \jmath \leq n, 
\end{align}
which are related to the equations for the Chebyshev coefficients. 
We start by computing a bound for \eqref{eq:D2c}. Since we will have to perform a similar analysis for the Taylor coefficients 
in the next paragraph, we first state a general result. 
\begin{lemma}
	\label{lemma:D2g}
	Suppose $\left( X, \ast \right)$ is a Banach algebra. Let $1 \leq j \leq n$
	and define $G_{j} : \bigoplus_{l=1}^{n} X \rightarrow X$ and $\tilde g_{j} : \RR^{n} \rightarrow \RR$
	by 
	\begin{align*}
		G_{j}(x) := \sum_{\alpha \in \mathcal{A}} \left[ g_{\alpha} \right]_{j} x^{\alpha},  \quad
		\tilde g_{j}(x) := \sum_{\alpha \in \mathcal{A} } \left \vert \left[ g_{\alpha} \right]_{j} \right \vert x^{\alpha},
	\end{align*}
	where $\left \{ g_{\alpha} : \alpha \in \mathcal{A} \right \} \subset \RR^{n}$ are the coefficients of
	$g$ in the monomial basis. Then 
	\begin{align*}
		\left \Vert D^{2} G_{j} \left( x + \tau r z \right) \left[y,z\right] \right \Vert \leq
		D^{2} \tilde g_{j} \left( \left \Vert x_{1} \right \Vert_{X} + r^{\ast}, \ldots, \left \Vert x_{n} \right \Vert_{X} + r^{\ast} \right) \left[ \bold{1}_{n}, \bold{1}_{n} \right],
	\end{align*}
	for any $x = \left( x_{l} \right)_{l=1}^{n}, y = \left( y_{l} \right)_{l=1}^{n}, z = \left( z_{l} \right)_{l=1}^{n} \in \bigoplus_{l=1}^{n} X$ such that 
	$\left \Vert y \right \Vert, \left \Vert z \right \Vert \leq 1$, $\tau \in [0,1]$ and $0<r \leq r^{\ast}$. Here $\left \Vert \cdot \right \Vert$ denotes the max-norm on 
	$\bigoplus_{l=1}^{n} X$, and $\bold{1}_{n}$ denotes the vector of length $n$ containing $1$-s. 
	\begin{proof}
		Note that we may use the ``usual'' rules of calculus to differentiate polynomials on $\bigoplus_{l=1}^{n} X$, 
		since $\left( X, \ast \right)$ is a Banach algebra. In particular, a straightforward
		computation shows that 
		\begin{align*}
			D^{2}G_{j} \left( x + \tau r z \right) \left[y,z \right] = 
			\sum_{i,l=1}^{n} \sum_{ \alpha \in \mathcal{A} } \alpha_{i} \left( \alpha_{l} - \delta_{il} \right) \left[ g_{\alpha} \right]_{j}
			\left( x + \tau rz \right)^{\alpha - e_{i} - e_{l}} \ast y_{i} \ast z_{l}. 
		\end{align*}
		Hence, by the Banach algebra estimate, 
		\begin{align*}
			\left \Vert D^{2}G_{j} \left( x + \tau r z \right) \left[y,z\right] \right \Vert &\leq 
			\sum_{i,l=1}^{n} \sum_{ \alpha \in \mathcal{A} }
			\alpha_{i} \left( \alpha_{l} - \delta_{il} \right) \left \vert \left[ g_{\alpha} \right]_{j} \right \vert 
			\prod_{k=1}^{n} \left( \left \Vert x_{k} \right \Vert_{X} + r^{\ast} \right)^{ \left( \alpha - e_{i} - e_{l} \right)_{k} } \\[2ex] &= 
			D^{2} \tilde g_{j} \left( \left \Vert x_{1} \right \Vert_{X} + r^{\ast}, \ldots, \left \Vert x_{n} \right \Vert_{X} + r^{\ast} \right) \left[ \bold{1}_{n}, \bold{1}_{n} \right],
		\end{align*}
		since $\left \Vert y \right \Vert, \left \Vert z \right \Vert \leq 1$, $\tau \in [0,1]$ and $0<r \leq r^{\ast}$.
	\end{proof}
\end{lemma}
\begin{remark}
	Note that the convolution mappings $\left[C\right]_{j}$ and $\left[c\right]_{j}$
	are of the form $G_{j}$. 
\end{remark}

Next, we use the above result to compute bounds for \eqref{eq:Z2_cheb}. The key observation is 
stated in the next lemma.
\begin{lemma}
	\label{lemma:Dc2_bound}
	Let $1 \leq \tilde \imath \leq m$, $1 \leq \tilde \jmath \leq n$ and $0<r \leq r^{\ast}$, then 
	\begin{align*}
		\sup_{\tau \in [0,1]} \left \Vert \Pi^{ \tilde \imath \tilde \jmath \NN}_{a} \left( y \left( \tau \right) \right) \right \Vert_{\nu_{\tilde \imath}} \leq 
		\dfrac{ L \nu_{\tilde \imath} \left( t_{\tilde \imath} - t_{\tilde \imath-1} \right) }{ 2 }
		D^{2} \tilde g_{\tilde \jmath} \left( 
			\left \Vert \left[ \hat a^{\tilde \imath} \right]_{1} \right \Vert_{\nu_{\tilde \imath}} + r^{\ast}, \ldots, 
			\left \Vert \left[ \hat a^{\tilde \imath} \right]_{n} \right \Vert_{\nu_{\tilde \imath}} + r^{\ast}
			\right) \left[ \boldsymbol{1}_{n}, \boldsymbol{1}_{n} \right].
	\end{align*}
	\begin{proof}
		Define linear operators
		$\sigma_{1} : \ell^{1}_{\nu_{1}} \rightarrow \ell^{1}_{\nu_{1}} / \CC$ and 
		$\left \{ \sigma_{\tilde \imath} : \ell^{1}_{\nu_{\tilde \imath}} \rightarrow  \ell^{1}_{\nu_{\tilde \imath}} \right \}_{\tilde \imath=2}^{m}$ 
		by 
		\begin{align*}
			\sigma_{1} \left (a^{1} \right) := \left( a^{1}_{k-1} - a^{1}_{k+1} \right)_{k \in \NN}, \quad 	
			\sigma_{\tilde \imath} \left(a^{\tilde \imath} \right) := 
			\begin{bmatrix}
				0 \\[2ex]
				\left( a^{\tilde \imath}_{k-1}  - a^{\tilde \imath}_{k+1} \right)_{k \in \NN} \quad 
			\end{bmatrix}.
		\end{align*}
		A direct application of Proposition \ref{prop:l1_operator_norm} shows that these operators are bounded and that 
		\begin{align}
			\label{eq:norm_sigma}
			\left \Vert \sigma_{1} \right \Vert_{\mathcal{B} \left( \ell^{1}_{\nu}, \ \ell^{1}_{\nu_{1}} / \CC \right) } = 2 \nu_{1}, \quad 
			\left \Vert \sigma_{\tilde \imath} \right \Vert_{\mathcal{B} \left( \ell^{1}_{\nu_{\tilde \imath}},  \ell^{1}_{\nu_{\tilde \imath}} \right) } = 
			2 \nu_{\tilde \imath}, \quad 2 \leq \tilde \imath \leq m.
		\end{align}
		Therefore, since 
		\begin{align*}
			\Pi^{\tilde \imath \tilde \jmath \NN}_{a} \left( y \left( \tau \right) \right) = 
			\dfrac{L\left( t_{ \tilde \imath } - t_{\tilde \imath -1} \right)}{4}
			\sigma_{\tilde \imath} D^{2} \left[ c \right]_{ \tilde \jmath} 
			\left( \hat a^{\tilde \imath} + \tau r  \dot{a}^{\tilde \imath} \right)\left[\dot{a}^{\tilde \imath}, \tilde a^{\tilde \imath} \right]
		\end{align*}
		by \eqref{eq:D2c}, the result follows directly from \eqref{eq:norm_sigma} and Lemma \ref{lemma:D2g}. 
	\end{proof}
\end{lemma}

It is now a straightforward task to compute bounds for \eqref{eq:Z2_cheb}. 
\begin{lemma}[Scalar and Taylor projections]
	\label{lemma:Z2_cheb_scal_Taylor}
	Let $1 \leq  \tilde \imath \leq m$, $1 \leq  \tilde \jmath  \leq n$, $0<r \leq r^{\ast}$ and  
	 $\Pi \in \left \{ \Pi^{j}_{\theta}, \Pi^{j}_{\phi}, \Pi_{\lambda^{u}}^{j}, \Pi^{j}_{\lambda^{s}} 
	  \Pi^{j}_{P}, \Pi^{j}_{Q} \right \}$, then 
	\begin{align*}
		& \sup_{\tau \in [0,1]} \left \Vert \Pi A \Pi^{\tilde \imath \tilde \jmath \NN}_{a} \left( y \left( \tau \right) \right) 
		\right \Vert_{\Pi \left( \mathcal{X}_{\nu} \right)}  \leq 
		\left \Vert \Pi A_{NK} \Pi^{\tilde \imath \tilde \jmath \left[1:N_{\tilde \imath}-1 \right]}_{a} 
		\right \Vert_{ \mathcal{B} \left( \ell^{1}_{\nu_{\tilde \imath} }, \Pi \left( \mathcal{X}_{\nu} \right) \right)}		
		\\[2ex] & \quad \cdot \dfrac{ L \nu_{\tilde \imath} \left( t_{\tilde \imath} - t_{\tilde \imath-1} \right) }{ 2 }
		D^{2} \tilde g_{\tilde \jmath} \left( 
			\left \Vert \left[ \hat a^{\tilde \imath} \right]_{1} \right \Vert_{\nu_{\tilde \imath}} + r^{\ast}, \ldots, 
			\left \Vert \left[ \hat a^{\tilde \imath} \right]_{n} \right \Vert_{\nu_{\tilde \imath}} + r^{\ast}
			\right) \left[ \boldsymbol{1}_{n}, \boldsymbol{1}_{n} \right].
	\end{align*}
	\begin{proof}
		It suffices to observe that 
		\begin{align*}
			\Pi A \Pi^{\tilde \imath \tilde \jmath \NN}_{a} \left( y \left( \tau \right) \right) = 
			\Pi A_{NK} \Pi^{\tilde \imath \tilde \jmath \left[1:N_{\tilde \imath}-1 \right]}_{a} 
			\left( y \left( \tau \right) \right), \quad \Pi \in \left \{ \Pi^{j}_{\theta}, \Pi^{j}_{\phi}, \Pi_{\lambda^{u}}^{j}, \Pi^{j}_{\lambda^{s}} 
	  		\Pi^{j}_{P}, \Pi^{j}_{Q} \right \},
		\end{align*}
		by construction of the approximate inverse $A$. 
		Therefore, the result follows directly from Lemma 
		\ref{lemma:Dc2_bound}.
	\end{proof}
\end{lemma}
\begin{lemma}[Chebyshev projections]
	\label{lemma:Z2_cheb_cheb}
	Let $1 \leq \tilde \imath, i \leq m$, $1 \leq \tilde \jmath, j \leq n$ and $0<r \leq r^{\ast}$, then 
	\begin{align*}
		& \sup_{\tau \in [0,1]} \left \Vert \Pi^{ij}_{a} A \Pi^{\tilde \imath \tilde \jmath \NN}_{a} \left( y \left( \tau \right) \right)
		 \right \Vert_{ \nu_{i} }
		\\[2ex] & \quad \leq \dfrac{ L \nu_{\tilde \imath} \left( t_{\tilde \imath} - t_{\tilde \imath-1} \right) }{ 2 }
		D^{2} \tilde g_{\tilde \jmath} \left( 
			\left \Vert \left[ \hat a^{\tilde \imath} \right]_{1} \right \Vert_{\nu_{\tilde \imath}} + r^{\ast}, \ldots, 
			\left \Vert \left[ \hat a^{\tilde \imath} \right]_{n} \right \Vert_{\nu_{\tilde \imath}} + r^{\ast}
			\right) \left[ \boldsymbol{1}_{n}, \boldsymbol{1}_{n} \right] \\[2ex] & \qquad 
		\cdot \begin{cases}
			\left \Vert \Pi_{a}^{ij} A_{NK} \Pi^{\tilde \imath \tilde \jmath \left[ 1 : N_{\tilde \imath}-1\right]}_{a} \right \Vert
			_{ \mathcal{B} \left( \ell^{1}_{\nu_{\tilde \imath}}, \ell^{1}_{\nu_{i}} \right) }, 
			& \left(\tilde \imath, \tilde \jmath \right) \not = \left(i,j\right), \\[4ex]
			\max \left \{
				\left \Vert \Pi_{a}^{ij} A_{NK} \Pi_{a}^{ij \left[ 1 : N_{i}-1\right]} \right \Vert
				_{ \mathcal{B} \left( \ell^{1}_{\nu_{i}}, \ell^{1}_{\nu_{i}} \right) }, \ \dfrac{1}{N_{i}}
				\right\}, 
				&\left(\tilde \imath, \tilde \jmath \right) =  \left(i,j\right) .
		\end{cases}	
	\end{align*}
	\begin{proof}
		It follows from the definition of the approximate inverse that
		\begin{align*}
			\Pi^{ij}_{a}  A \Pi^{\tilde \imath \tilde \jmath \NN}_{a} &= 
			\Pi^{ij}_{a} A_{NK} \Pi^{\tilde \imath \tilde \jmath \left[ 1: N_{\tilde \imath} - 1 \right]}_{a} + 
			\Pi^{ij}_{a} A \Pi^{\tilde \imath \tilde \jmath \left[ N_{\tilde \imath} : \infty \right)}_{a}
			\\[2ex] &= 			
			\begin{cases}
				\Pi_{a}^{ij}  A_{NK} \Pi^{\tilde \imath \tilde \jmath \left[ 1 : N_{\tilde \imath}-1\right]}_{a}, 
				& \left(\tilde \imath, \tilde \jmath \right) \not = \left(i,j\right), \\[3ex]
				\begin{bmatrix}
					\Pi^{ij}_{a} A_{NK} \Pi^{i j \left[ 1 : N_{ i}-1\right]}_{a} & & & \\[2ex]
					& \dfrac{1}{N_{i}} & & \\
					& & \dfrac{1}{N_{i}+1}  & \\
					& &  & \ddots
				\end{bmatrix}, 
				&  \left(\tilde \imath, \tilde \jmath \right) = \left(i,j\right).
			\end{cases}
		\end{align*}
		Hence
		\begin{align*}
			\left \Vert \Pi^{ij}_{a} A \Pi^{\tilde \imath \tilde \jmath \NN}_{a} \right \Vert
			_{ \mathcal{B} \left( \ell^{1}_{\nu_{\tilde \imath}}, \ell^{1}_{\nu_{i}} \right) } = 
			\begin{cases}
				\left \Vert \Pi_{a}^{ij} A_{NK} \Pi^{\tilde \imath \tilde \jmath \left[ 1 : N_{\tilde \imath}-1\right]}_{a} \right \Vert
				_{ \mathcal{B} \left( \ell^{1}_{\nu_{\tilde \imath}}, \ell^{1}_{\nu_{i}} \right) }, 
				& \left(\tilde \imath, \tilde \jmath \right) \not = \left(i,j\right), \\[4ex]
				\max \left \{
					\left \Vert \Pi_{a}^{ij} A_{NK} \Pi_{a}^{ij \left[ 1 : N_{i}-1\right]} \right \Vert
					_{ \mathcal{B} \left( \ell^{1}_{\nu_{i}}, \ell^{1}_{\nu_{i}} \right) }, \ \dfrac{1}{N_{i}}
					\right\}, 
					&\left(\tilde \imath, \tilde \jmath \right) =  \left(i,j\right),
			\end{cases}		
		\end{align*}
		by Proposition \ref{prop:l1_operator_norm}. The result now follows from Lemma \ref{lemma:Dc2_bound}. 
	\end{proof}
\end{lemma}

\paragraph{Taylor coefficients}
Finally, we consider the terms associated to 
\begin{align}
		\label{eq:Z2_Taylor}
		\int_{0}^{1} A \Pi \left (y \left( \tau \right) \right) \mbox{d} \tau, 
		\quad \Pi \in \left \{ \Pi^{ \tilde \jmath }_{P}, \Pi^{ \tilde \jmath }_{Q} \right \},
\end{align}
which are related to the equations for the Taylor coefficients of the (un)stable manifolds. 
Observe that $\Pi^{ \tilde \jmath}_{P} \left( y \left( \tau \right) \right) \not \in W^{1}_{\nu_{u}}$, 
due to the presence of the terms 
$\left \langle \tilde \lambda^{u}, k \right \rangle \dot{p}_{k} + \left \langle \dot{\lambda}^{u}, k \right \rangle \tilde p_{k}$,
see \eqref{eq:Z2_Pi_Taylor}. For this reason, in order to facilitate the analysis of \eqref{eq:Z2_Taylor}, we introduce the term 
$y^{\tilde \jmath}_{P} \left( \tau \right) := \Pi^{ \tilde \jmath}_{P}\left( y \left( \tau \right) \right) - y^{\tilde \jmath}_{P,\lambda} \left( \tau \right)$, where 
\begin{align*}
	\left( y^{\tilde \jmath}_{P,\lambda} \left( \tau \right) \right)_{k} := 
	\begin{cases}
		0, & k \in \mathcal{K}^{u}, \\[2ex]
		\left \langle \tilde \lambda^{u}, k \right \rangle \left[ \dot{p}_{k} \right]_{\tilde \jmath} + 
		\left \langle \dot{\lambda}^{u}, k \right \rangle \left[ \tilde p_{k} \right]_{\tilde \jmath}, 
		& k \not \in \mathcal{K}^{u},
	\end{cases} 
\end{align*}
i.e., $y^{\tilde \jmath}_{P} \left( \tau \right)$ is defined by removing the linear terms 
$\left \langle \tilde \lambda^{u}, k \right \rangle \left[ \dot{p}_{k} \right]_{\tilde \jmath} 
+ \left \langle \dot{\lambda}^{u}, k \right \rangle \left[ \tilde p_{k} \right]_{\tilde \jmath}$
from the tail of $ \Pi^{ \tilde \jmath}_{P}\left( y \left( \tau \right) \right)$. 
Hence $y^{\tilde \jmath}_{P} \left( \tau \right) \in W^{1}_{\nu_{u}}$. 
An analogous decomposition $\Pi^{ \tilde \jmath}_{Q}\left( y \left( \tau \right) \right)  = y^{\tilde \jmath}_{Q} \left( \tau \right) + y^{\tilde \jmath}_{Q,\lambda} \left( \tau \right)$
is used to analyze the terms associated to the stable manifold. 
\begin{lemma}
	\label{lemma:D2C_bound}
	Let $1 \leq \tilde \jmath \leq n$ and $0 < r \leq r^{\ast}$, then 
	\begin{align*}
		\sup_{\tau \in [0,1]} \left \Vert y^{\tilde \jmath}_{P} \left( \tau \right) \right \Vert_{\nu_{u}} \leq \sup \boldsymbol{\sigma}^{\tilde \jmath}_{u}, \quad 
		\sup_{\tau \in [0,1]} \left \Vert y^{\tilde \jmath}_{Q} \left( \tau \right) \right \Vert_{\nu_{s}} \leq \sup \boldsymbol{\sigma}^{\tilde \jmath}_{s},
	\end{align*}
	where
	\begin{align*}
		\boldsymbol{\sigma}^{\tilde \jmath}_{u} :&= 
			\sum_{i_{3}=1}^{n} \left[ \sum_{i=1}^{n_{u}} \left \vert \left[ \hat p_{e_{i}} \right]_{i_{3}} \right \vert \nu_{u} + r^{\ast} \right] 
			\sum_{i_{1},i_{2}=1}^{n} \left \vert \dfrac{ \partial^{3} g_{\tilde \jmath} }{ \partial x_{i_{1}} \partial x_{i_{2}} \partial x_{i_{3}} } 
			\left( \boldsymbol{ \hat p_{0}} \right) \right \vert 
			\\[2ex] & \quad + 		
			3 \sum_{i_{1}, i_{2} = 1}^{n} 
			\left \vert \frac{ \partial^{2} g_{ \tilde \jmath } }{ \partial x_{i_{1}} \partial x_{i_{2}} } \left( \boldsymbol{ \hat p_{0}} \right) \right \vert + 
			D^{2} \tilde g_{\tilde \jmath} \left( \left \Vert \left[ \hat p \right]_{1} \right \Vert_{\nu_{u}} + r^{\ast}, \ldots, 
			\left \Vert \left[ \hat p \right]_{n} \right \Vert_{\nu_{u}} + r^{\ast} \right) \left[ \bold{1}_{n}, \bold{1}_{n} \right] \\[2ex] & \quad + 
			2 \left( \left \vert K^{u} \right \vert + 1 \right).
	\end{align*}
	 The bounds $\left \{ \boldsymbol{\sigma}^{\tilde \jmath}_{s} : 1 \leq \tilde \jmath \leq n \right \}$ associated to the
	 stable manifold are defined analogously. 
	\begin{proof}
		Let $\tau \in [0,1]$ and $1 \leq \tilde \jmath \leq n$ be arbitrary. It follows directly from \eqref{eq:Z2_Pi_Taylor}
		that
		\begin{align}
			\label{eq:Pi_k0}
			\left \vert \left( y^{\tilde \jmath}_{P} \left( \tau \right) \right)_{0} \right \vert \leq 
			\sum_{i_{1}, i_{2} = 1}^{n} 
			\left \vert \frac{ \partial^{2} g_{ \tilde \jmath } }{ \partial x_{i_{1}} \partial x_{i_{2}} } \left( \hat p_{0} + \tau r \dot{p}_{0} \right) \right \vert,
		\end{align}
		since $\left \Vert \dot{p} \right \Vert_{\nu_{u},n},  \left \Vert \tilde p \right \Vert_{\nu_{u},n} \leq 1$.	
		Next, we consider the first order components of $y^{\tilde \jmath}_{P} \left( \tau \right)$ by separately analyzing
		the terms in \eqref{eq:Z2_Pi_Taylor} for $\left \vert k \right \vert =1$. To this end, observe that 
		\begin{align*}
			& \sum_{i=1}^{n_{u}} \left \vert D^{3} g_{\tilde \jmath} \left( \hat p_{0} + \tau r \dot{p}_{0} \right) 
			\left[ \dot{p}_{0}, \tilde {p}_{0}, \hat p_{e_{i}} + \tau r \dot{p}_{e_{i}} \right]  \right \vert \nu_{u} \\[2ex] & \quad \leq 
			\sum_{i_{3}=1}^{n} \left[ \sum_{i=1}^{n_{u}} \left \vert \left[ \hat p_{e_{i}} \right]_{i_{3}} \right \vert \nu_{u} + r^{\ast} \right] 
			\sum_{i_{1},i_{2}=1}^{n} \left \vert \dfrac{ \partial^{3} g_{\tilde \jmath} }{ \partial x_{i_{1}} \partial x_{i_{2}} \partial x_{i_{3}} } 
			\left( \hat p_{0} + \tau r \dot{p}_{0} \right) \right \vert,
		\end{align*}
		and
		\begin{align*}
			\sum_{i=1}^{n_{u}} \left \vert D^{2} g_{\tilde \jmath} \left( \hat p_{0} + \tau r \dot{p}_{0} \right) \left[ \tilde {p}_{0}, \dot{p}_{e_{i}} \right]  			\right \vert \nu_{u} &\leq \sum_{i_{1}, i_{2} = 1}^{n} 
			\left \vert \frac{ \partial^{2} g_{ \tilde \jmath } }{ \partial x_{i_{1}} \partial x_{i_{2}} } \left( \hat p_{0} + \tau r \dot{p}_{0} \right) \right 					\vert, 
			\\[2ex]
			\sum_{i=1}^{n_{u}} \left \vert \tilde \lambda^{u}_{i} \left[ \dot{p}_{e_{i}} \right]_{\tilde \jmath} + 
			\dot{\lambda}^{u}_{i} \left[ \tilde p_{e_{i}} \right]_{\tilde \jmath} \right \vert \nu_{u} &\leq 2,
		\end{align*}
		where we used that $r \leq r^{\ast}$, $\left \Vert \tilde \lambda^{u} \right \Vert_{\infty}, 
		\left \Vert \dot{\lambda}^{u} \right \Vert_{\infty}, \left \Vert \dot{p} \right \Vert_{\nu_{u}.n},  \left \Vert \tilde p \right \Vert_{\nu_{u},n} \leq 1$. 
		It now follows from the expression in \eqref{eq:Z2_Pi_Taylor} that 
		\begin{align}
			\sum_{\substack{ k \in \NN_{0}^{n_{u}} \\ \left \vert k \right \vert = 1} } 
			\left \vert \left( y^{\tilde \jmath}_{P} \left( \tau \right) \right)_{k} \right \vert \nu_{u} &\leq 
			\sum_{i_{3}=1}^{n} \left[ \sum_{i=1}^{n_{u}} \left \vert \left[ \hat p_{e_{i}} \right]_{i_{3}} \right \vert \nu_{u} + r^{\ast} \right] 
			\sum_{i_{1},i_{2}=1}^{n} \left \vert \dfrac{ \partial^{3} g_{\tilde \jmath} }{ \partial x_{i_{1}} \partial x_{i_{2}} \partial x_{i_{3}} } 
			\left( \hat p_{0} + \tau r \dot{p}_{0} \right) \right \vert \nonumber \\[2ex] & \quad + 2 \sum_{i_{1}, i_{2} = 1}^{n} 
			\left \vert \frac{ \partial^{2} g_{ \tilde \jmath } }{ \partial x_{i_{1}} \partial x_{i_{2}} } \left( \hat p_{0} + \tau r \dot{p}_{0} \right) \right \vert 			+ 2.
			\label{eq:Pi_k1}
		\end{align}

		Finally, we consider the remainder of $y^{\tilde \jmath}_{P}$, i.e., the components associated to $\left \vert k \right \vert \geq 2$. 
		First, note that 
		\begin{align*}
			\sum_{ \substack{ k \in \mathcal{K}^{u} \\ \left \vert k \right \vert \geq 2 } }
			\left \vert \left \langle \tilde \lambda^{u}, k \right \rangle \left[ \dot{p}_{k} \right]_{\tilde \jmath} + 
			\left \langle \dot{\lambda}^{u}, k \right \rangle \left[ \tilde p_{k} \right]_{\tilde \jmath} \right \vert \nu_{u}^{k}
			\leq 2 \left \vert K^{u} \right \vert,
		\end{align*}
		since $\left \Vert \tilde \lambda^{u} \right \Vert_{\infty}, \left \Vert \dot{\lambda}^{u} \right \Vert_{\infty}, \left \Vert \dot{p} \right 
				\Vert_{\nu_{u}.n},  \left \Vert \tilde p \right \Vert_{\nu_{u},n} \leq 1$ and $\left \vert k \right \vert \leq \left \vert K^{u} \right \vert$
		for $k  \in \mathcal{K}^{u}$. 
		Furthermore, 
		\begin{align*}
			\sum_{ \substack{ k \in \NN_{0}^{n_{u}} \\ \left \vert k \right \vert \geq 2 } } \left \vert D^{2}C_{k} \left( \hat p + \tau r \dot{p} \right) \left[ \dot{p}, \tilde p \right] \right \vert \nu_{u}^{k}
			\leq D^{2} \tilde g_{\tilde \jmath} \left( \left \Vert \left[ \hat p \right]_{1} \right \Vert_{\nu_{u}} + r^{\ast}, \ldots, 
			\left \Vert \left[ \hat p \right]_{n} \right \Vert_{\nu_{u}} + r^{\ast} \right) \left[ \bold{1}_{n}, \bold{1}_{n} \right]
		\end{align*}
		by Lemma \ref{lemma:D2g}. Therefore, 
		\begin{align}
			\sum_{ \substack{ k \in \NN_{0}^{n_{u}} \\ \left \vert k \right \vert \geq 2 } } 
			\left \vert  \left( y^{\tilde \jmath}_{P} \left( \tau \right) \right)_{k} \right \vert \nu_{u}^{k} \leq
			2 \left \vert K^{u} \right \vert + 
			D^{2} \tilde g_{\tilde \jmath} \left( \left \Vert \left[ \hat p \right]_{1} \right \Vert_{\nu_{u}} + r^{\ast}, \ldots, 
			\left \Vert \left[ \hat p \right]_{n} \right \Vert_{\nu_{u}} + r^{\ast} \right) \left[ \bold{1}_{n}, \bold{1}_{n} \right].
			\label{eq:Pi_k2}
		\end{align} 
		Altogether, the sum of \eqref{eq:Pi_k0}, \eqref{eq:Pi_k1} and \eqref{eq:Pi_k2} yields an upper bound for 
		$\left \Vert y^{\tilde \jmath}_{P} \left( \tau \right) \right \Vert_{\nu_{u}}$, for any $\tau \in [0,1]$ and $r \leq r^{\ast}$, and is
		contained in $\boldsymbol{\sigma}^{\tilde \jmath}_{u}$. This proves the result. 
	\end{proof}
\end{lemma}

It is now a straightforward task to compute bounds for \eqref{eq:Z2_Taylor}. 
\begin{lemma}
	\label{lemma:Z2_Taylor_scal_cheb}
	Let $1 \leq \tilde \jmath \leq n$ and $0 < r \leq r^{\ast}$, then  
	\begin{align}
		\label{eq:Z2_P_scal_cheb}
		\sup_{\tau \in [0,1] } \left \Vert \Pi A \Pi^{ \tilde \jmath }_{P} \left( y \left( \tau \right) \right) \right \Vert_{ \Pi \left( \mathcal{X}_{\nu} \right) } 
		&\leq 
		\left \Vert \Pi A_{NK} \Pi^{ \tilde \jmath }_{P} \right \Vert_{ \mathcal{B} \left( W^{1}_{\nu_{u}}, \Pi \left( \mathcal{X}_{\nu} \right) \right)}	
		\sup \boldsymbol{\sigma}^{\tilde \jmath}_{u}, 
	\end{align}
	for $\Pi \in \left \{ \Pi_{\theta}^{j}, \Pi_{\phi}^{j}, \Pi_{\lambda^{u}}^{j}, \Pi_{\lambda^{s}}^{j}, \Pi^{ij}_{a}, \Pi^{j}_{Q} \right \}$ and
	\begin{align}
		\label{eq:Z2_Q_scal_cheb}
		\sup_{\tau \in [0,1] } \left \Vert \Pi A \Pi^{ \tilde \jmath }_{Q} \left( y \left( \tau \right) \right) \right \Vert_{ \Pi \left( \mathcal{X}_{\nu} \right) } 
		&\leq 
		\left \Vert \Pi A_{NK} \Pi^{ \tilde \jmath }_{Q} \right \Vert_{ \mathcal{B} \left( W^{1}_{\nu_{s}}, \Pi \left( \mathcal{X}_{\nu} \right) \right)}	
		\sup \boldsymbol{\sigma}^{\tilde \jmath}_{s}	
	\end{align}
	for $\Pi \in \left \{ \Pi_{\theta}^{j}, \Pi_{\phi}^{j}, \Pi_{\lambda^{u}}^{j}, \Pi_{\lambda^{s}}^{j}, \Pi^{ij}_{a}, \Pi^{j}_{P} \right \}$.
	\begin{proof}
		It suffices to observe that 
		\begin{align*}
			\Pi A \Pi^{ \tilde \jmath }_{P} \left( y \left( \tau \right) \right)= 
			\Pi A_{NK} \Pi^{ \tilde \jmath \mathcal{K}^{u} }_{P} \left( y \left( \tau \right) \right), 
			\quad \Pi \in \left \{ \Pi_{\theta}^{j}, \Pi_{\phi}^{j}, \Pi_{\lambda^{u}}^{j}, \Pi_{\lambda^{s}}^{j}, 
					\Pi^{ij}_{a}, \Pi^{j}_{Q} \right \},
		\end{align*}
		by construction of the approximate inverse $A$. Therefore, since
		\begin{align*}
			\Pi^{ \tilde \jmath \mathcal{K}^{u} }_{P} \left( y \left( \tau \right) \right) = \Pi^{ \tilde \jmath \mathcal{K}^{u} }_{P} \left( y^{\tilde \jmath}
			_{P} \left( \tau \right) \right), 
		\end{align*}
		the result follows from Lemma \ref{lemma:D2C_bound}. 
	\end{proof}
\end{lemma}
\begin{lemma}
	\label{lemma:Z2_Taylor_Taylor}
	 Let $1 \leq \tilde \jmath, j \leq n$ and $0<r \leq r^{\ast}$, then 
	 \begin{align*}
	 	& \sup_{\tau \in [0,1] } \left \Vert \Pi^{j}_{P} A \Pi^{ \tilde \jmath }_{P} \left( y \left( \tau \right) \right) \right \Vert_{\nu_{u}} \\[2ex] & \quad \leq 
		\begin{cases}
			\left \Vert \Pi^{j}_{P} A_{NK} \Pi^{ \tilde \jmath }_{P} \right \Vert_{ \mathcal{B} \left( W^{1}_{\nu_{u}},  W^{1}_{\nu_{u}} \right)}
			\sup \ \boldsymbol{\sigma}^{\tilde \jmath}_{u}, &
			\tilde \jmath \not = j, \\[3ex]	
			\normalfont
			\max \left \{ \left \Vert \Pi^{j}_{P} A_{NK} \Pi^{j}_{P} \right \Vert_{\mathcal{B} \left( W^{1}_{\nu_{u}}, W^{1}_{\nu_{u}} \right) },  
			\left( \displaystyle\min_{1 \leq i \leq n_{u}} \left \vert \text{Re} \left( \hat \lambda^{u}_{i} \right) \right \vert  
			\left( K^{u}_{i} + 1 \right) \right)^{-1}\right \}
			\cdot \normalfont \sup \boldsymbol{\sigma}^{\tilde \jmath}_{u}
			\\[4ex] \quad  + 
			2  \left( \displaystyle\min_{1 \leq i \leq n_{u}} \left \vert \text{Re} \left( \hat \lambda^{u}_{i} \right) \right \vert \right)^{-1},
			& \tilde \jmath = j.
		\end{cases}
	 \end{align*}
	 The statement and corresponding bound for the stable manifold is analogous. 
	 \begin{proof}
		If $\tilde \jmath \not = j$, then 
		\begin{align*}
			\Pi^{j}_{P} A \Pi^{ \tilde \jmath }_{P} \left( y \left( \tau \right) \right) = 
			\Pi^{j}_{P} A_{NK} \Pi^{ \tilde \jmath \mathcal{K}^{u} }_{P}\left( y^{\tilde \jmath}_{P} \left( \tau \right) \right)
		\end{align*}
		by definition of the approximate inverse and $y^{\tilde \jmath}_{P} \left( \tau \right)$. 
		Therefore, in this case, the result follows immediately from Lemma \ref{lemma:D2C_bound}. 
		Now, assume that $\tilde \jmath = j$ and observe that  
		\begin{align*}
			\Pi^{j}_{P} A \Pi^{ j }_{P} \left( y \left( \tau \right) \right) =
			\Pi^{j}_{P} A \Pi^{j}_{P} \left( y^{j}_{P} \left( \tau \right) \right) + 
			\Pi^{j}_{P} A \Pi^{j}_{P} \left( y^{j}_{P,\lambda} \left( \tau \right) \right). 
		\end{align*}
		In particular, since $y^{j}_{P} \left( \tau \right) \in W^{1}_{\nu_{u}}$, we may bound the first term in the
		above expression by  
		\begin{align*}
			\left \Vert \Pi^{j}_{P} A \Pi^{j}_{P} \left( y^{j}_{P} \left( \tau \right) \right) \right \Vert_{\nu_{u}} &\leq 
			\left \Vert \Pi^{j}_{P} A \Pi^{j}_{P} \right \Vert_{\mathcal{B} \left( W^{1}_{\nu_{u}}, W^{1}_{\nu_{u}} \right) } 
			\left \Vert y^{j}_{P} \left( \tau \right) \right \Vert_{\nu_{u}} \\[2ex] &\leq 
			\left \Vert \Pi^{j}_{P} A \Pi^{j}_{P} \right \Vert_{\mathcal{B} \left( W^{1}_{\nu_{u}}, W^{1}_{\nu_{u}} \right) } 
			\sup \ \boldsymbol{ \sigma}^{j}_{u},
		\end{align*}
		where in the last line we used Lemma \ref{lemma:D2C_bound} again. 
		
		Next, we derive a more explicit expression for 
		$\left \Vert \Pi^{j}_{P} A \Pi^{j}_{P} \right \Vert_{\mathcal{B} \left( W^{1}_{\nu_{u}}, W^{1}_{\nu_{u}} \right) }$. To this end, 
		recall the definition of $\Pi^{j}_{P} A \Pi^{j}_{P}$ (see Definition \ref{def:approximate_inverse}) and observe that 
		\begin{align*}
			\left \vert \left \langle \hat \lambda^{u}, k \right \rangle \right \vert \geq 
			\min_{1 \leq i \leq n_{u}} \left \vert \text{Re} \left( \hat \lambda^{u}_{i} \right) \right \vert \left( K^{u}_{i} + 1 \right)
		\end{align*} 
		for any $k \in \mathcal{K}^{u}_{c}$. Therefore, it follows from Proposition \ref{prop:operator_norm_mult} that 
		\begin{align*}
			& \left \Vert \Pi^{j}_{P} A \Pi^{j}_{P} \right \Vert_{\mathcal{B} \left( W^{1}_{\nu_{u}}, W^{1}_{\nu_{u}} \right) } \\[2ex] & \quad \leq
			\max \left \{ \left \Vert \Pi^{j}_{P} A_{NK} \Pi^{j}_{P} \right \Vert_{\mathcal{B} \left( W^{1}_{\nu_{u}}, W^{1}_{\nu_{u}} \right) },  
			\left( \min_{1 \leq i \leq n_{u}} \left \vert \text{Re} \left( \hat \lambda^{u}_{i} \right) \right \vert \left( K^{u}_{i} + 1 \right) \right)^{-1}\right 			\}.
		\end{align*}		
		Finally, a straightforward computation shows that 
		\begin{align*}
			\left \Vert \Pi^{j}_{P} A \Pi^{j}_{P} \left( y^{j}_{P,\lambda} \left( \tau \right) \right) \right \Vert_{\nu_{u}} &=
			\sum_{k \in \mathcal{K}^{u}_{c}} 
			\left \vert \left \langle \tilde \lambda^{u}, k \right \rangle \left[ \dot{p}_{k} \right]_{j} + 
			\left \langle \dot{\lambda}^{u}, k \right \rangle \left[ \tilde p_{k} \right]_{j} \right \vert 
			\left \vert \left \langle \hat \lambda^{u}, k \right \rangle \right \vert^{-1}
			\nu_{u}^{k} \\[2ex] &\leq 2  
			\left( \min_{1 \leq i \leq n_{u}} \left \vert \text{Re} \left( \hat \lambda^{u}_{i} \right) \right \vert \right)^{-1},
		\end{align*}
		since 		
		$\left \vert \left \langle \lambda^{u}, k \right \rangle \right \vert \leq \left \Vert \lambda^{u} \right \Vert_{\infty} \left \vert k \right \vert$
		and $\left \vert \left \langle \hat \lambda^{u}, k \right \rangle \right \vert \geq \min_{1 \leq i \leq n_{u}}
			\left \vert \text{Re} \left( \hat \lambda^{u}_{i} \right) \right \vert \left \vert k \right \vert$
		for any $k \in \NN_{0}^{n_{u}}$ and $\lambda^{u} \in \CC^{n_{u}}$.
		Altogether, this proves the result. 
	 \end{proof}
\end{lemma}

\paragraph{Second order coefficients of $Z_{\Pi}(r)$}
We are now ready to finish the construction of the quadratic polynomials
$Z_{\Pi}(r)$ for $\Pi \in \mathbb{P}$. As before, we first introduce some additional
notation. We will denote the bounds in \eqref{eq:Z2_t0}, \eqref{eq:Z2_tm},  
\eqref{eq:Z2_P_scal_cheb}, \eqref{eq:Z2_Q_scal_cheb} and 
Lemmas \ref{lemma:Z2_cheb_scal_Taylor}, \ref{lemma:Z2_cheb_cheb}, \ref{lemma:Z2_Taylor_Taylor} by 
$Z^{2, \tilde \jmath t_{0} }_{\Pi}$, $Z^{2, \tilde \jmath t_{m}}_{\Pi}$, 
$Z^{2, \tilde \jmath P}_{\Pi}$, $Z^{2, \tilde \jmath Q}_{\Pi}$,
$Z^{2, \tilde \imath \tilde \jmath \NN}_{\Pi}$, 
$Z^{2, \tilde \imath \tilde \jmath \NN}_{\Pi^{ij}_{a}}$
and $Z^{2, \tilde \jmath P}_{\Pi^{j}_{P}}$, $Z^{2, \tilde \jmath Q}_{\Pi^{j}_{Q}}$,
respectively. 
Finally, we set
\begin{align*}
		Z^{2}_{\Pi} := \sum_{ \tilde \jmath =1}^{n} \left( 
		Z^{2, \tilde j t_{0} }_{\Pi} + Z^{2, \tilde j t_{m} }_{\Pi} + 
		\sum_{\tilde \imath =1}^{m} Z^{2, \tilde \imath \tilde \jmath \NN}_{\Pi}
		+ Z^{2, \tilde \jmath P}_{\Pi} + Z^{2, \tilde \jmath Q}_{\Pi}
		\right), \quad \Pi \in \mathbb{P},
\end{align*}
and define
\begin{align*}
		Z_{\Pi}(r) := Z^{1}_{\Pi} + Z^{2}_{\Pi}r, \quad \Pi \in \mathbb{P}. 
\end{align*}
Then $Z_{\Pi}(r)$ satisfies \eqref{eq:Z} by construction.

\section{Applications: traveling fronts in parabolic PDEs}
\label{sec:applications}
In this section we use our method to prove the 
existence of connecting orbits in systems of ODEs which arise from the study of 
traveling fronts in scalar parabolic PDEs. In addition, we perform discrete
continuation (discrete in the sense that we rigorously validate the solution for many parameter values, but we do not attempt to obtain a continuous parametrized branch of solutions).
This also demonstrates the effectiveness of the phase condition introduced
in Definition \ref{def:time_parameterization}. 
Before we proceed to the applications, we first give a rough outline of our main procedure for validating 
connecting orbits. We have tried to automate as many steps as possible, but there 
are still certain steps which are based on experimentation.

\paragraph{Step 1: Compute parameterizations of the local (un)stable manifolds}
\begin{itemize}
	\item[$1.1$] 	Compute numerical approximations $\tilde p_{0}$ and $\tilde q_{0}$ of the equilibria 
				of interest. 
	\item[$1.2$] 	Compute numerical approximations $\left \{ \left( \tilde \lambda_{k}^{u}, \tilde p_{k} \right) : \left \vert k \right \vert =1 \right \}$
				and $\left \{ \left( \tilde \lambda_{k}^{s}, \tilde q_{k} \right) : \left \vert k \right \vert =1 \right \}$ of the eigendata associated to
				$Dg \left( \tilde p_{0} \right)$ and $Dg \left( \tilde q_{0} \right)$, respectively. In this step we set the length of the
				approximate eigenvectors to one. 
	\item[$1.3$]	Choose the number of Taylor coefficients $K^{u} \in \NN_{0}^{n_{u}}$ and $K^{s} \in \NN_{0}^{n_{s}}$ 
	and compute approximate zeros  $\left( \hat \lambda^{u}, \hat p \right)$,
				$\left( \hat \lambda^{s}, \hat q \right)$
				of the mappings 
				\begin{align*}
					\left( \lambda^{u}, p \right) \mapsto \begin{bmatrix}
					\left[ \Pi_{P}^{j \mathcal{K}^{u} }F_{P} \left( \lambda^{u}, p \right) \right]_{j=1}^{n} \\
					\left[ \left \langle p_{e_{k}}, \tilde p_{e_{k}} \right \rangle \right]_{k=1}^{n_{u}}
					\end{bmatrix}, \quad 
					\left( \lambda^{s}, q \right) \mapsto \begin{bmatrix}
					\left[ \Pi_{Q}^{j \mathcal{K}^{s} }F_{Q} \left( \lambda^{s}, q \right) \right]_{j=1}^{n} \\
					\left[ \left \langle q_{e_{k}}, \tilde q_{e_{k}} \right \rangle \right]_{k=1}^{n_{s}}
					\end{bmatrix}
				\end{align*}
				by Newton's method. 
	\item[$1.4$]	If necessary, increase the truncation parameters and rescale 
				the eigenvectors so that validation is feasible, see Remark \ref{remark:manifold_parameters} below.
\end{itemize}

\paragraph{Step 2: Compute an accurate approximation of a connecting orbit}
\begin{itemize}
	\item[$2.1$] 	Compute a numerical approximation of a connecting orbit. This step is based on solving
				the truly nonlinear part of the problem and involves experimentation. It is obviously 
				problem dependent. 
	\item[$2.2$] 	Use the domain decomposition algorithm developed in \cite{domaindecomposition} to 
				compute a grid $\left(t_{i} \right)_{i=0}^{m}$ and an accurate approximate connecting orbit
				\begin{align*}
					\hat u = \sum_{i=1}^{m} \bold{1}_{ \left[ t_{i-1}, t_{i} \right]} \left( \hat a_{0} + 2 \sum_{k=1}^{N_{i}-1} \hat a^{i}_{k} T^{i}_{k} \right),
				\end{align*}
				so that the decay rates of the Chebyshev coefficients $\hat a^{i}$ are equidistributed over the subdomains $\left[t_{i-1}, t_{i} \right]$. 
				The number of modes $N_{i}$ is chosen in such a way that $\left \vert \hat a^{i}_{N_{i}-1} \right \vert \approx 10^{-16}$. 
				The number of subdomains $m$ is determined by experimentation. In general, we use as many subdomains as necessary in order to 
				ensure high decay rates of the Chebyshev coefficients. 
	\item[$2.3$] 	Use $\hat u$ as a reference orbit to fix the time parameterization of the connecting orbit 
				(see Definition \ref{def:time_parameterization}). 
\end{itemize}

\paragraph{Step 3: Validate the connecting orbit and (un)stable manifolds}
\begin{itemize}
	\item [$3.1$]	Combine the results from the previous two steps to construct a \emph{symmetric} approximate zero 
				$\hat x = \left( \hat \theta, \hat \phi, \hat \lambda^{u}, \hat \lambda^{s}, \hat a, \hat p, \hat q \right)$ of $F_{NK}$
				(see Remark \ref{remark:symmetry}). 
	\item[$3.2$] 	Set $r^{\ast} = 10^{-5}$, $\nu_{u} = \nu_{s} = 1$ (see Remark \ref{remark:manifold_parameters}) and 
				compute the weights $\left( \nu_{i} \right)_{i=1}^{m}$ as explained in Remark \ref{remark:nu_i} below. 	
	\item[$3.3$]	Initialize the numerical data with interval arithmetic and construct the $n(m+4)+2$ radii-polynomials
				\begin{align*}
					p_{\Pi}(r) := Z^{2}_{\Pi} r^{2} + \left( Z^{1}_{\Pi}-1 \right) r + Y_{\Pi}, \quad 
					\Pi \in \mathbb{P}. 
				\end{align*}
	\item[$3.4$]	Determine an interval $I$ on which all the radii polynomials are negative. 
\end{itemize}
If we fail to find an interval $I$ on which all the radii polynomials are negative,
we try to determine which parameters (the truncation parameters, the weights $\nu_{i}$ or the
``scalings'' of the coefficients $\hat p$ and $\hat q$) need to be modified by ``visual'' inspection and try again.

\begin{remark}[Scaling and the number of Taylor coefficients]
	\label{remark:manifold_parameters}
	Observe that for any $\nu_{u}, \tilde \nu_{u} >0$ it holds that 
	$p \in W^{1}_{\nu_{u}}$ if and only if $\frac{ \nu_{u} }{ \tilde \nu_{u}}p \in W^{1}_{\tilde \nu_{u}}$ with 
	the scaling notation introduced in Remark \ref{remark:scaling}. 
	Therefore, since the parameterization mappings $F_{P}$ and $F_{Q}$ are invariant under the rescaling $p \mapsto \mu p$
	(see Remark \ref{remark:scaling}), we have chosen to set $\nu_{u} = \nu_{s} = 1$ and search for appropriate scalings which
	ensure that $\hat p$ and $\hat q$ decay sufficiently fast to zero. To be more precise, we explain in detail how we choose 
	the scalings of the eigenvectors and the number of Taylor coefficients for the unstable manifold 
	(the procedure for the stable manifold is analogous). 
			
	The main idea is to choose the scalings and number of Taylor coefficients in such a way that the bound 
	for $\Pi_{P}^{j} \left( DF \left( \hat x \right) - \widehat{DF} \right)$ is below
	some prescribed tolerance. More precisely, in light of \eqref{eq:Pi_parm}, \eqref{eq:Dc_taylor} and Lemma \ref{lemma:Lambda}, we aim
	to find a truncation parameter $K^{u} \in \NN_{0}^{n_{u}}$ and a scaling factor $\mu \in\left(0, \infty \right)^{n_{u}}$ such that  
	\begin{align}
		\label{eq:mu_K}
		\left[ \min_{1 \leq i \leq n_{u}} \left( K^{u}_{i} + 1 \right)
		\min_{1 \leq i \leq n_{u}} \left \vert \normalfont \text{Re} \left( \hat \lambda_{i} \right) \right \vert \right]^{-1} 
		\left \Vert  \mu \hat G^{jl} \right \Vert_{1} \leq \varepsilon_{u}, 
		\quad 1 \leq j, l \leq n,
	\end{align}
	where $\varepsilon_{u} >0$ (in practice we set $\varepsilon_{u} = \frac{1}{2}$). 	
	We start by determining $K^{u} \in \NN_{0}^{n_{u}}$. To this end, observe that the scaling factor $\mu$ has no effect
	on $\left \vert \hat G^{jl}_{0} \right \vert$. For this reason, we set $K^{u} := \left( \max_{1 \leq j,l \leq n_{u} } K^{u}_{jl} \right) \bold{1}_{n_{u}}$, where
	$K^{u}_{jl} \in \NN$ is the smallest integer such that 
	\begin{align*}
		K^{u}_{jl} > \dfrac{\left \vert \hat G^{jl}_{0} \right \vert}{ \zeta \varepsilon_{u} \min_{1 \leq i \leq n_{u}} \normalfont \left \vert \text{Re} \left( \hat \lambda^{u}_{i} \right) \right \vert}-1,
		\quad 1 \leq j,l \leq n, \quad \zeta \in (0,1]. 
	\end{align*}
	Here $\zeta \in (0,1]$ is an additional parameter chosen through experimentation (in practice we use $\zeta = \frac{3}{4}$). 
	
	Next, we determine an appropriate scaling factor $\mu$. Let $1 \leq j,l \leq n$ and approximate 
	$\left \vert \hat G^{jl}_{k} \right \vert \approx \left \vert \hat G^{jl}_{0} \right \vert \rho_{jl}^{-\left \vert k \right \vert}$, where $\rho_{jl} =  e^{- s_{jl} }$ and
	$s_{jl}$ is the slope of the best line through the points
	\begin{align}
		\label{eq:decay_mult_array}
		\left \{ \left(d, \log \left( \sum_{\left \vert k \right \vert = d} \left \vert \hat G^{jl}_{k} \right \vert \right) \right)
		: \   \sum_{\left \vert k \right \vert = d} \left \vert \hat G^{jl}_{k} \right \vert  > 10^{-16}, \ 0 \leq d \leq \left \vert M_{jl}K^{u} \right \vert \right \}.
	\end{align} 
	Recall that $M_{jl} = \normalfont \text{order} \left( \dfrac{ \partial g_{j} }{ \partial x_{l} } \right)$. 
	Now, if $\mu_{\imath} < \rho_{jl}$ for all $1 \leq \imath \leq n_{u}$ and $\left \vert K^{u} \right \vert$ is sufficiently large, then 
	\begin{align*}
		\left \Vert \mu \hat G^{jl} \right \Vert_{1} \approx \left \vert \hat G^{jl}_{0} \right \vert \prod_{\imath=1}^{n_{u}} 
		\dfrac{ \rho_{jl} }{ \rho_{jl} - \mu_{\imath}}. 
	\end{align*}
	Motivated by this observation and the inequality in \eqref{eq:mu_K}, we set $\mu_{\imath} = \mu$ for $1 \leq \imath \leq n_{u}$ and require that 
	\begin{align*}
		\dfrac{ \rho_{jl} }{\rho_{jl} - \mu} \leq 
		\left[ \dfrac{ \varepsilon_{u} \min_{1 \leq i \leq n_{u}} \left( K^{u}_{i} + 1 \right)
		\min_{1 \leq i \leq n_{u}} \left \vert \normalfont \text{Re} \left( \hat \lambda^{u}_{i} \right) \right \vert }{ \left \vert \hat G^{jl}_{0} \right \vert } \right]^{\frac{1}{d}} =: \xi_{jl}
	\end{align*}
	for all $1 \leq j, l \leq n$. Therefore, we set
	\begin{align*}
		\mu := \min_{1 \leq j,l \leq n_{u} } \rho_{jl} \dfrac{ \xi_{jl} -1}{ \xi_{jl} }.
	\end{align*}
	
	We remark that one could determine a more ``refined'' scaling factor $\mu$, which need not
	be the same in each direction, by taking the decay rates of $\hat G^{jl}$ in each separate direction 
	into account (as opposed to using the ``uniform'' rate in \eqref{eq:decay_mult_array} which ignores
	the different directions of the array). In addition, one could take the different sizes of the eigenvalues into account in the definition of 
	$\xi_{jl}$. 
\end{remark}

\begin{remark}
	\label{remark:nu_i}
	To determine the weights $\left( \nu_{i} \right)_{i=1}^{m}$, we use a heuristic procedure slightly more refined
	than the one used in \cite{domaindecomposition}. Namely, we try to ensure that the bound for the tail 
	of $\Pi_{a}^{ij} \left( DF \left( \hat x \right) - \widehat{DF} \right)$ is below some prescribed
	tolerance  (rather than requiring the residual to be below some tolerance as in \cite{domaindecomposition}). 
	More precisely, in light of \eqref{eq:Pi_cheb}, we require that 
	\begin{align}
		\label{eq:nu_i1}
		\frac{L \left( t_{i} - t_{i-1} \right)}{2N_{i}} \left( \nu_{i} + \nu^{-1}_{i} \right) \left \Vert \hat g^{ijl} \right \Vert_{\nu_{i}} \leq \varepsilon,
		\quad 1 \leq i \leq m, \ 1 \leq j, l \leq n. 
	\end{align}
	where $\varepsilon >0$ is some prescribed tolerance (in practice we set $\varepsilon = \frac{1}{2}$). 
	
	We use the rough approximation  
	$\left \vert \hat g_{k}^{ijl} \right \vert \approx \left \vert \hat g^{ijl}_{0} \right \vert \rho^{-k}_{ijl}$, where
	$\rho_{ijl} = e^{-s_{ijl}}$ and $s_{ijl}$ is the slope of the
	best line through the points
	\begin{align*}
		\left \{ \left(k,\log \left( \left \vert \hat g^{ijl}_{k} \right \vert \right) \right): \ 0 \leq k \leq M_{jl} \left( N_{i} - 1 \right), \
	 	\left \vert \hat g^{ijl}_{k} \right \vert > 10^{-16} \right \}.
	\end{align*} 
	In practice, $\rho_{ijl}$ is roughly the same for all $1 \leq i \leq m$ and $1 \leq j,l \leq n$ due to the choice of the grid,
	and we therefore write
	$\rho = \rho_{ijl}$. If $\nu_{i} < \rho$ and $N_{i}$ is sufficiently large, then 
	\begin{align}
		\label{eq:nu_i2}
		\left \Vert \hat g^{ijl} \right \Vert_{\nu_{i}} \approx 
		\left \vert \hat g^{ijl}_{0} \right \vert \left( 1 + 2 \sum_{k=1}^{M_{jl} \left( N_{i} - 1 \right)} \left( \dfrac{ \nu_{i} }{ \rho } \right)^{k} \right) \approx
		\left \vert \hat g^{ijl}_{0} \right \vert \left( \dfrac{ 2 \rho  }{ \rho - \nu_{i}} -1 \right).
	\end{align}
	Altogether, \eqref{eq:nu_i1} and \eqref{eq:nu_i2} yield the constraint 
	\begin{align}
		\label{eq:nu_i3}
		\nu_{i}^{3} + \left( \alpha_{ijl} + \rho \right) \nu_{i}^{2} + \left( 1 - \alpha_{ijl} \rho  \right) \nu_{i} + \rho \leq 0, \quad
		\alpha_{ijl} := \dfrac{ 2N_{i} \varepsilon}{ L \left( t_{i} - t_{i-1} \right) \left \vert \hat g^{ijl}_{0} \right \vert }.
	\end{align}
	Finally, we determine an interval $\left[ \nu_{\min}, \nu_{\max} \right] \subset \RR_{>0}$ on which \eqref{eq:nu_i3} is satisfied
	for all $1 \leq i \leq m$ and $1 \leq j,l\leq n$. Then, if $\nu_{\max} >1$, we choose a weight $\hat \nu \in \left[ \nu_{\min}, \nu_{\max} \right]$
	such that $1 < \hat \nu < \rho$ and set $\nu_{i} = \hat \nu$ on each subdomain (in practice we set
	$\hat \nu = \frac{1}{2} \left( \max \left \{1, \nu_{\min} \right \} + \nu_{\max} \right)$). If $\nu_{\max} \leq 1$, we increase the 
	number of subdomains (to increase $\rho$) or use a higher number of Chebyshev coefficients $N_{i}$ 
	and then try again. 
\end{remark}

\subsection{Lotka-Volterra}
\label{sec:LK}
We have proven the existence of connecting orbits from 
$\left(b,0,1-b,0\right)$ to $\left(1,0,0,0\right)$ in \eqref{eq:LK_ODE} for 
$a = 5, \ b = \frac{1}{2}$, $D = 3$ and different values of $\kappa$. Recall that 
these orbits correspond to traveling fronts of \eqref{eq:LK} with wave speed $\kappa$. 
The choices for these parameter values were somewhat arbitrary and were obtained by experimenting 
with the parameter values considered in \cite{MR760975}. In particular, we chose the parameters in such a way that
the stable eigenvalues associated to $(1,0,0,0)$ consisted of one complex conjugate pair of eigenvalues and one
real eigenvalue. 

\begin{figure}[tp] 
	\centering
	\subfloat[$\kappa = -1$]{\label{fig:LK_chebcoeffs}\includegraphics[width=0.5\textwidth]{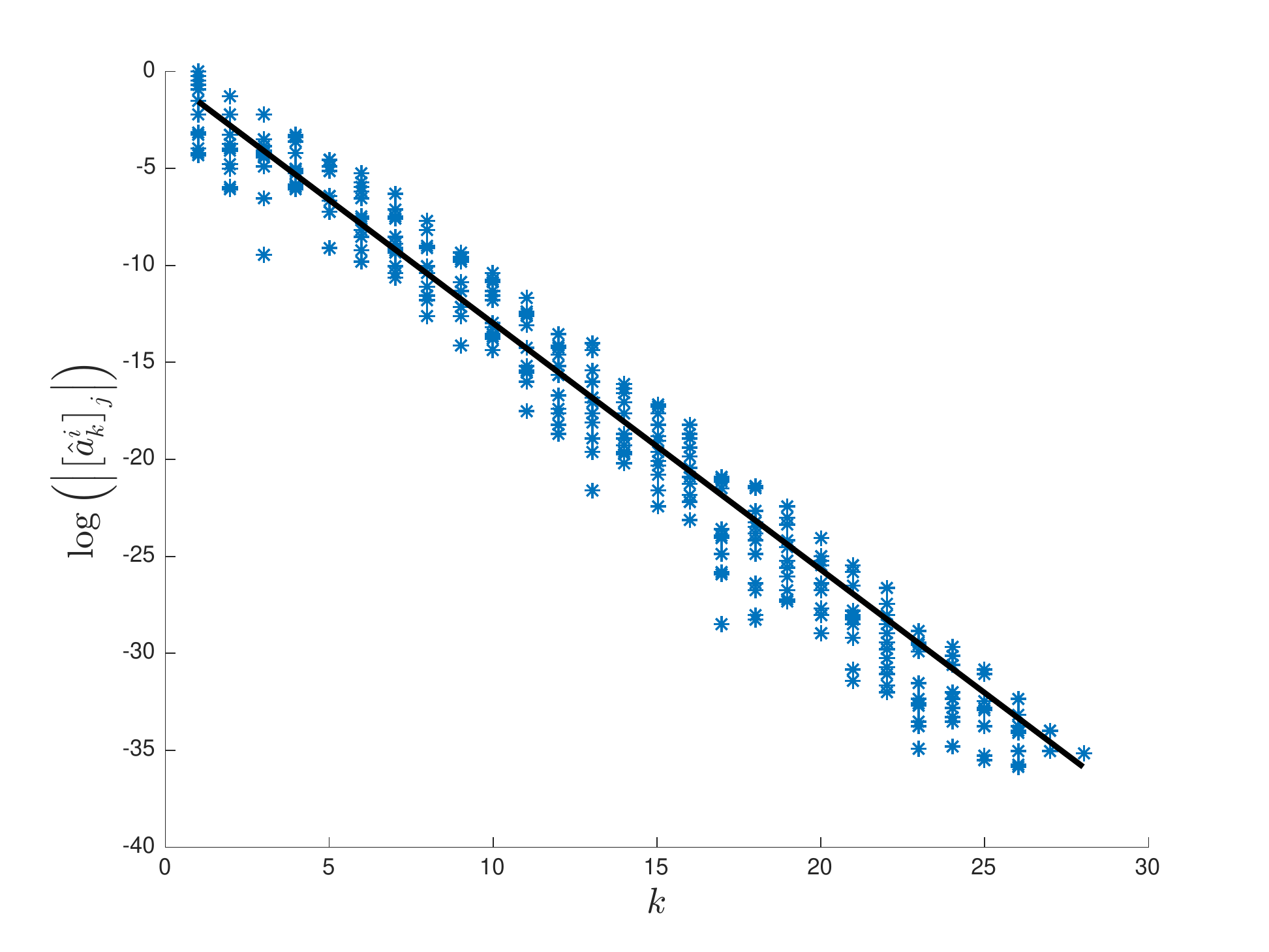}} 
	\subfloat[$\kappa \approx -0.7861$]{\label{fig:LK_a_decay_cont}\includegraphics[width=0.5\textwidth]{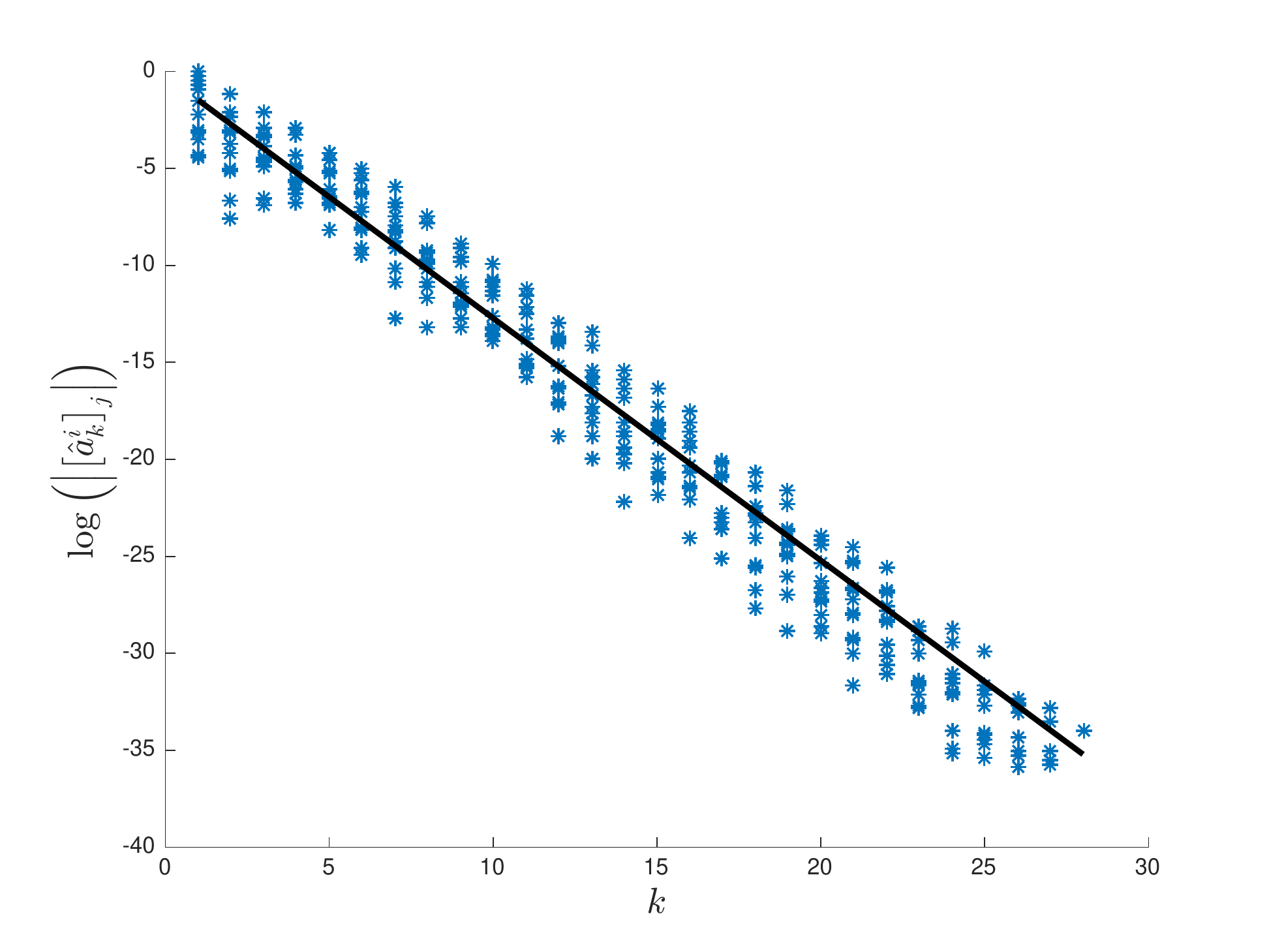}}
	\caption{\label{fig:LK_a_decay} 
		      A semi-logarithmic plot of the (nonzero) Chebyshev coefficients of the connecting orbits at
		      $\kappa \in \left\{-1, -0.7861 \right \}$ on all subdomains for all four components. The black lines
		      correspond to the best line through the points 
		      $ \left \{ \left(k, \log \left \vert \left[ \hat a^{i}_{k} \right]_{j} \right \vert \right) : \
		       \left \vert \left[ \hat a^{i}_{k} \right]_{j} \right \vert \geq 10^{-16}, \ 
		       0 \leq k \leq N_{i} - 1, \ 1 \leq i \leq m, \ 1 \leq j \leq n \right \}$.
		      The results show that the decay rates of the Chebyshev coefficients remained roughly the same  
		      for $\kappa \in \left[-1, -0.7861 \right]$.}
\end{figure}

\paragraph{Connecting orbit at $\kappa=-1$}
We started with a numerical approximation of a connecting orbit at $\kappa=-1$ and used the steps outlined in the previous
section to obtain the following computational parameters: 
\begin{itemize}
	\item \emph{Parameterization mappings:} we used $K^{u} = \begin{bmatrix} 13 & 13 \end{bmatrix}$ and 
		$K^{s} = \begin{bmatrix} 9 & 9 & 9 \end{bmatrix}$ Taylor coefficients
		for approximating the local (un)stable manifolds. The length of the stable and unstable eigenvectors was 
		set to $\epsilon_{u,e_{k}} = 0.0565$ and $\epsilon_{s,e_{k}} = 0.0635$, respectively.
		The truncation parameters and the scalings of the eigenvectors were obtained via the procedure 
		in Remark \ref{remark:manifold_parameters}. 
		The scalings of the eigenvectors were relatively small, since the procedure in Remark \ref{remark:manifold_parameters} was 
		designed to use as little Taylor coefficients as possible to ensure that validation is feasible. If we would allow for larger 
		truncation parameters, the scalings of the eigenvectors (and hence the ``size'' of the charts on the local (un)stable manifolds)
		could be increased substantially. However, since it is computationally cheaper to increase the integration time in comparison 
		to increasing the truncation parameters for the (un)stable manifolds, 
		we have chosen to keep the truncation parameters $K^{u}$ and (especially) $K^{s}$ small. 
	\item {\emph{Chebyshev approximations}}: 
		we used $m=3$ subdomains and $N = \begin{bmatrix} 50 &47& 50 \end{bmatrix}$ Chebyshev modes. 
		The integration time was set to $L=15$. 
		The Chebyshev coefficients are shown in Figure \ref{fig:LK_chebcoeffs}. This figure 
		shows that the decay rates of the Chebyshev coefficients were approximately the same on each subdomain
		(hence the domain decomposition was successful). 
	\item {\emph{Validation parameters}}: we used $\nu_{i} = 1.1967$ on each subdomain. This value
		was obtained from the procedure in Remark \ref{remark:nu_i}.
\end{itemize} 
The dimension of the Galerkin projection was $\dim \left( \mathcal{X}^{NK} \right) = 5382$. With the above choices for the
computational parameters, we successfully validated a connecting orbit at $\kappa =-1$ and proved that
the radii-polynomials were negative for $r \in \left[8.6070 \cdot 10^{-11}, r^{\ast} \right]$.
We remark that it is possible to validate the connecting orbit with a smaller number of Chebyshev coefficients as well. The reason why we used
more Chebyshev coefficients than strictly necessary was to get the bounds $Z^{1}_{\Pi^{ij}_{a}}$ as small as possible with an eye towards future work, namely in order to make
continuation with large step-sizes feasible. 

\begin{figure}[tp]
	\centering
	\subfloat[$\kappa \in {\left[-1, -0.7861 \right]}$]{\label{fig:LK_bif_1}\includegraphics[width=0.5\textwidth]{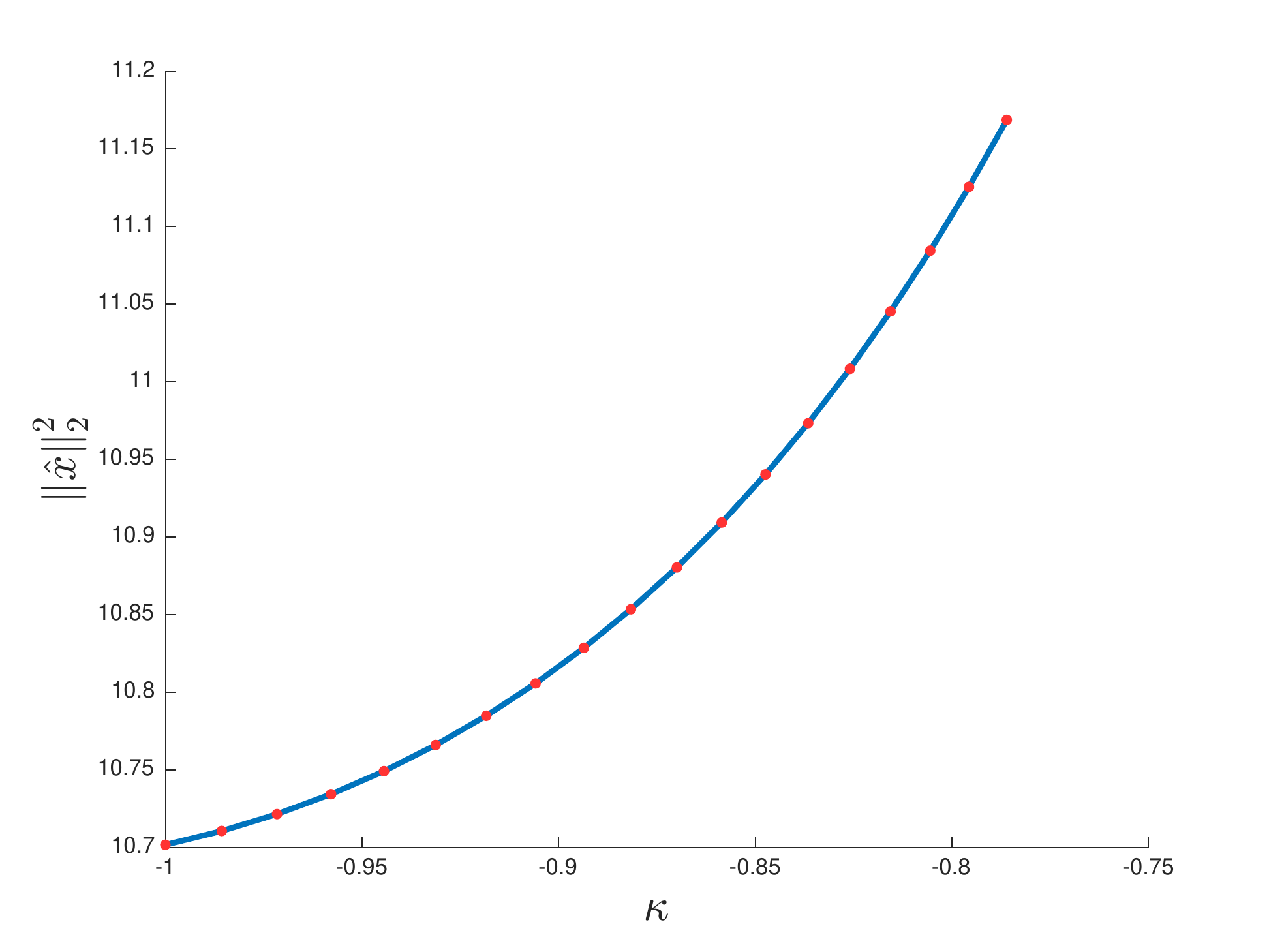}}
	\subfloat[$\kappa \in {\left[-0.7767, -0.5938 \right]}$]{\label{fig:LK_bif_2}\includegraphics[width=0.5\textwidth]{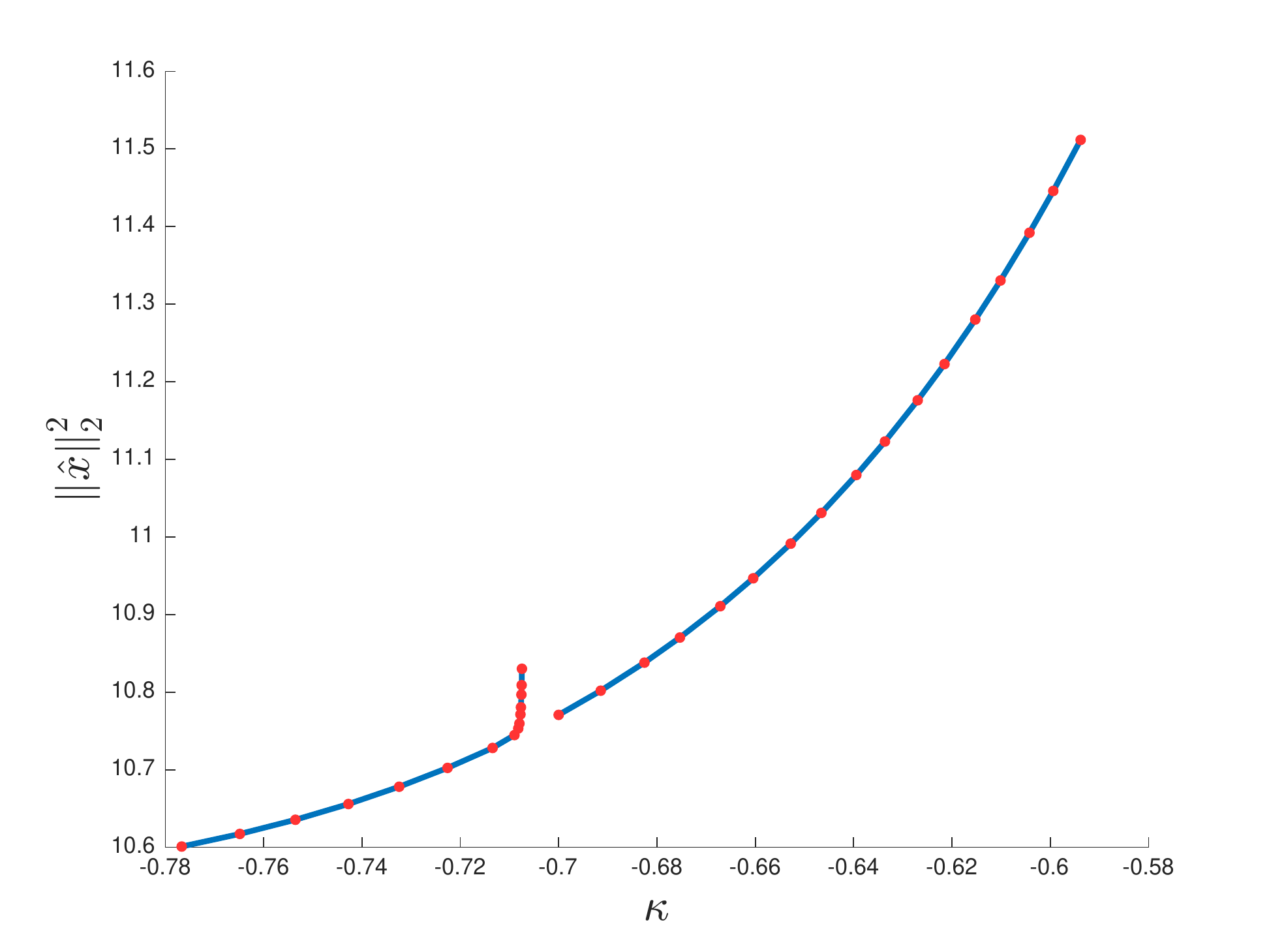}} 
	\caption[a bla]{\label{fig:LK_bif} 
		      A bifurcation diagram obtained by continuing the connecting orbit at $\kappa = -1$. The bifurcation
		      curves were computed by performing (non-rigorous) pseudo-arc length continuation in the parameter 
		      $\kappa$. The red points on the curves correspond to validated connecting orbits.
		      In all cases the validation radius $\hat r$ was bounded by $2.5347 \cdot 10^{-9}$.   
		      The curve in Figure \ref{fig:LK_bif_1} was computed with truncation parameters
		      $K^{u} = \begin{bmatrix} 13 & 13 \end{bmatrix}$, $K^{s} = \begin{bmatrix} 9 & 9 & 9 \end{bmatrix}$ 
		      and $N = \begin{bmatrix} 50 &47& 50 \end{bmatrix}$. The connecting orbits were validated 
		     by using $\nu_{i} = 1.1967$ on each subdomain. The curve in Figure \ref{fig:LK_bif_2}
		      was computed with truncation parameters $K^{u} = \begin{bmatrix} 13 & 13 \end{bmatrix}$, 
		      $K^{s} = \begin{bmatrix} 12 & 12 & 12 \end{bmatrix}$ and $N = \begin{bmatrix} 55 &52& 62 \end{bmatrix}$. 
		      The connecting orbits were validated by using $\nu_{i} =1.1627$ on each subdomain. The ``gap''  
		      at $\kappa \approx -0.7071$ corresponds to a bifurcation caused by the presence of a resonance
		      at $\kappa = -\frac{1}{2} \sqrt{2}$.
}
\end{figure}

\paragraph{Discrete continuation}
Next, we continued the connecting orbit at $\kappa =-1$ by performing pseudo-arc length continuation. At each
continuation step, we tried to validate the orbit with the same computational parameters. We succeeded in validating  a family
of connecting orbits in this way for a finite number of wave speeds $\kappa \in \left[-1, -0.7861 \right ]$,
see Figure \ref{fig:LK_bif_1}. 
The reader is referred to the code for the exact parameter values $\kappa$ at which the connecting orbits were validated.
Here we only give rounded values of $\kappa$ using four decimal places. 

\begin{figure}[tp]
	\centering
	\subfloat[$\hat \lambda^{u}$]{\label{fig:LK_eig_u}\includegraphics[width=0.5\textwidth]{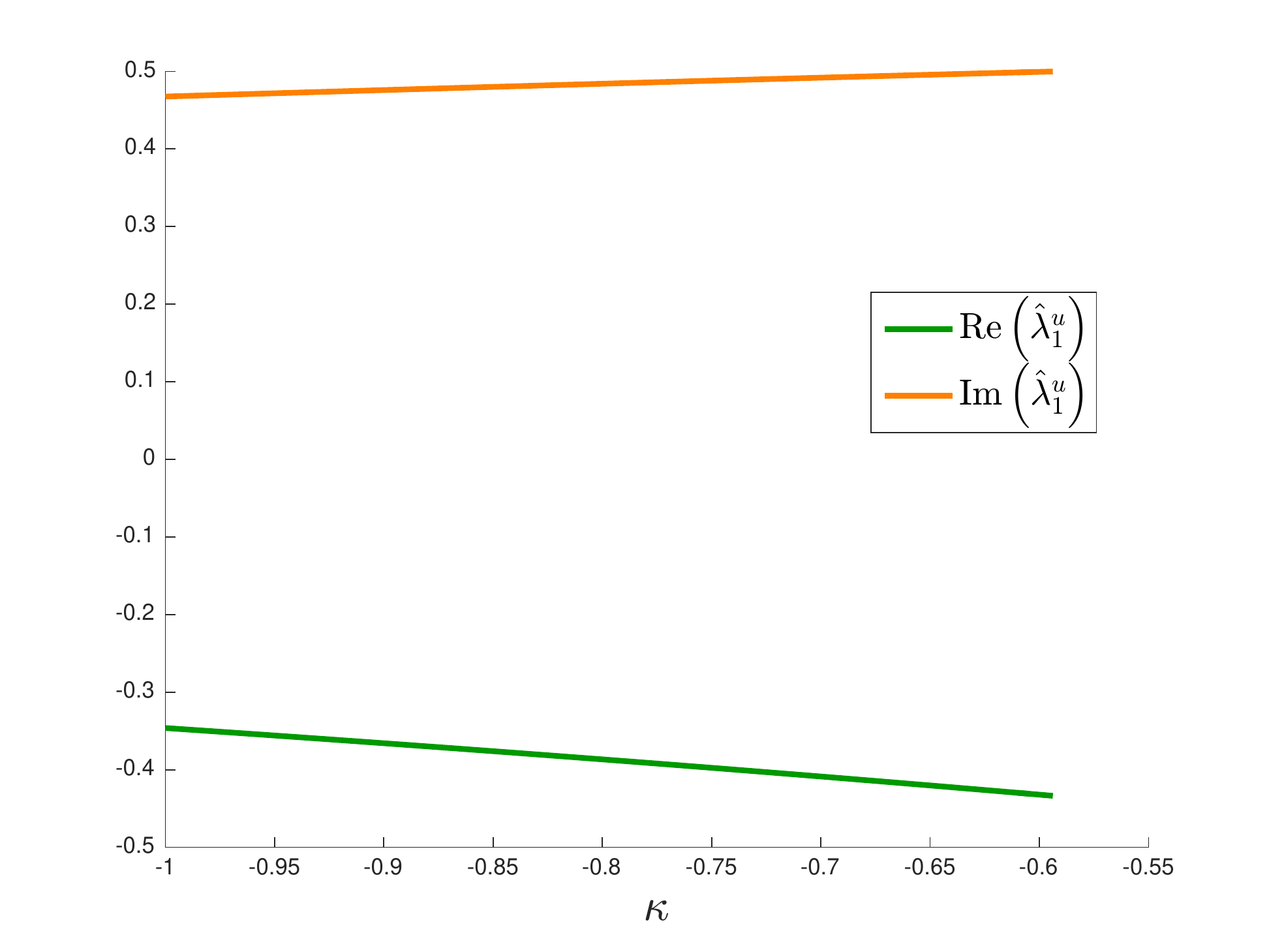}}
	\subfloat[$\hat \lambda^{s}$]{\label{fig:LK_eig_s}\includegraphics[width=0.5\textwidth]{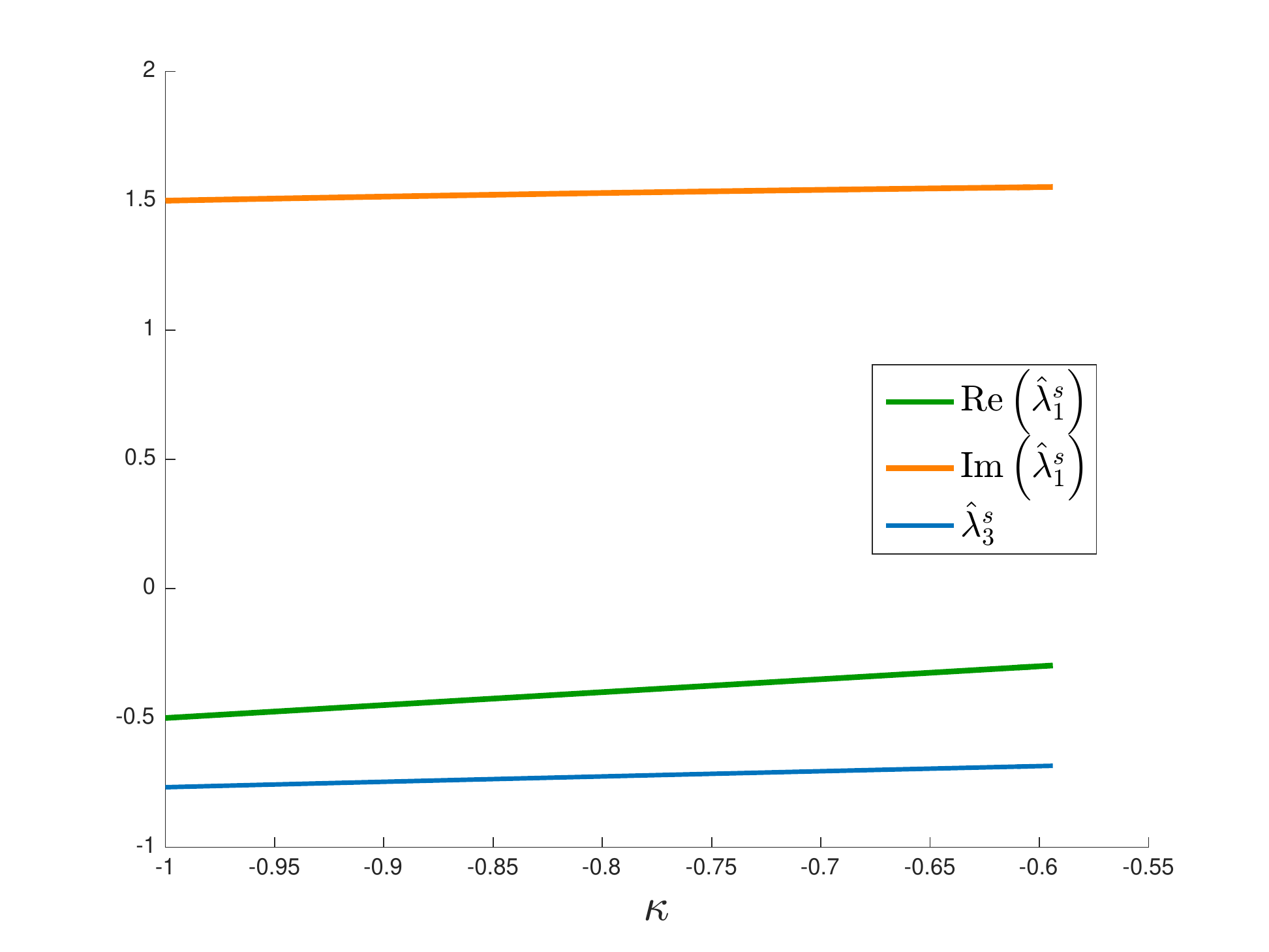}} 
	\caption[a bla]{\label{fig:LK_eig} The dependence of the stable and unstable eigenvalues on $\kappa$. 
}
\end{figure}

We were not able to validate the connecting orbit at the next continuation step $\kappa  \approx -0.7767$ with the same computational parameters. The reason for this was that the bound 
for $\Pi_{Q}^{4} \left( DF \left( \hat x \right) - \widehat{DF} \right)$ became too large, which 
was related to the fact that the decay rates of $\hat q$ decreased as $\kappa$ increased. 
In addition, the real part of the stable eigenvalues decreased as well
when $\kappa$ increased (see Figure \ref{fig:LK_eig}), which contributed to the deterioration of the bounds for the stable manifold. 
On the other hand, the bounds for $\Pi_{a}^{ij} \left( DF \left( \hat x \right) - \widehat{DF} \right)$ did not deteriorate at all during
the continuation process. The main reason for this is that the shape and time parameterization of the orbit remained roughly ``the same'' throughout
the continuation procedure. This caused the decay rates of the 
Chebyshev coefficients to remain roughly the same as well, as shown in Figure \ref{fig:LK_a_decay}. 

To validate connecting orbits for $\kappa \geq -0.7767$, we recomputed the parameterization of the local stable manifold as explained in  
Remark \ref{remark:manifold_parameters}. This resulted in a parameterization with $K^{s} = \begin{bmatrix} 12 & 12 &12 \end{bmatrix}$ 
Taylor coefficients. Furthermore, the length of the stable eigenvectors was set to $\epsilon_{s,e_{k}} = 0.02275$. The resulting parameterized part of the stable manifold
was  significantly smaller due to the new scaling of the eigenvectors. To ensure that the endpoint of the connecting orbit was
contained in the smaller chart, we integrated the connecting orbit forward in time (while keeping the ``initial'' starting point on the unstable manifold fixed)
and increased the integration time to $L=20$. 
We then used the procedures from the previous section again to refine the Chebyshev approximations. In particular, we used $m=3$ subdomains, 
$N = \begin{bmatrix} 55 &52& 62 \end{bmatrix}$ Chebyshev modes and used $\nu_{i} = 1.1627$ on each subdomain. 
The dimension of the Galerkin projection was $\dim \left( \mathcal{X}^{NK} \right) = 10258$. With these parameter values, we were able to successfully 
validate a finite number of connecting orbits for $\kappa \in \left[-0.7767, -0.7075\right]$, see Figure \ref{fig:LK_bif_2}. 

Figure \ref{fig:LK_bif_2} shows that there is a bifurcation at $\kappa \approx -0.7071$. 
To understand what caused this bifurcation, we consider the (approximate) stable eigenvalues 
$\hat \lambda^{s}$ at $\kappa \approx -0.7075$:
\begin{align*}
	\hat \lambda^{s} \approx 
	\begin{bmatrix}
		-0.3537 + 1.541i \\
		-0.3537 - 1.541i \\
		- 0.7072
	\end{bmatrix}.
\end{align*}
Note that $\hat \lambda^{s}_{1} + \hat \lambda^{s}_{2} \approx \hat \lambda^{s}_{3}$. This provides numerical evidence for the presence of an eigenvalue
resonance at $\kappa \approx -0.7071$ and explains the observed bifurcation. In this relatively simple case,
one can prove with pen and paper that there is in fact a resonance at $\kappa = -\frac{1}{2} \sqrt{2}$. 
Hence, in order to validate connecting orbits near $\kappa = -\frac{1}{2} \sqrt{2}$ (and in particular at the resonance point itself), 
one needs to modify the mapping $F_{Q}$ by conjugating to a nonlinear normal form (instead of just the linear one) as explained in \cite{ManifoldTaylor}. 
Here we do not pursue this issue any further and leave it as a future research project. 

\begin{figure}[tb]
	\begin{center}
		\includegraphics[scale = 0.55]{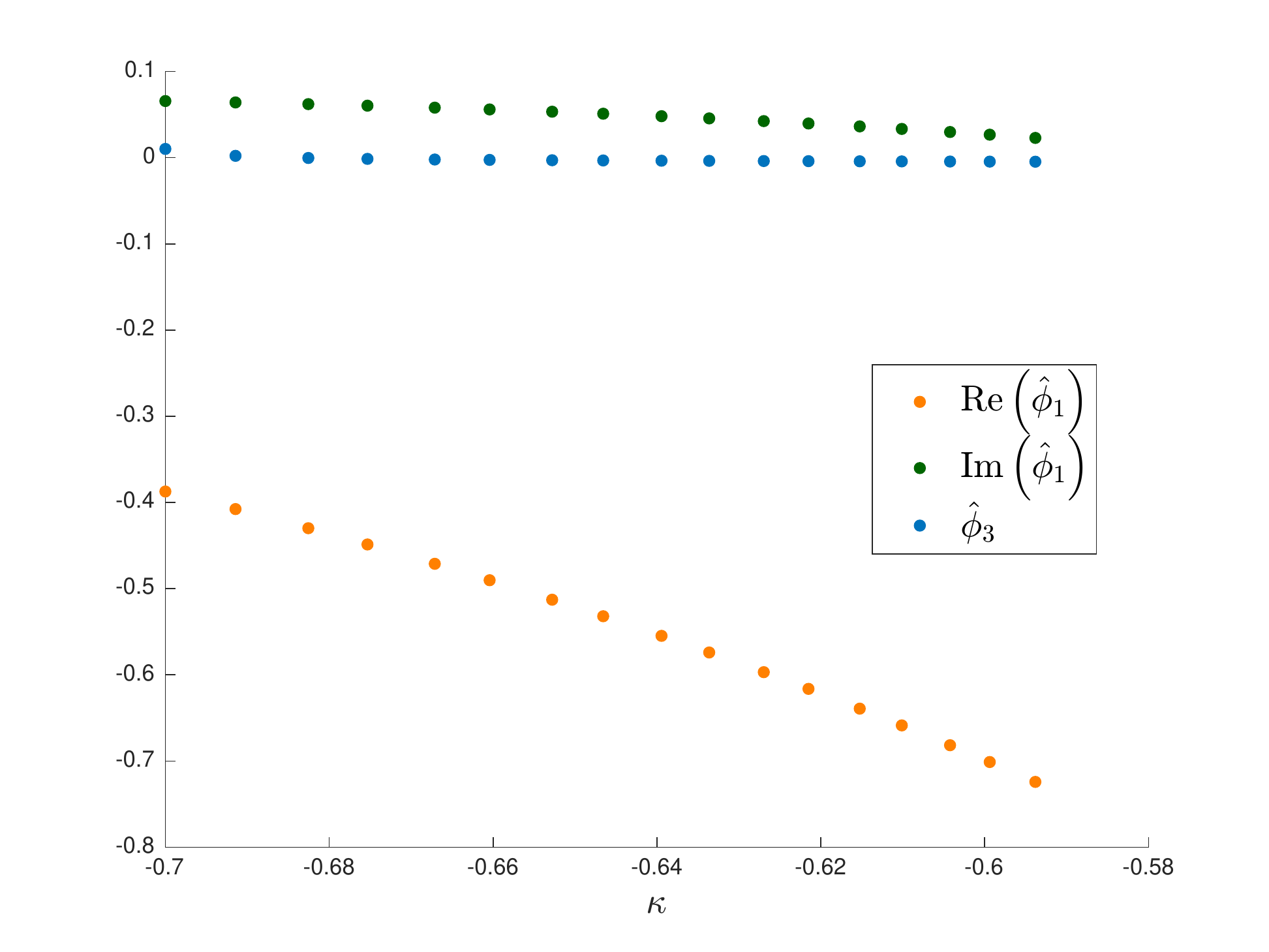}
	\end{center}
	\caption{\label{fig:LK_phi} The dependence of the stable parameterization variables $\hat \phi$ for $\kappa \in \left[-0.7,-0.5938\right]$. 
	The results show that $\left \vert \text{Re} \left( \hat \phi_{1} \right) \right \vert = \left \vert \text{Re} \left( \hat \phi_{2} \right) \right \vert $ increased as $\kappa$ increased
	and eventually became too large (in the sense that the bounds related to the equation $u(1) = Q \left( \phi \right)$ were too large). This issue can be resolved by either integrating the connecting orbit 
	further in time (increasing $L$) or by enlarging the chart on the local stable manifold in the ``directions'' of $\phi_{1}$ and $\phi_{2}$.  
}
\end{figure}

To continue the connecting orbit past the resonance, we set $\kappa = -0.7$ and then continued further from this point. We succeeded in validating connecting orbits
for $\kappa \in \left[-0.7, -0.5938\right]$ without changing the computational parameters, see Figure \ref{fig:LK_bif_2}. 
As $\kappa$ increased, the ``size'' of the chart on the local stable manifold kept decreasing. 
As a consequence, for $ \kappa > -0.5938$, the endpoint of the connecting orbit was too far away from the stable equilibrium 
in the sense that the bounds related to the equation $u(1) = Q \left( \phi \right)$ were too large (also see Figure \ref{fig:LK_phi}). 
Although we did not continue any further, we remark that validation for $\kappa > -0.5938$ is feasible by increasing the integration time and recomputing 
the Chebyshev approximations as before.

An interesting future research project would be to develop algorithms for automatically detecting when the integration time and/or manifolds need to be
modified. We believe that the heuristics in Remarks \ref{remark:manifold_parameters} and \ref{remark:nu_i} would be a good starting point for
developing such algorithms. 

\FloatBarrier

\subsection{Traveling fronts in a fourth order parabolic PDE}
\label{sec:Fourth_order}
We have proven the existence of connecting orbits in \eqref{eq:Fourth_order_ODE} from $(-1,0,0,0)$ to 
$(a,0,0,0)$ for $a = -0.1$, $\kappa = -2$ and various values of $\gamma$.
Recall that these orbits correspond to traveling fronts of \eqref{eq:Fourth_order_PDE} with wave speed $\kappa$. 
The values for $a$ 
and the wave speed $\kappa$ were chosen through experimentation. 
We started with a connecting orbit at $\gamma \approx 0.4557$ (see the code for exact parameter value) 
and then used the procedure as explained at the beginning
of this section to select the following computational parameters:
\begin{itemize}
	\item \emph{Parameterization mappings:} we used $K^{u} = \begin{bmatrix} 15 & 15 \end{bmatrix}$ and 
		$K^{s} = \begin{bmatrix} 12 & 12 & 12 \end{bmatrix}$ Taylor coefficients
		for approximating the local (un)stable manifolds. The length of the stable and unstable eigenvectors was 
		set to $\epsilon_{u,e_{k}} = 5.4476 \cdot 10^{-2}$ and $\epsilon_{s,e_{k}} = 5.3337 \cdot 10^{-3}$, respectively. 
	\item {\emph{Chebyshev approximations}}: 
		we used $m=2$ subdomains and $N = \begin{bmatrix} 62 & 61 \end{bmatrix}$ Chebyshev modes. 
		The integration time was set to $L=30$. The decay rates of the Chebyshev coefficients were approximately 
		the same on each subdomain (hence the domain decomposition was successful). 
	\item {\emph{Validation parameters}}: we used $\nu_{i} = 1.1491$ on each subdomain. 
\end{itemize} 
The dimension of the Galerkin projection was $\dim \left( \mathcal{X}^{NK} \right) = 10314$.
With the above choices for the computational parameters, we were able to validate the connecting orbit at $\gamma \approx 0.4557$
and proved that the radii-polynomials were negative for $r \in \left[4.8332 \cdot 10^{-10}, 2.8332 \cdot 10^{-6} \right]$. 

\begin{figure}[tp]
	\centering
	\subfloat[$\gamma \in {\left[0.4557,1.233 \right]}$]{\label{fig:FK_bif_12}\includegraphics[width=0.33\textwidth]{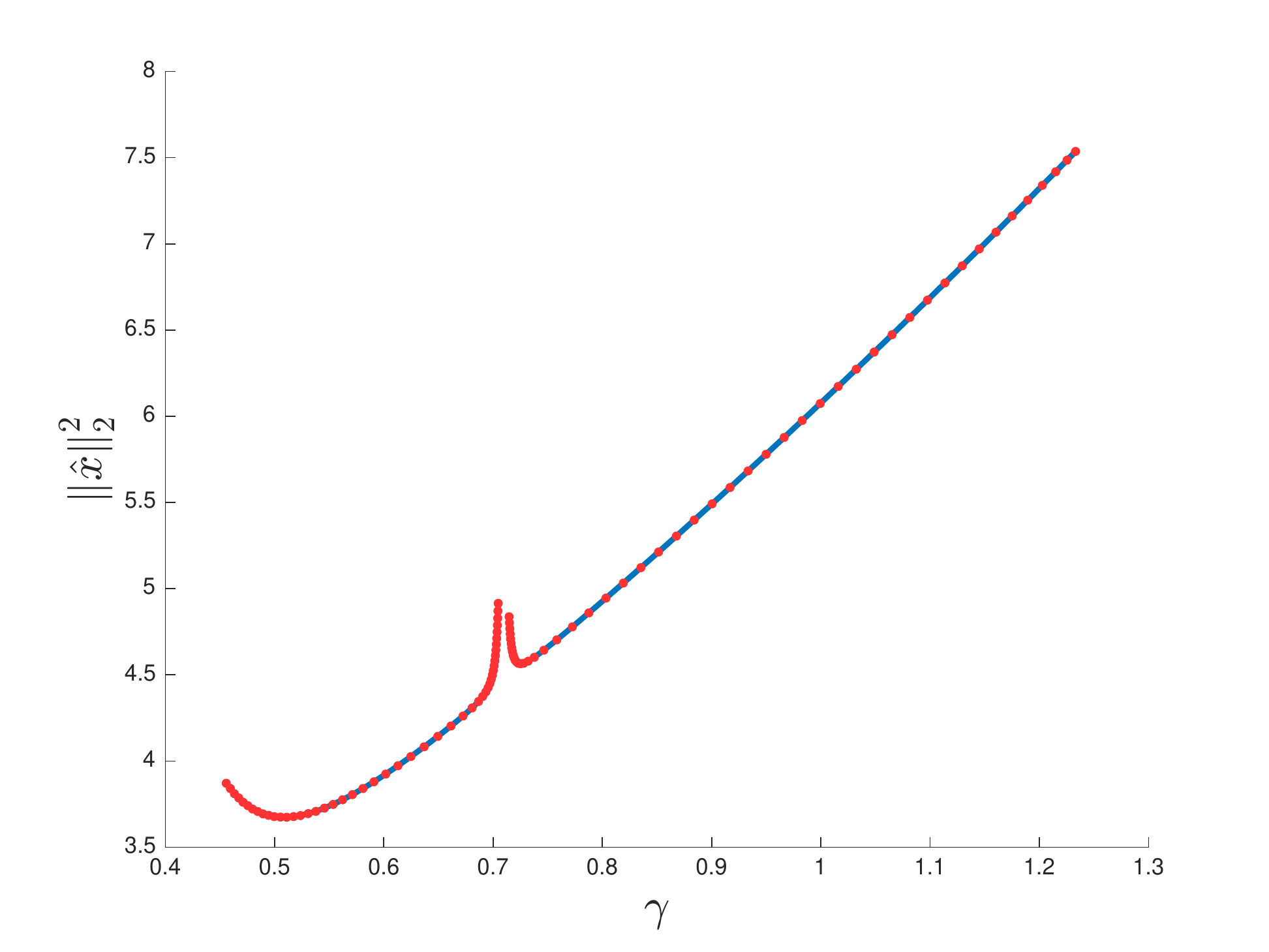}}
	\subfloat[$\gamma \in {\left[1.243, 4.089 \right]}$]{\label{fig:FK_bif_3}\includegraphics[width=0.33\textwidth]{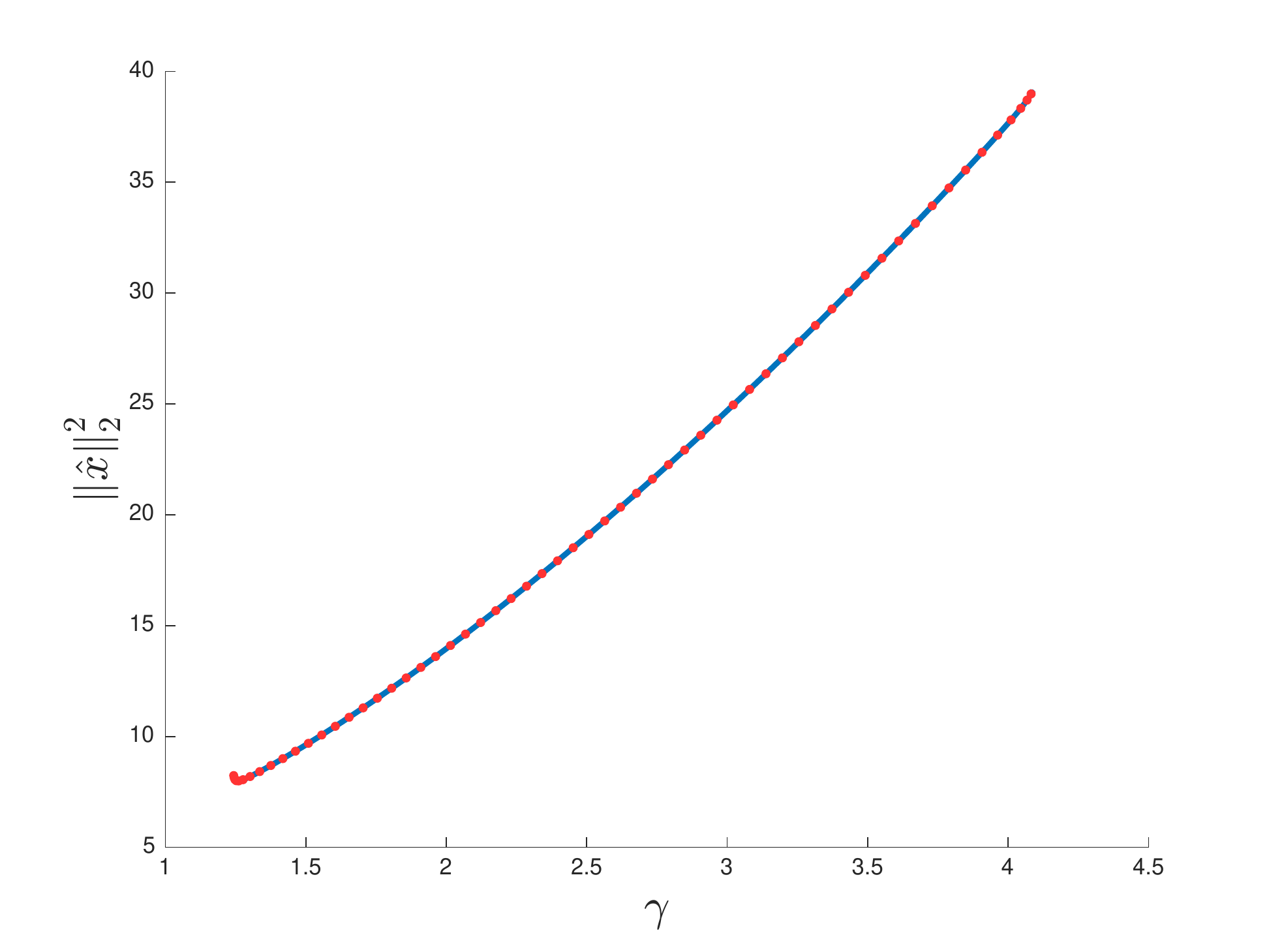}} 
	\subfloat[$\gamma \in {\left[4.202, 10.50 \right]}$]{\label{fig:FK_bif_4}\includegraphics[width=0.33\textwidth]{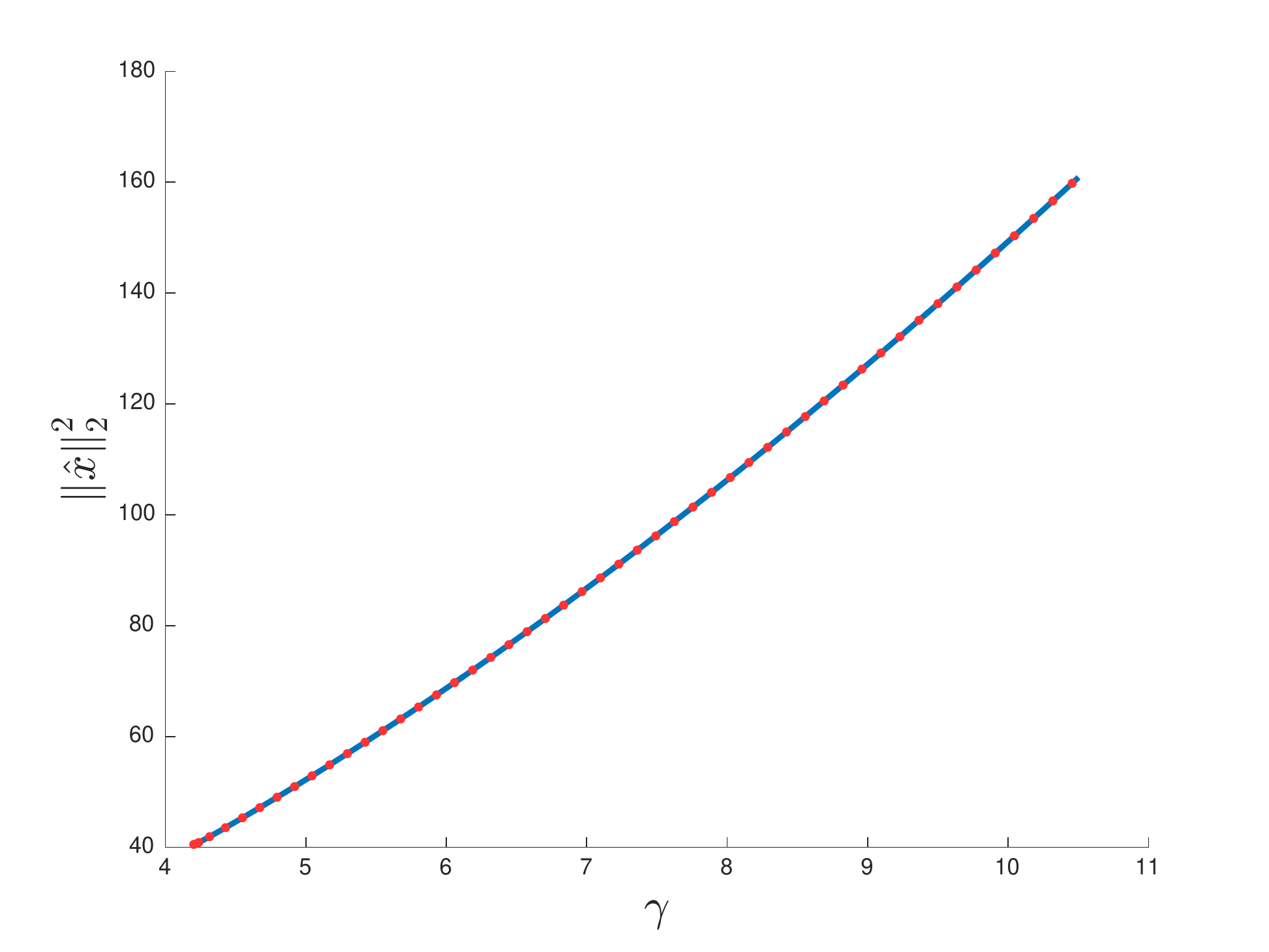}} 
	\caption[a bla]{\label{fig:FK_bif} 
		      Bifurcation diagrams obtained by continuing the connecting orbit at $\gamma \approx 0.4557$.
		      The bifurcation curves were computed by performing (non-rigorous) pseudo-arc length continuation in $\gamma$. 
		      The red points on the curves correspond to validated connecting orbits. In each case, the validation 
		      radii were bounded by: \protect~\subref{fig:FK_bif_12} $\hat r \leq 4.6324 \cdot 10^{-6}$, 
		      \protect~\subref{fig:FK_bif_3} $\hat r \leq 5.3757 \cdot 10^{-6}$, 
		       \protect~\subref{fig:FK_bif_4} $\hat r \leq 9.9475 \cdot 10^{-6}$. 
		      The values of the
		      computational parameters are reported in Table \ref{table:FK_results}. 
}
\end{figure}

\begin{figure}[t]
	\centering
	\subfloat[$\gamma \in {[0.4557,1.233]}$, $L =30$]{\label{fig:FK_u_gamma_12}\includegraphics[width=0.333\textwidth]{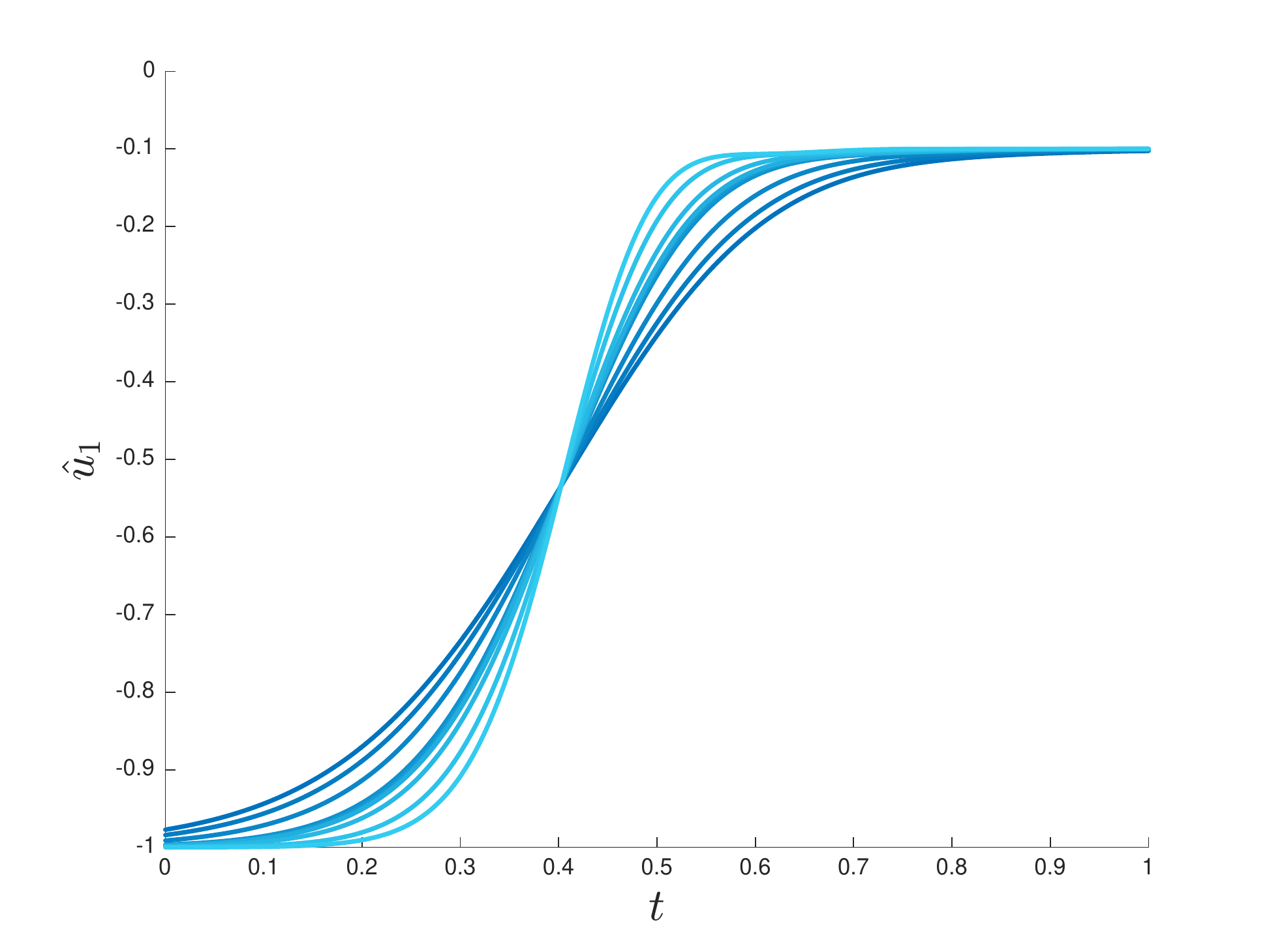}}
	\subfloat[$\gamma \in {[1.243,4.089]}$, $L=8$]{\label{fig:FK_u_gamma_3}\includegraphics[width=0.333 \textwidth]{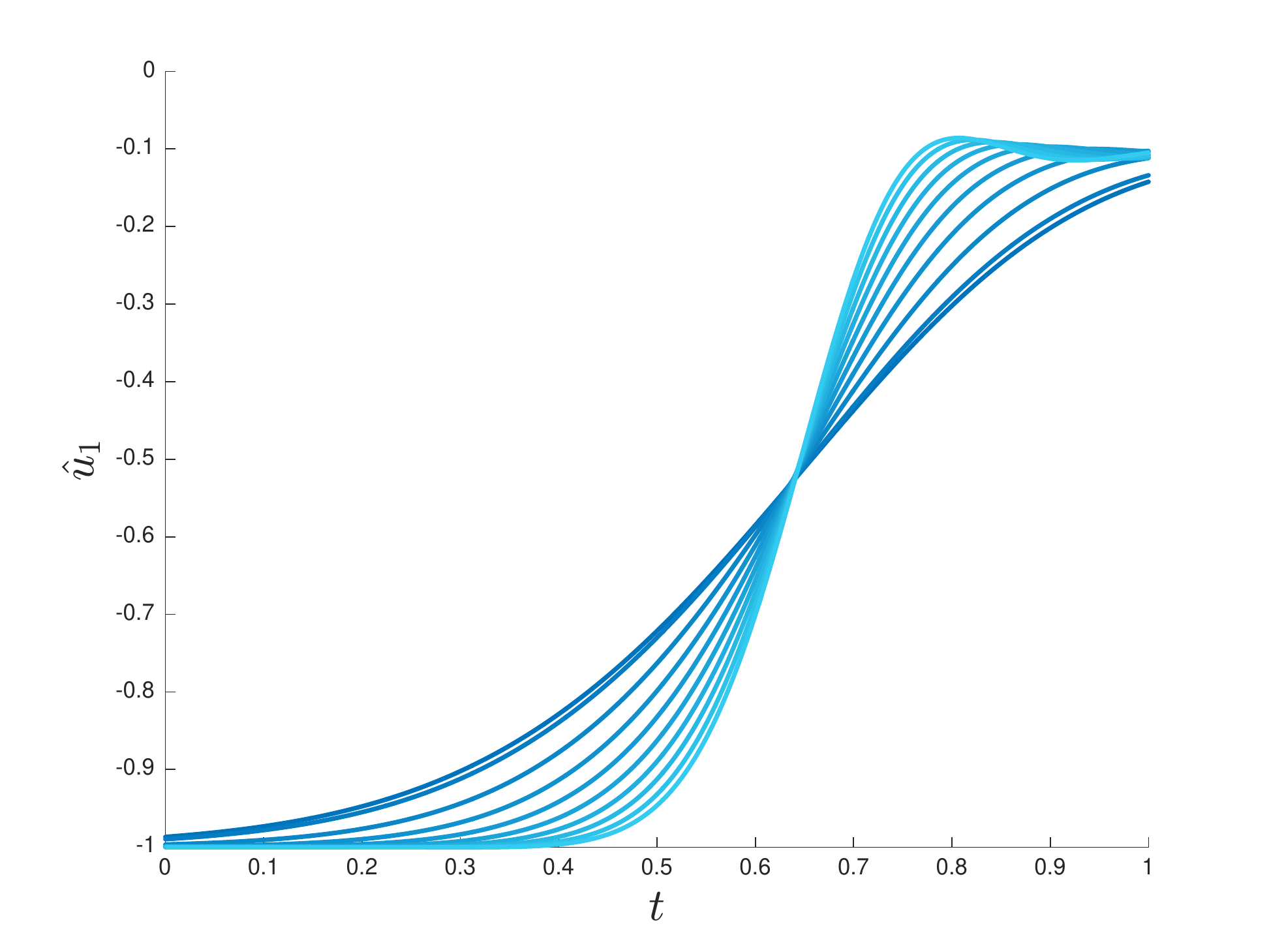}} 
	\subfloat[$\gamma \in {[4.202,10.50]}$, $L=4$]{\label{fig:FK_u_gamma_4}\includegraphics[width=0.333 \textwidth]{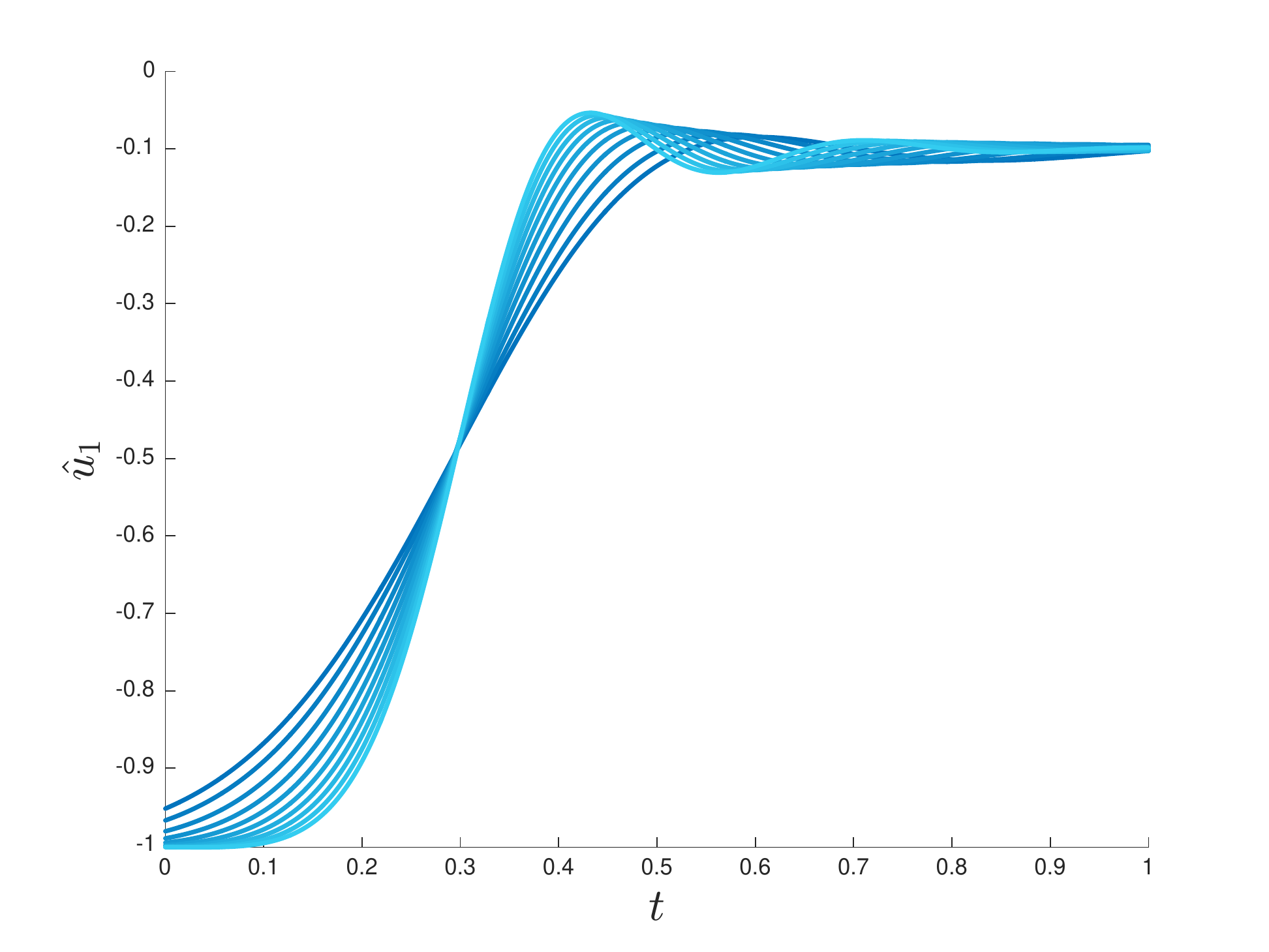}} 
	\caption[a bla]{\label{fig:FK_u_gamma} 
	The first component $u_{1}$, which corresponds to a traveling wave profile of~\eqref{eq:Fourth_order_PDE}, of the validated connecting orbits
	for $\gamma \in \left[0.4557, 10.50 \right]$. In each case, for $\gamma$ close to the left endpoint of the interval, we
	have colored the associated orbits in dark blue. As $\gamma$ increased, we used increasingly lighter shades
	of blue. Note the oscillations for larger values of $\gamma$. 
}
\end{figure}

Next, we performed (non-rigorous) pseudo-arclength continuation and tried to validate the orbits at each discrete
continuation step by using the same computational parameters. If the validation failed at a particular continuation step, 
we determined the cause (as in the previous section) and resolved the issue by modifying the computational parameters. In addition, 
we also checked in the case of failure whether the dimension of the Galerkin-projection could be significantly reduced by 
decreasing the integration time and the number of Chebyshev or Taylor coefficients (with special emphasis on reducing
the number of Taylor coefficients associated to the stable manifold). 
The results are summarized in Table \ref{table:FK_results}. The corresponding bifurcation curves are shown 
in Figure \ref{fig:FK_bif} and the associated traveling wave profiles are depicted in Figure \ref{fig:FK_u_gamma}. 

\begin{table}[tp]
	\begin{center}
	\renewcommand{\arraystretch}{1.2}
	\begin{tabular}{|c|c|c|c|c|c|c|c|c|} 
		\hline
		$\gamma$ 				& 				$L$ 					& 				$N$ 												&
		$K^{s}$ 					& 	 			$\epsilon_{u,k}$ 		& 				$\epsilon_{s,k}$ 									& 	
		$\nu_{i}$					&				Obstruction														 \\
		\hline	
		$\left[0.4557,0.7046 \right]$ 	& 				$30$ 				& 				$\left[ 62 \quad 61  \right]$ 							&
		$12$ 					& 	 			$5.4 \cdot 10^{-2}$ 		& 				$5.3 \cdot 10^{-3}$ 									&	 	
		$1.15$					&				$\hat \lambda_{1}^{u} \approx 2 \hat \lambda^{u}_{2}$														 \\
		\hline	
		$\left[0.7146, 1.233 \right]$ 	& 				$30$ 				& 				$\left[ 62 \quad 61  \right]$  							&
		$9$ 						& 	 			$5.4 \cdot 10^{-2}$ 		& 				$5.3 \cdot 10^{-3}$ 									&	 	
		$1.15$					&				$\hat \lambda_{1}^{u} \approx \hat \lambda^{u}_{2}$													 \\
		\hline	
		$\left[1.243, 4.089 \right]$ 		& 				$8$ 					& 				$63$ 											&
		$9$ 						& 	 			$0.12$ 				& 				$0.17$ 											&	 	
		$1.11$					&				$\hat \lambda_{1}^{s} + \hat \lambda^{s}_{2} \approx \hat \lambda^{s}_{3}$							\\
		\hline	
		$\left[4.202, 10.50 \right]$ 	& 				$4$ 					& 				$51$ 											&
		$7$ 						& 	 			$0.17$ 								& 				$0.12$ 											&	 	
		$1.15$					&				$\hat r > r^{\ast}$							\\
		\hline									
		\end{tabular}
	\end{center}
\caption[]{The (approximate) computational parameters used to validate connecting orbits 
	      for $\gamma \in \left[0.4557,10.50 \right]$. Each row in the table
	      corresponds to an interval on which we performed (non-rigorous) pseudo-arclength
	      continuation and validated rigorously a finite number of connecting orbits with the same
	      computational parameters. In particular, we used $K^{u} = \begin{bmatrix} 15 & 15 \end{bmatrix}$
	      on each interval (though validation with less Taylor coefficients is feasible). 
	      We were not able to validate connecting orbits past
	      the right endpoints of the intervals without modifying the computational parameters. 
	      In each case, we have indicated the obstruction for validating connecting orbits 
	      near the (right) endpoint of the interval.}
\label{table:FK_results}
\end{table}

\begin{figure}[t]
	\centering
	\subfloat[$\hat \lambda^{u}$, $\gamma \in {\left[ 0.4557, 1.233 \right]}$]{\label{fig:FK_eig_u_12}\includegraphics[width=0.333\textwidth]{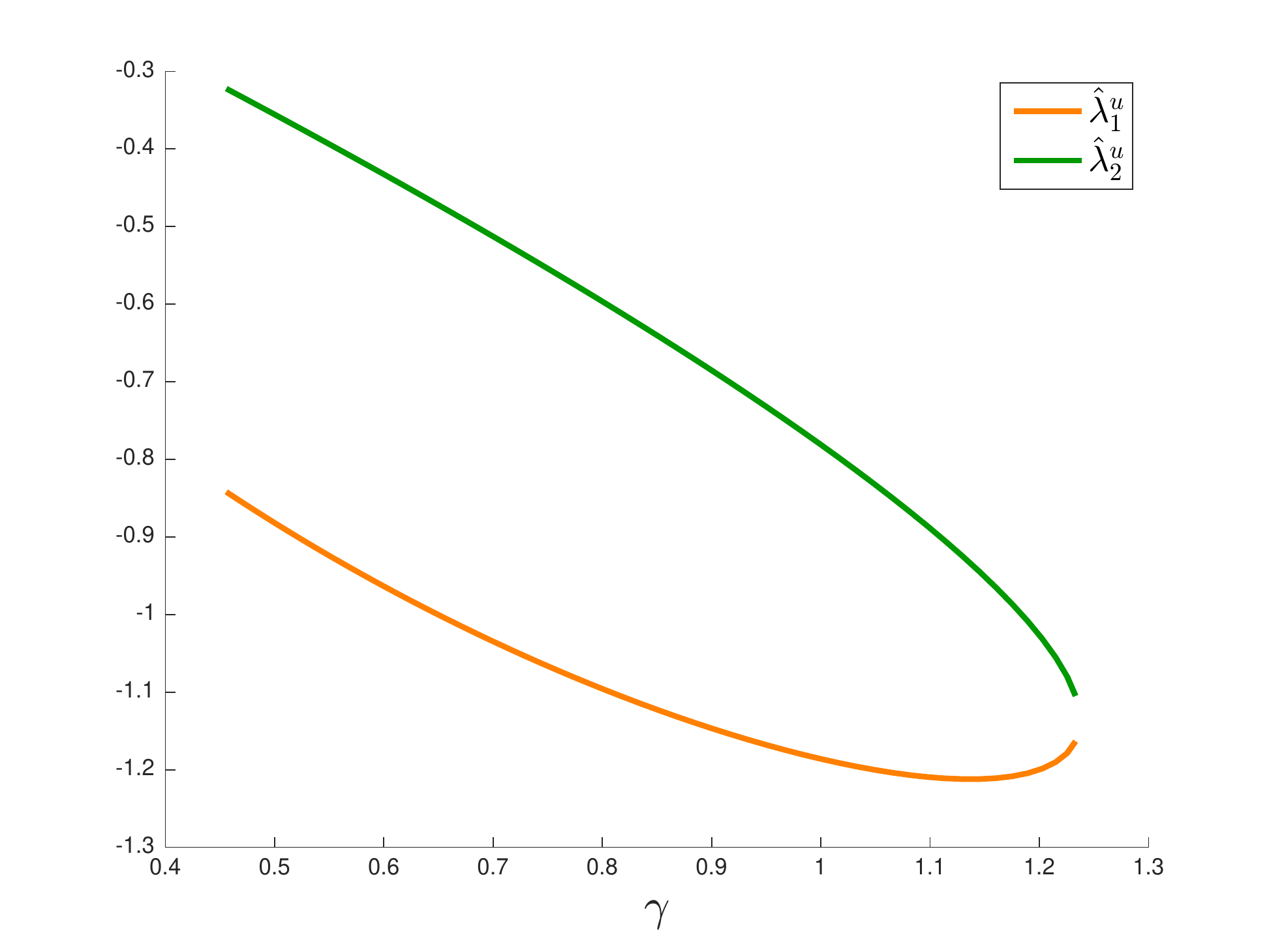}}
	\subfloat[$\hat \lambda^{u}$, $\gamma \in {\left[ 1.243, 10.50 \right]}$]{\label{fig:FK_eig_u_34}\includegraphics[width=0.333 \textwidth]{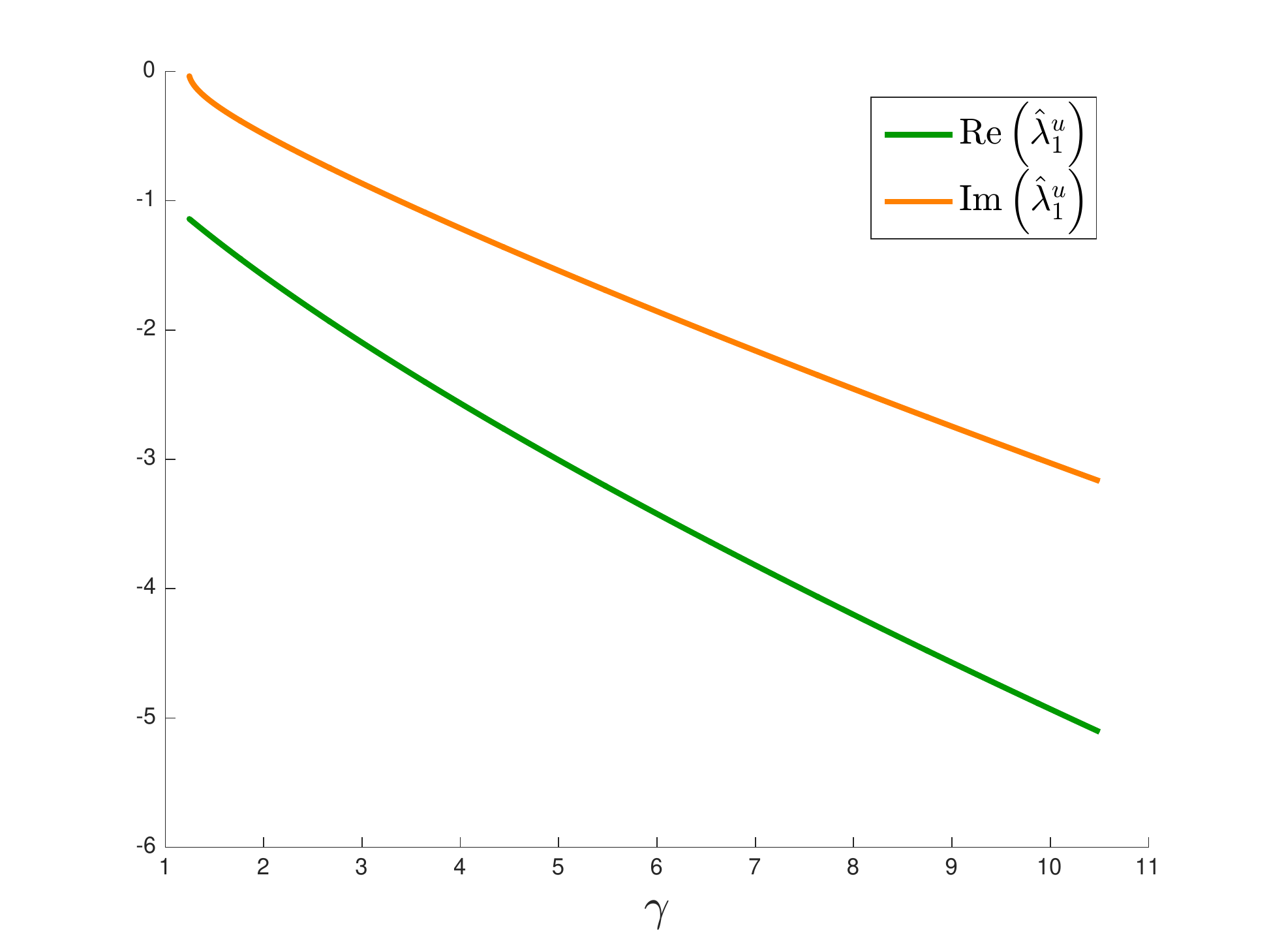}} 
	\subfloat[$\hat \lambda^{s}$, $\gamma \in {\left[0.4557, 10.50 \right]}$]{\label{fig:FK_eig_s}\includegraphics[width=0.333 \textwidth]{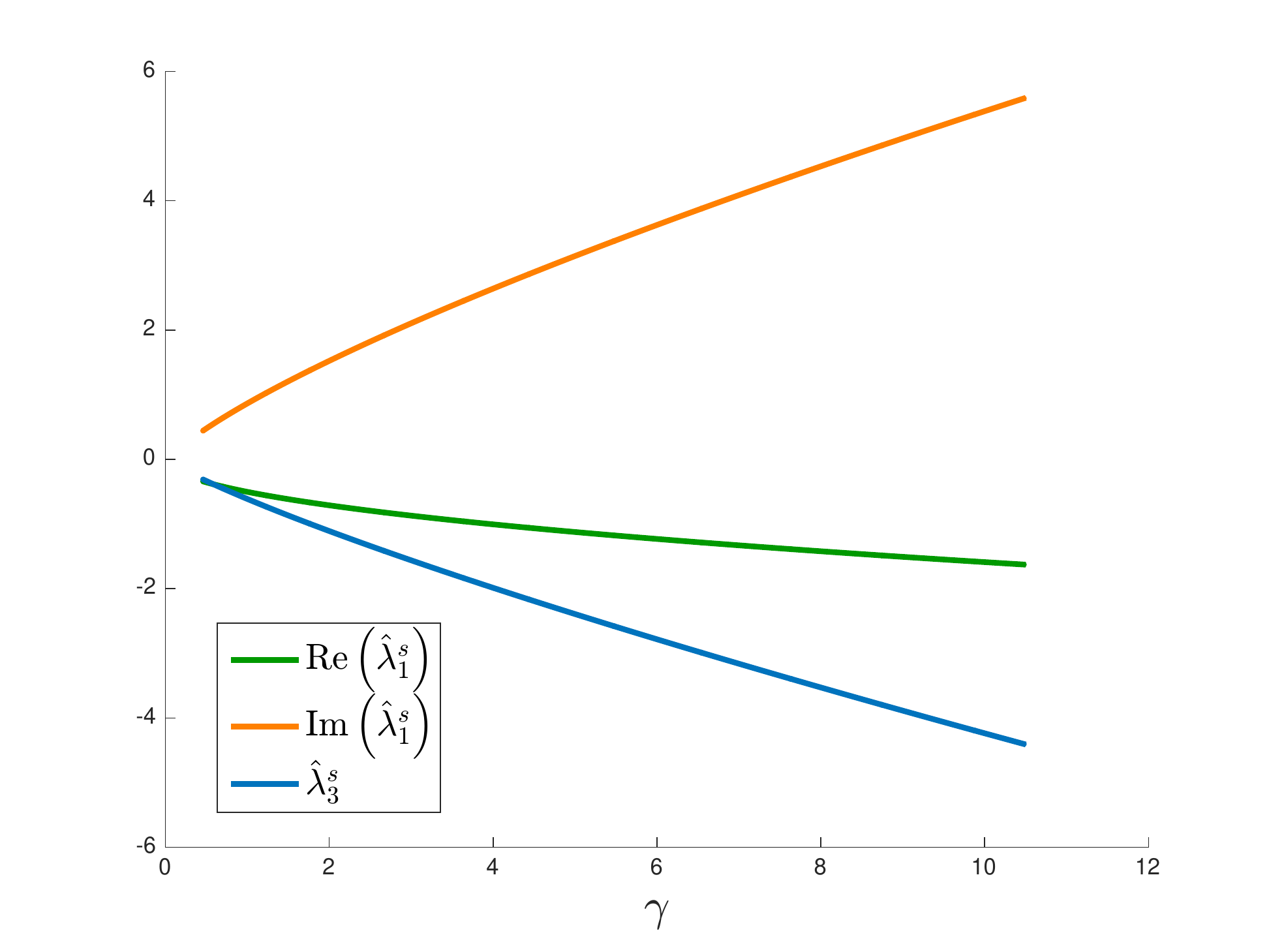}} 
	\caption[a bla]{\label{fig:FK_eig} The dependence of the stable and unstable eigenvalues on $\gamma$. 
}
\end{figure}

We remark that validation of connecting orbits for $\gamma > 10.50$ is feasible; the reason for the ``obstruction'' $\hat r> r^{\ast}$ 
was that the decay rates of the Chebyshev coefficients decreased as $\gamma$ increased, which eventually
resulted in a bound for the residual (i.e. the bound $Y^{ij}_{a}$) that was too large.
This issue can be easily resolved by using domain decomposition or increasing 
the number of Chebyshev coefficients (or by just increasing $r^{\ast}$). Similarly, validation of connecting orbits for $\gamma < 0.4557$ (but
sufficiently far away from $0$) is feasible as well; 
the bottleneck for small $\gamma$ is the validation of the local stable manifold. Indeed, as $\gamma$ decreases, the real parts of the stable eigenvalues
decrease (see Figure \ref{fig:FK_eig}) and the number of needed Taylor coefficients increases. An interesting future research project would be
to determine how close one can get to the ``singular'' case $\gamma =0$ with the current method. 

\begingroup
\let\itshape\upshape
\bibliographystyle{plain}
\bibliography{Bibliography} 

\end{document}